\documentclass[11pt]{article}

\usepackage{amssymb,latexsym,amsmath,amsthm,verbatim}
\usepackage{graphicx,epsfig,epstopdf,amssymb,color,subcaption}
\usepackage[normalem]{ulem}

\newtheorem{theorem}{Theorem}[section]

\newtheorem{remark}[theorem]{Remark}

\catcode`\@=11
\@addtoreset{equation}{section}

\def\beq{\begin{equation}}
\def\eeq{\end{equation}}
\def\beqa{\begin{eqnarray}}
\def\eeqa{\end{eqnarray}}

\def\eps{\varepsilon}

\def\RR{\mathbb{R}}
\def\CC{\mathbb{C}}

\def\LL{\mathbb{L}}

\def\C{{\cal C}}

\def\K{{\cal K}}
\def\L{{\cal L}}
\def\M{{\cal M}}

\def\O{{\cal O}}

\def\V{{\cal V}}
\def\W{{\cal W}}

\def\al{\alpha}
\def\be{\beta}
\def\de{\delta}
\def\ga{\gamma}
\def\la{\lambda}

\def\Ga{\Gamma}

\def\tc{\tilde{c}}

\def\tF{\tilde{F}}
\def\tG{\tilde{G}}
\def\tI{\tilde{I}}

\def\tHH{\tilde{{\cal H}}}

\def\tM{\tilde{M}}

\def\tq{\tilde{q}}

\def\tu{\tilde{u}}
\def\tU{\tilde{U}}
\def\tv{\tilde{v}}
\def\tV{\tilde{V}}

\def\tx{\tilde{x}}
\def\ty{\tilde{y}}

\def\tal{\tilde{\alpha}}

\def\tla{\tilde{\lambda}}

\def\ttau{\tilde{\tau}}

\def\tLL{\tilde{{\mathbb{L}}}}

\def\bF{\bar{F}}
\def\bG{\bar{G}}

\def\bu{\bar{u}}
\def\bU{\bar{U}}
\def\bv{\bar{v}}
\def\bV{\bar{V}}
\def\bW{\bar{W}}
\def\bw{\bar{w}}
\def\bal{\bar{\alpha}}

\def\vmu{\vec{\mu}}

\def\ca{c_{\ast}}

\def\Fa{F_{\ast}}
\def\Ga{G_{\ast}}

\def\tMa{\tilde{M}_{\ast}}

\def\qa{q_{\ast}}
\def\ua{u_{\ast}}
\def\va{v_{\ast}}

\def\tua{\tilde{u}_\ast}
\def\tva{\tilde{v}_\ast}
\def\tqa{\tilde{q}_\ast}

\def\hc{\hat{c}}

\def\hu{\hat{u}}
\def\hv{\hat{v}}
\def\hq{\hat{q}}
\def\htau{\hat{\tau}}

\def\stackdef{\stackrel{\rm def}{=}}

\begin{document}

\title{Criteria for the (in)stability of planar interfaces in singularly perturbed 2-component reaction-diffusion equations}

\author{Paul Carter\footnotemark[1], Arjen Doelman\footnotemark[2], Kaitlynn Lilly\footnotemark[3], Erin Obermayer\footnotemark[4], \& Shreyas Rao\footnotemark[5]}

\maketitle
\renewcommand{\thefootnote}{\fnsymbol{footnote}}
\footnotetext[1]{Department of Mathematics, University of California, Irvine, USA, pacarter@uci.edu}
\footnotetext[2]{Mathematisch Instituut, Universiteit Leiden, the Netherlands, doelman@math.leidenuniv.nl}
\footnotetext[3]{University of Maryland, Baltimore County, Baltimore, USA}
\footnotetext[4]{The College of New Jersey, Ewing, USA}
\footnotetext[5]{Brown University, Providence, USA}
\renewcommand{\thefootnote}{\arabic{footnote}}

\begin{abstract}
We consider a class of singularly perturbed 2-component reaction-diffusion equations which admit bistable traveling front solutions, manifesting as sharp, slow-fast-slow, interfaces between stable homogeneous rest states. In many example systems, such as models of desertification fronts in dryland ecosystems, such fronts can exhibit an instability by which the interface destabilizes into fingering patterns. Motivated by the appearance of such patterns, we propose two versions of a 2D stability criterion for (transversal) long wavelength perturbations along the interface of these traveling slow-fast-slow fronts. The fronts are constructed using geometric singular perturbation techniques by connecting slow orbits on two distinct normally hyperbolic slow manifolds through a heteroclinic orbit in the fast problem. The associated stability criteria are expressed in terms of the nonlinearities of the system and the slow-fast-slow structure of the fronts. We illustrate and further elaborate the general set-up by explicitly working out the existence and transversal (in)stability of traveling fronts in a number of example systems/models. We analytically establish the instability of invading bare soil/vegetation interfaces against transversal long wavelength perturbations in several dryland ecosystem models and numerically recover fingering vegetation patterns counter-invading an invading desertification front. 
\end{abstract}

\section{Introduction}
\label{s:Intro}
In this paper we consider a general system of 2-component singularly perturbed reaction-diffusion equations,
\begin{equation}
\label{e:RDE}
\left\{	
\begin{array}{rcrcl}
\tau U_t &=& \Delta U & + & F(U,V;\vmu)\\
V_t &=& \frac{1}{\varepsilon^2} \Delta V & + & G(U,V;\vmu)
\end{array}
\right.
\end{equation}
for $(x,y) \in \RR^2$, i.e. on the full unbounded 2-dimensional plane, with $U(x,y,t), V(x,y,t): \RR^2 \times \RR^+ \to \RR$, $F(U,V;\vmu)$ and $G(U,V;\vmu)$ sufficiently smooth, and $\tau > 0$, $\vmu \in \RR^m$, parameters.  Originally motivated by the fact that fronts connecting stable homogeneous states, also called coexistence states in ecology, have a direct ecological interpretation as invasion fronts -- see for instance \cite{Eig21,FOTM19} -- we investigate the existence and stability of planar interfaces between stable states in systems of reaction-diffusion equations. We focus on model (\ref{e:RDE}) and assume throughout this paper that system (\ref{e:RDE}) has (at least) two stable homogeneous background states, $(U(x,y,t),V(x,y,t)) \equiv (\bU^\pm, \bV^\pm)$. Moreover, we assume that it is singularly perturbed: we consider $0 < \eps \ll 1$ so that the diffusive spreading speed of the $V$-component in (\ref{e:RDE}) is much larger than that of the $U$-component. Since pattern formation in ecosystems is typically driven by counteracting feedback mechanisms on widely different spatial scales \cite{RvdK08}, this is a very natural assumption in the setting of ecosystem models. Although ecosystem models certainly are not restricted to 2-component models (see section \ref{s:Disc}), there is a large variety of ecological models that are covered by (\ref{e:RDE}) -- see 
\cite{BCD19,Basetal18,Doe22,Eig21,ES20,FOTM19,JDCBM20,Kla99,Mer18,SDERRS15,Sitetal14,vLetal03,ZMB15} and the references therein. Moreover, models of type (\ref{e:RDE}) occur throughout the literature in a wide variety of (non-ecological) settings -- see \cite{Doe19,DIN04,Fife88,HM94,HM97a,HM97b,NF87,NMIF90,Tan03,TN94,TK88,Ward18} and the references therein. 
\\ \\
In the case of a scalar reaction-diffusion equation,
\beq
\label{e:RDE-scalar}
W_t = \Delta W  + H(W;\vmu),
\eeq
the assumption that (\ref{e:RDE-scalar}) has two stable homogeneous rest states, $W(x,y,t) \equiv \bW^\pm$ immediately implies by phase plane techniques that there must be a monotonic front solution $W(x,y,t) = W_h(x-\ca t)$ that travels with a well-defined, unique, critical speed $\ca$ in the (longitudinal) direction -- which we identify with $x$ -- and that can be seen as a connection between the stable states $\bW^\pm$: $\lim_{\xi \to \pm \infty} W_h(\xi) =  \bW^\pm$ (with $\xi = x-\ca t)$. Its longitudinal (spectral) stability, i.e. its stability with respect to perturbations that only depend on $x$, is determined by the linearized Sturm-Liouville problem,
\[
\L_w(\xi) \bw = \bw_{\xi \xi} + \ca \bw + H_w(W_h(\xi)) \bw = \la \bw.
\]
It follows by the monotonicity of $W_h(\xi)$ that the translational eigenvalue is the critical eigenvalue, i.e. $\la_c = 0$, and thus that the front is stable \cite{KP13}. The stability of the front of (\ref{e:RDE-scalar}) as a traveling planar interface in $\RR^2$ is determined by allowing for transversal perturbations through the inclusion of Fourier modes $e^{ i \ell y}$ in the linearization Ansatz, which yields the 1-parameter family of spectral problems
\[
\bw_{\xi \xi} - \ell^2 \bw + \ca \bw + H_w(W_h(\xi)) \bw = (\L_w - \ell^2)\bw = \la(\ell) \bw, \; \; \ell \in \RR
\]
that can be reduced to the longitudinal problem by replacing $\la$ by $\la + \ell^2$. Since the parabolic curve $\la_c(\ell) = -\ell^2$ has its maximum at $\ell = 0$ it follows that the planar interface is also stable with respect to transversal perturbations and thus as traveling planar interfacial solution of (\ref{e:RDE-scalar}) on $\RR^2$. In this paper we consider the question whether similar general insights on the existence and stability of planar interfaces between stable states can be deduced for non-scalar reaction-diffusion equations.
\\
\begin{figure}[t]
\hspace{.02\textwidth}
\begin{subfigure}{.52 \textwidth}
\centering
\includegraphics[width=1\linewidth]{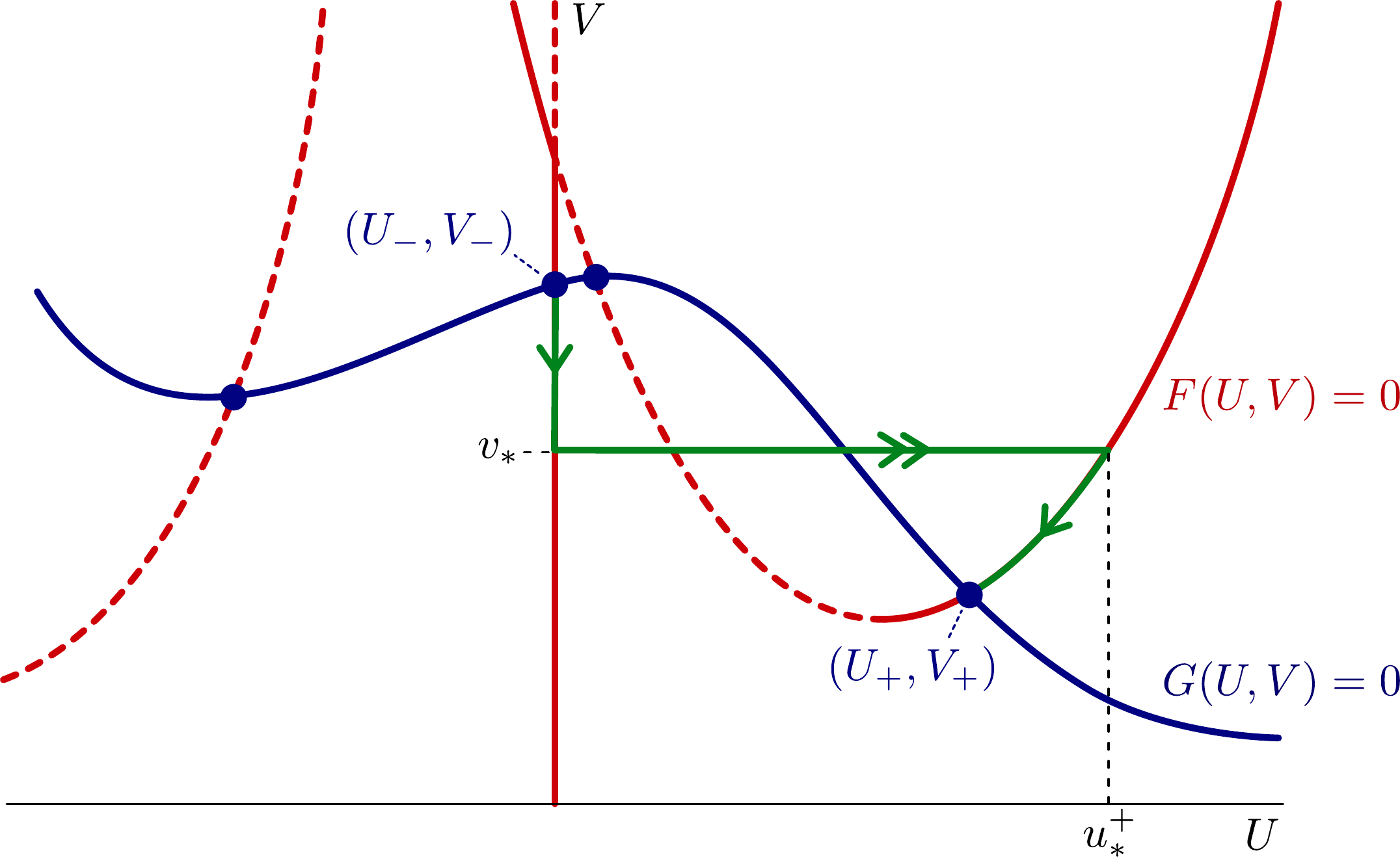}
\caption{}
\end{subfigure}
\hspace{.05\textwidth}
\begin{subfigure}{.4 \textwidth}
\centering
\includegraphics[width=1\linewidth]{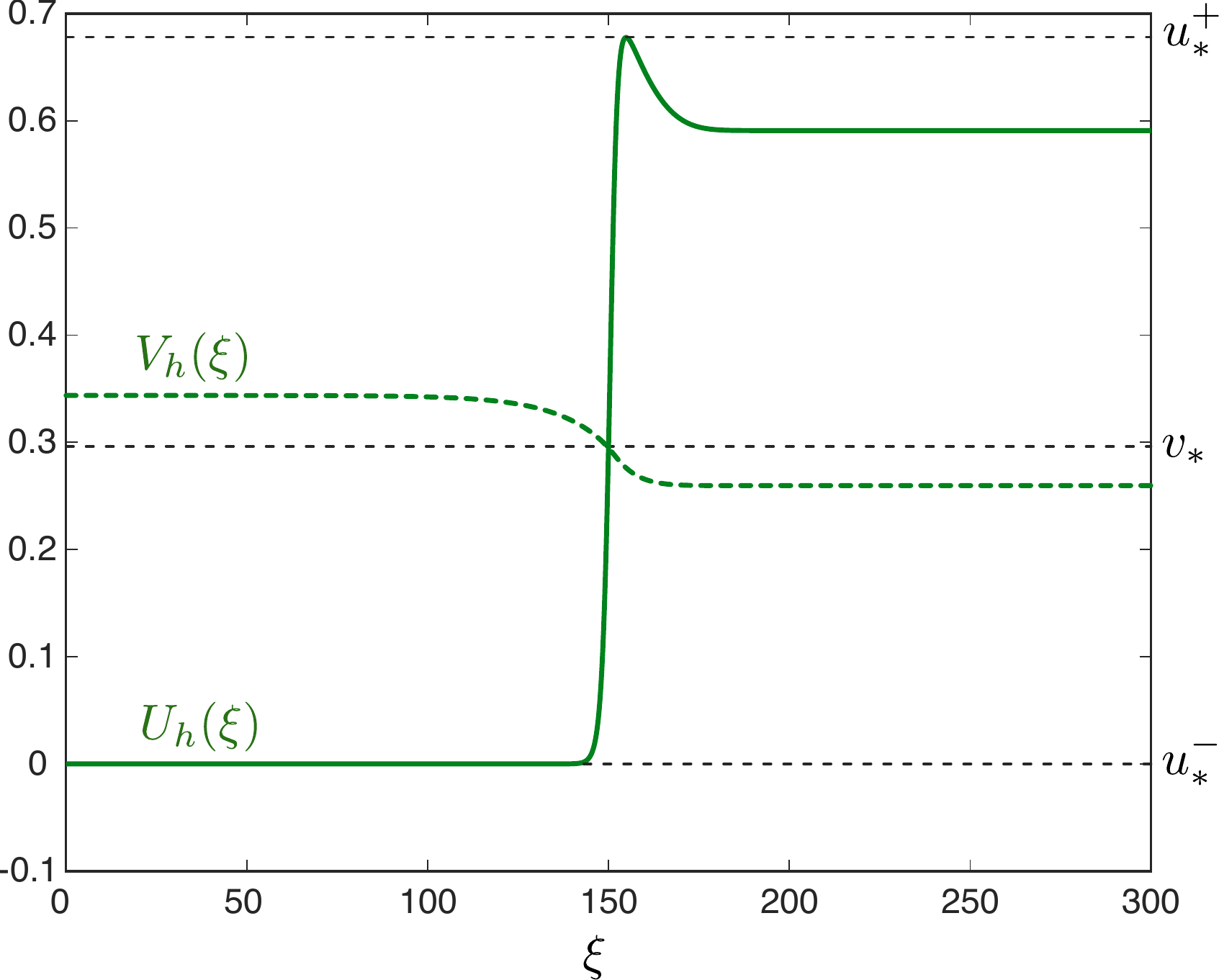}
\caption{}
\end{subfigure}
\caption{\small{A sketch of the construction and the nature of the slow-fast-slow traveling front solution in the setting of ecosystem model (\ref{e:RDE-FOTM}) of \cite{FOTM19} (a) The manifold $\{F(u,v) = 0\}$ (in red) with normally hyperbolic submanifolds $\M^\pm_0$ indicated by continuous red curves (and non-normally hyperbolic parts dashed), the stable background states $(\bU^\pm,\bV^\pm)$ as intersections of the normally hyperbolic slow manifolds and the manifold $\{G(u,v) = 0\}$ (in blue) and the projection of the singular slow-fast-slow front (in green) that makes a fast jump  between two well-defined points $(\ua^-,\va)$ -- with here $\ua^- = 0$ -- and $(\ua^+,\va)$ along the heteroclinic solution $\ua(\xi)$ of (\ref{e:fast-red}). (b) The components $U_h(\xi)$ and $V_h(\xi)$ of the front. Note that the fast jump $\ua(\xi)$ is monotonic between its endpoints $\ua^\pm$, but that $U_h(\xi)$ typically is not monotonic and that $V_h(\xi)$ only varies slowly (in $X = \eps \xi$).}}
\label{f:manifold&interface}
\end{figure}
\\
In the 2-component setting of (\ref{e:RDE}), introducing traveling coordinate $\xi = x - ct$ and $u(\xi) = U(x,y,t)$, $ v(\xi) = V(x,y,t)$, reduces (\ref{e:RDE}) to a system of coupled second order (ordinary differential) equations in $\xi$,
\begin{equation}
\label{e:ODE}
\left\{	
\begin{array}{rclcc}
u_{\xi\xi} & = & - c \tau u_\xi & - & F(u,v) \\
v_{\xi\xi} & = & -\eps^2(c v_\xi & + & G(u,v))
\end{array}
\right.
\end{equation}
which has, by taking $\eps \to 0$, the 2-parameter family
\beq
\label{e:fast-red}
u_{\xi\xi} + c \tau u_\xi +  F(u,v_0)  =  0, \; \;
v = v_0, q = q_0
\eeq
as fast reduced limit. The critical points of (\ref{e:fast-red}) are determined by the manifold $\{(u,v): F(u,v)=0\}$. Away from the (degenerate) transition points at which $F_u(u,v) = 0$, $F(u,v) = 0$ determines $J \geq 1$ branches that can be written as graphs $\{u=f^j(v)\}$ ($j = 1,2, ..., J$) with $f^j(v)$ such that $F(f^j(v),v) \equiv 0$. The background states $(\bU^\pm,\bV^\pm)$ of (\ref{e:RDE}) must correspond to points on one of these branches (since $F(\bU,\bV) = G(\bU,\bV) = 0$).  We define two specific branches $\{u=f^\pm(v)\}$ by selecting $j$'s such that $\bU^\pm = f^\pm(\bV^\pm)$. Next to our assumption that the background states $(\bU^\pm,\bV^\pm)$ are stable as solutions of (\ref{e:RDE}), we impose a second assumption that underlies our analysis: we assume that $(\bU^\pm,\bV^\pm)$ are on different branches of $\{(u,v): F(u,v)=0\}$, i.e. that $f^+(v) \not\equiv f^-(v)$. In the setting of the singularly perturbed 4-dimensional dynamical system associated to (\ref{e:ODE}) these assumptions together imply that the background states correspond to saddle points for the slow flows on two different 2-dimensional normally hyperbolic slow manifolds $\M^\pm$ that are at leading order in $\eps$ given by $\M^\pm_0 = \{ (u,p): u=f^\pm(v), p=0 \}$ -- see section \ref{ss:construction}. The traveling interfaces $(U_h(\xi), V_h(\xi))$ we study here correspond to heteroclinic orbits between the saddles associated to $(\bU^\pm,\bV^\pm)$ on these slow manifolds. The fact that $\M^- \neq \M^+$ implies that these interfaces/orbits must consist of 3 parts: a slow part along $\M^-$, followed by a fast jump from $\M^- \neq \M^+$, concluded by a final slow part along $\M^+$ (see Fig. \ref{f:manifold&interface}).
\\ \\
The details of the construction of these slow-fast-slow orbits are given in section \ref{ss:construction} for $\tau = \O(1)$ (and in section \ref{ss:ExStab-SmallTau} for $\tau = \O(\eps)$ -- see below). For the formulation and interpretation of the upcoming (in)stability criteria we need to provide somewhat more information on the (leading order) fast jump. The analysis of the (slow) flows on $\M^\pm_0$ yields an explicit critical value $\va$ (that does not depend on $c$) at which the jump must occur -- see Fig. \ref{f:manifold&interface}. Planar fast reduced system (\ref{e:fast-red}) -- with $v_0 = \va$ -- has (by construction) two saddle points $(\ua^\pm,0)$ with $\ua^\pm = f^\pm(\va)$, these are connected by a monotonic heteroclinic orbit $\ua(\xi)$ -- see again Fig. \ref{f:manifold&interface} -- for a uniquely determined value $(c \tau)_\ast$, which thus determines the speed $c = \ca = (c \tau)_\ast/\tau$: $\ua(\xi)$ is the leading order approximation of the fast jump and $\ca$ determines the (leading order) speed of the full interface $(U_h(\xi), V_h(\xi))$. Note that thus neither $\ua(\xi)$, $\ua^\pm$ nor the product $\ca \tau$ depends directly on $\tau$ (at leading order in $\eps$): these expressions are solely determined by $\va$.
\\ \\
In this paper we do not consider the analytical details of the stability of the front $(U_h(\xi), V_h(\xi))$ with respect to longitudinal perturbations, i.e. perturbations that only depend on $x \in \RR$ (or $\xi \in \RR$). Although this may be a nontrivial technical endeavour for a given system, we note that the stability can be established by the methods developed in \cite{BCD19,CdRS16,DIN04,DV15,NF87,NMIF90,Ward18} and the references therein. Thus, we assume that the front is longitudinally stable and thus that the translational eigenvalue is the most critical eigenvalue: $\la_c = 0$ (as in the scalar case). The central question considered in this paper is:
\\ \\
{\it Under what condition(s) are longitudinally stable slow-fast-slow fronts of system (\ref{e:RDE}) (un)stable as planar traveling interfaces in $\RR^2$, i.e. when are these fronts (un)stable with respect to transversal perturbations?}
\\ \\
We do not pursue this question in its full generality, although this again can also be done for a given explicit system: we focus on obtaining general (in)stability conditions against long wavelength transversal perturbations. More precise: given a longitudinally stable interface, we introduce $\ell \in \RR$ by the factor $e^{i \ell y}$ in the spectral stability Ansatz as the wavenumber of transversal perturbations so that $\la_c$ becomes a function of $\ell$ and consider $|\ell| \ll 1$. Our main general (in)stability criteria concern the sign of $\la_{2,c}$, the leading order coefficient of the expansion of $\la_c(\ell)$ for $|\ell| \ll 1$: $\la_c(\ell) = \la_{2,c} \ell^2 + \O(\ell^4)$ (since $\la_c$ is a function of $\ell^2$ by the reflection symmetry of (\ref{e:RDE}) in $y$). Our approach is based on the explicit construction of the $\la=0$-eigenfunction of the adjoint of the (matrix) operator associated to the spectral stability of the fronts $(U_h(\xi), V_h(\xi))$ -- which is based on the singular slow-fast-slow structure of the fronts (see section \ref{ss:stability}).
\\
\begin{figure}[t]
\centering
\begin{minipage}{.245\textwidth}
		\centering
		\includegraphics[width =\linewidth]{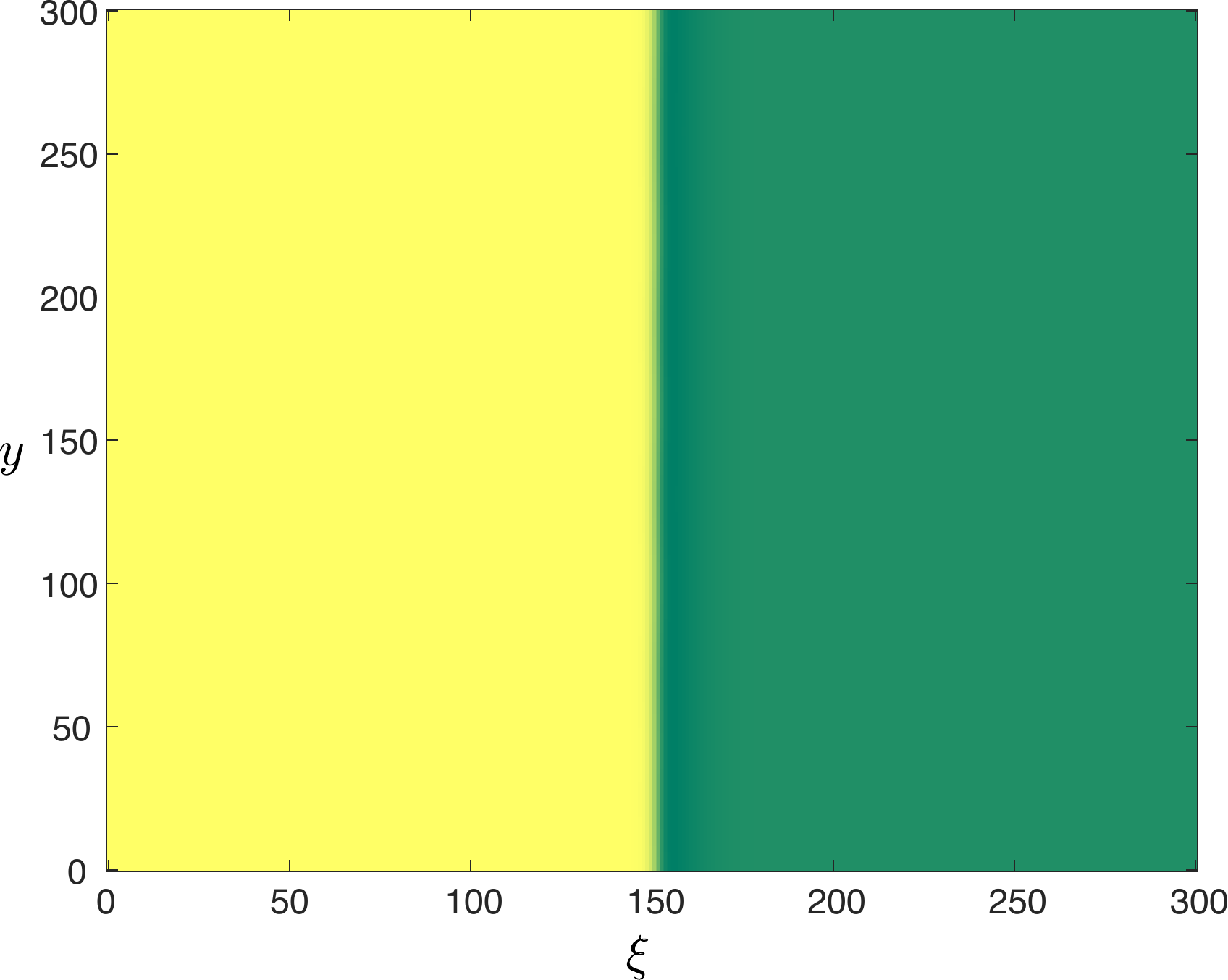}
\end{minipage}
%\hspace{.05cm}
\begin{minipage}{.245\textwidth}
		\centering
		\includegraphics[width =\linewidth]{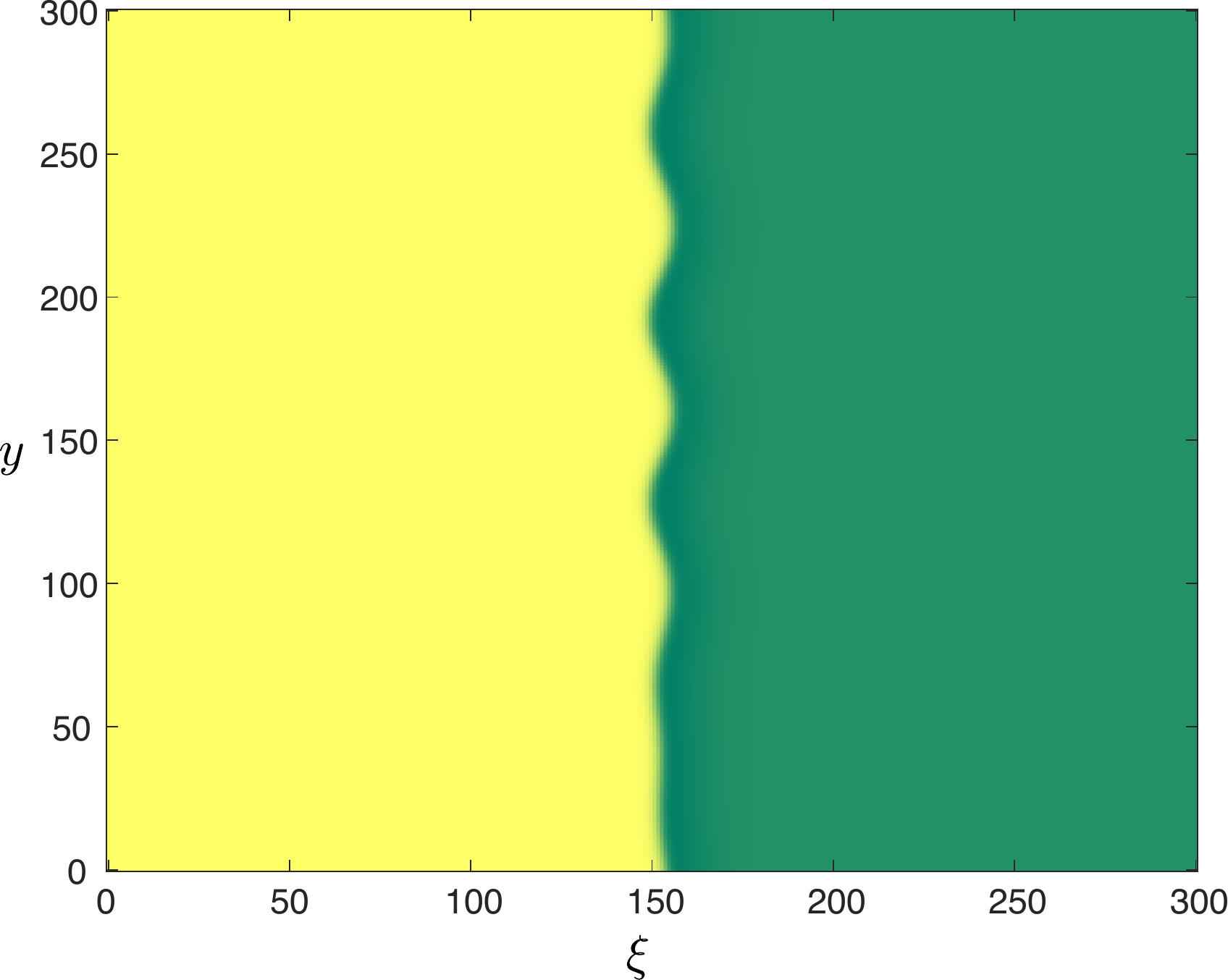}
\end{minipage}
%\hspace{.05cm}
\begin{minipage}{.245\textwidth}
		\centering
		\includegraphics[width =\linewidth]{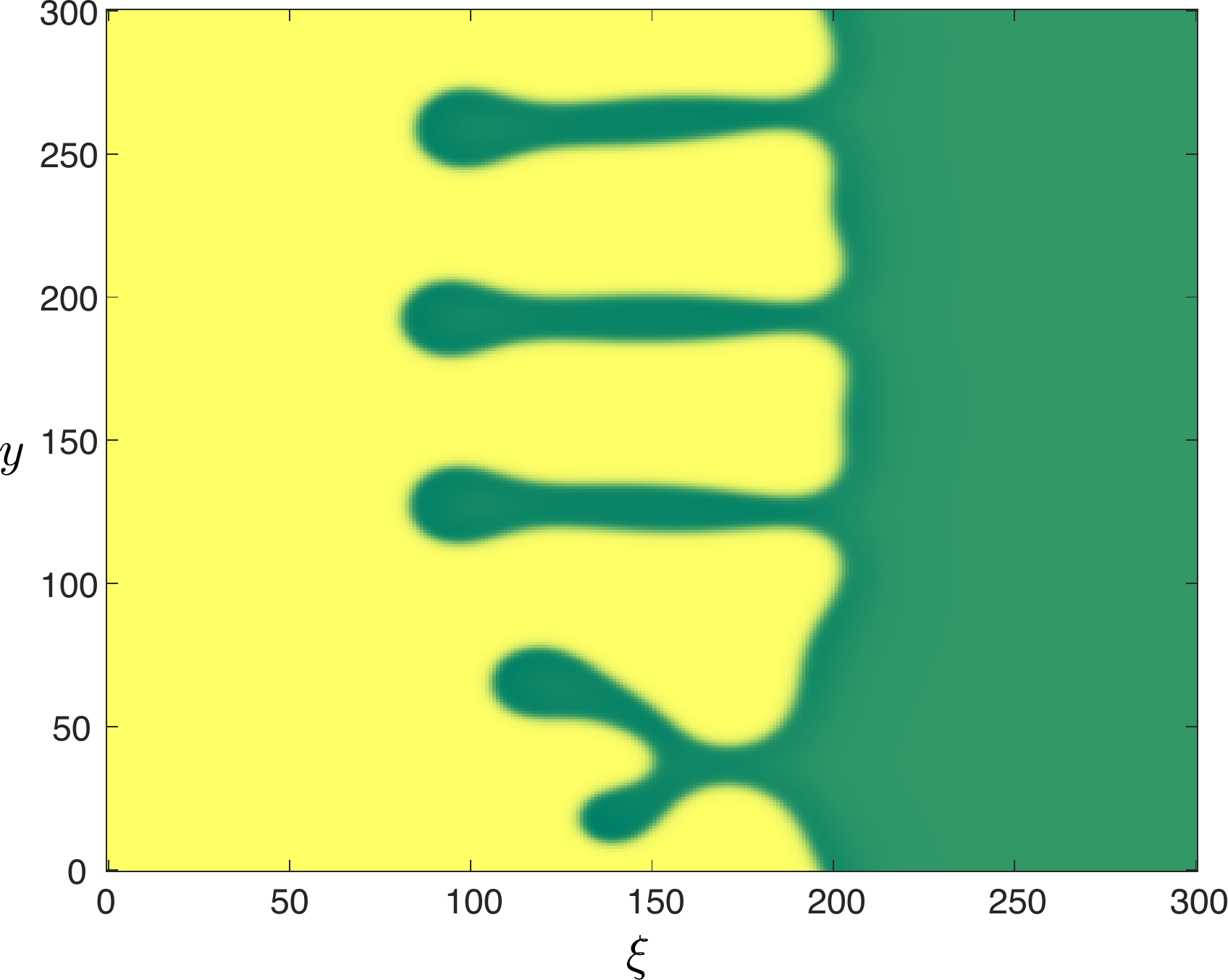}
\end{minipage}
%\hspace{.05cm}
\begin{minipage}{.245\textwidth}
		\centering
		\includegraphics[width =\linewidth]{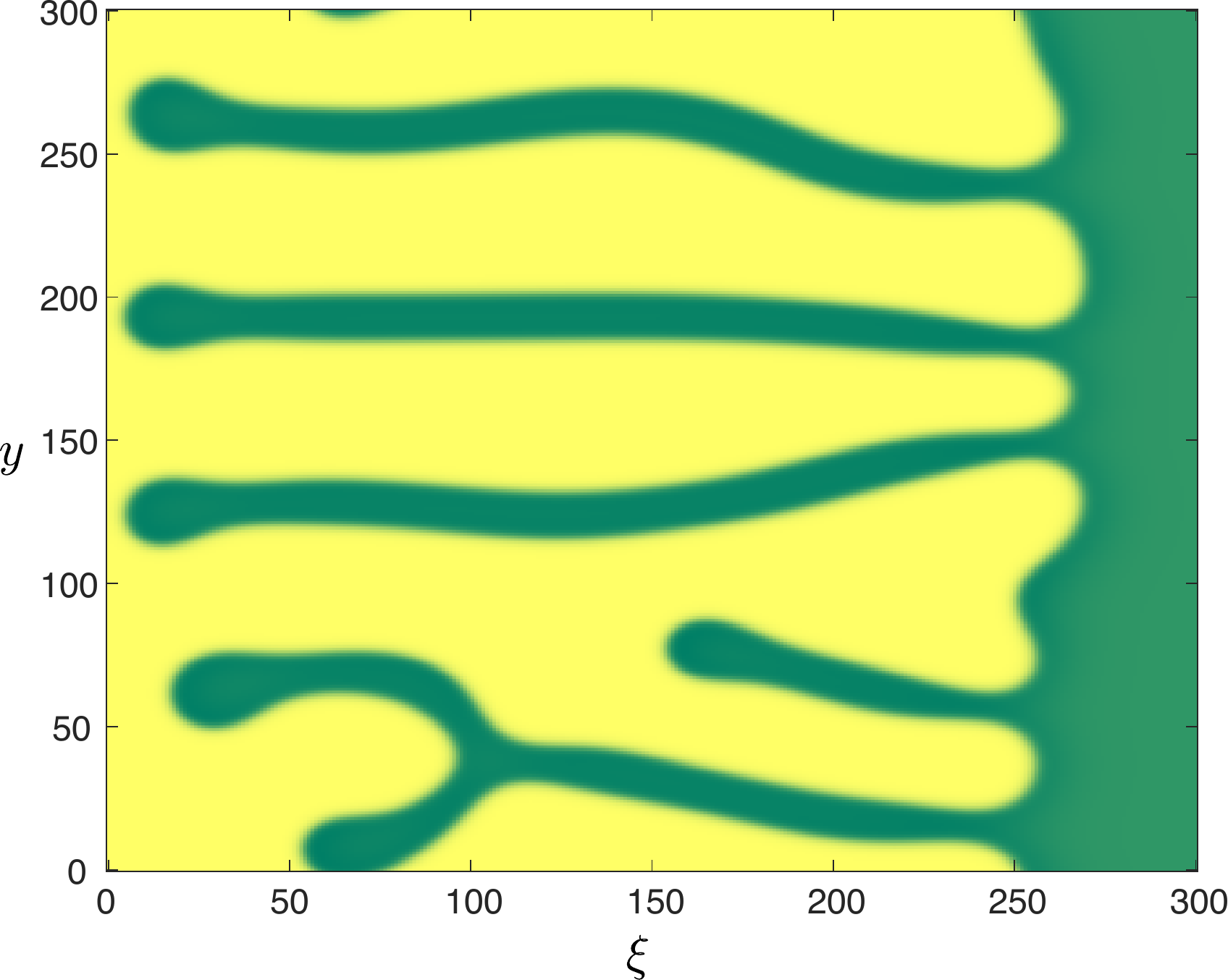}
\end{minipage}
\caption{\small{Four snapshots (from left to right, $t=500,1500,2500,3500$) of a numerical simulation of an invading planar desertification front between a bare soil state (in yellow) and a homogeneously vegetated state (in green) in ecosystem model (\ref{e:RDE-FOTM}) of \cite{FOTM19} (with $\vmu \in \RR^5$ and $\vmu =(3.5,1.1,3.2,1.0,0.4)$, $\eps = 0.04$). The simulations employed finite differences for spatial discretization and MATLAB's ode15s routine for time integration. The interface is unstable against transversal long wavelength perturbations by criterion (\ref{e:cond-longtrans}) for all ecologically realistic parameter combinations (and $\eps$ sufficiently small) -- see section \ref{sss:FOTM}. As first observed in \cite{FOTM19}: the interface is counter-invaded by fingering vegetation patterns.}}
\label{f:FingIntro}
\end{figure}
\\
In the so far considered case $\tau = \O(1)$ we define,
\beq
\label{d:FGast}
F_\ast(\vmu) = \int_\RR F_v(u_\ast(\xi),v_\ast) u_{\ast, \xi} (\xi) \, e^{c_\ast \tau \xi} \, d\xi, \;
G_\ast(\vmu) = G(u_\ast^+,v_\ast) - G(u_\ast^-,v_\ast).
\eeq
Our first (in)stability criterion reads,
\beq
\label{e:cond-longtrans}
{\rm sign}(\la_{2,c}(\vmu)) =  - \, {\rm sign}(F_\ast(\vmu)) \times {\rm sign}(G_\ast(\vmu)),
\eeq
where we note that this condition does not depend on the value of $\tau$ (as long as it is $\O(1)$ with respect to $\eps$). Based on the monotonicity of $\ua(\xi)$ and the fact that one typically knows the signs of $\va$ and $\ua(\xi)$, the simplicity of this criterion allows us to immediately draw conclusions on the instability of traveling planar interfaces in the dryland ecosystem models considered in \cite{BCD19,Eig21,FOTM19,JDCBM20}, and thus on the counter-invasion of fingering vegetation patterns against an invading desertification front \cite{JDCBM20} -- see Fig. \ref{f:FingIntro} and section \ref{ss:fingering}. In fact, if (\ref{e:RDE}) is of activator-inhibitor type \cite{Mur02} with (for instance) $V(x,y,t)$ as activator  -- i.e. $F_v(u,v) > 0$ -- and $U(x,y,t)$ as inhibitor -- $G_u(u,v) < 0$ -- it follows from (\ref{d:FGast}) and (\ref{e:cond-longtrans}) that the (in)stability of (invading) planar interfaces against transversal long wavelength perturbations typically is decided by the signs of $\bU^+-\bU^-$ and $\bV^+-\bV^-$.
\\ \\
Clearly, $\la_{2,c} = -1 < 0$ in the scalar case of (\ref{e:RDE-scalar}): the longitudinal perturbations are potentially the most `dangerous' perturbations in (\ref{e:RDE-scalar}). The `superficial intuition' that interfaces in 2-component model (\ref{e:RDE}) typically behave like the scalar case thus is incorrect: it follows from (\ref{e:cond-longtrans}) -- and likewise from (\ref{e:cond-longtrans-t-Intro}) below -- that $\la_{2,c} > 0$ occurs as natural as $\la_{2,c} < 0$: a longitudinally stable front is unstable against transversal long wavelength perturbations under a simple, explicitly verifiable, condition on the vector field $(F,G)$ of (\ref{e:RDE}) -- see section \ref{s:Models}. In this paper we do not go into the nonlinear nature of this instability. Typically, it gives rise to `fingering patterns' \cite{FOTM19,HM94,TN94}. This type of growth pattern originates from fluid dynamics and is associated to the sideband mechanism \cite{Cou00}, i.e. the instability with respect to perturbations with wave numbers $|k| \ll 1$, while $k=0$ itself is (marginally) stable \cite{CH93,Doe19} -- as is the case here (with $\ell$ in the role of $k$). Although understanding the onset of counter-invading fingers -- \cite{FOTM19} and Fig. \ref{f:FingIntro} -- formed the original motivation for the research presented in this paper, we note that $\la_{2,c} > 0$ does not necessarily induce the formation of fingers: see the simulation of Fig. \ref{f:FingCusp-BCDE} in which a transversally unstable invading vegetation front in an ecosystem model based on \cite{BCD19,Eig21} exhibits evolving bounded cusps reminiscent of Kuramoto-Sivashinsky dynamics associated to a sideband instability \cite{DSSS09,HN86} (however, in this model invading desertification fronts also exhibit counter-invading vegetation fingers as in \cite{FOTM19} and Fig. \ref{f:FingIntro} -- see again Fig. \ref{f:FingCusp-BCDE}).
\\ \\
By rescaling space -- $(\tx,\ty) = (\eps x, \eps y)$ -- and defining $D = \eps^2$, we may interpret the parameters $(\tau, D)$ as measures for the relative magnitudes of the evolution rates of $U$ and $V$ in time ($\tau$) and space ($D$). The  choices sofar made -- $\tau = O(1)$ and $D \ll 1$ -- are based on the ecosystem models that originally motivated our investigations. However, in the non-ecological literature on models of the type (\ref{e:RDE}) often the scaling $\tau = O(\sqrt{D})$ is considered -- see for instance \cite{Fife88,HM94,HM97a,HM97b,NMIF90,Tan03,TN94} and the references therein. In fact, in \cite{HM94,HM97a,HM97b} the parameters $(\nu, \eta)$ are introduced, with $\nu = \eps^2$, $\eta = \tau/\eps$ in the notation of (\ref{e:RDE}). In the notation of \cite{HM94,HM97a,HM97b} we have so far thus considered $0 < \nu \ll 1$, $\eta \gg 1$ , while $\eta = \O(1)$ (with respect to $\eps$ or $\nu$) is the critical -- richer -- scaling. In other words, $\tau = O(\sqrt{D})$ is a `significant degeneration' in the terminology of \cite{Eck79}. Therefore, we introduce $\ttau$ by setting $\tau = \eps \ttau$ and set up an existence/stability analysis similar to the $\tau = \O(1)$ case considered so far.
\\ \\
For $\tau = \O(\eps)$, the existence of the planar fronts can be approached by the same geometrical set-up as for $\tau = \O(1)$. However, the `algebra' becomes more involved: for $\tau = \O(1)$ one can first determine $\va$ by only considering the slow flows on the manifolds $\M^\pm_0$; $\ca$ is determined by the fast field, using $\va$ as (given) input. As was already noted, the fast heteroclinic jump $\ua(\xi)$ in (\ref{e:fast-red}) occurs for a given value $(\tau c)_\ast = \tau \ca$; thus as $\tau$ decreases, $\ca$ increases: $c$ must be scaled as $\tc/\eps$ for $\tau = \eps \ttau$. However, in that case $\tc$ will appear as $\O(1)$ term in the slow reduced flows: as for $\tau = \O(1)$, $\tva$ can still be determined, however, it becomes a function of $\tc$. The same (geometric) procedure by which $\ca$ is determined for given $\va$, now determines $\tc$ as function of $\va(\tc)$: $\tva$ and $\tc$ are determined simultaneously by algebraic relations (see section \ref{ss:ExStab-SmallTau}). Taking $\tau = \O(\eps)$ has a similar impact on the stability analysis: the essence of the method is not affected, but working out the details becomes more involved. As a consequence, the criterion for the sign of $\tla_{2,c}$, i.e. for the (in)stability of the associated planar interface against transversal long wavelength perturbations, changes from (\ref{e:cond-longtrans}) into,
\beq
\label{e:cond-longtrans-t-Intro}
{\rm sign}(\tla_{2,c}(\vmu; \ttau)) =  - \, {\rm sign}(\tF_\ast(\vmu;\ttau)) \times {\rm sign}(\tM_\ast(\vmu; \ttau))
\eeq
where $\tF_\ast$ is completely similar to $\Fa$ (cf. (\ref{d:FGast}) and (\ref{d:tFGast})), but where $\tM_\ast(\vmu; \ttau)$ is a somewhat more involved Melnikov-type expression (see (\ref{d:tMa})).
\\ \\
There is a significant difference between criteria (\ref{e:cond-longtrans}) and (\ref{e:cond-longtrans-t-Intro}): unlike (\ref{e:cond-longtrans}), (\ref{e:cond-longtrans-t-Intro}) depends explicitly on $\ttau$. To better understand the impact of $\ttau$, we deduce a geometrical interpretation of the implication of a change in sign of Melnikov expression $\tM_\ast(\vmu; \ttau)$. To do so, we restrict ourselves for simplicity to stationary fronts -- as is common in studies of (\ref{e:RDE}) (see for instance \cite{HM94,NF87,TN94}). We show in section \ref{ss:BifTrav-SmallTau} that a sign change of $\tM_\ast(\vmu; \ttau)$ initiates a bifurcation into traveling waves: if $\ttau_\ast$ is such that $\tM_\ast(\vmu; \ttau_\ast) = 0$, then (generically) a pair of counter-propagating traveling fronts bifurcates off from the stationary front as $\ttau$ passes through $\ttau_\ast$. Note that this implies that an expression obtained in the context of the stability analysis of planar interfaces against transversal perturbations has an interpretation in terms of the existence analysis of one-dimensional fronts -- see section \ref{s:Disc}.
\\ \\
In sections \ref{ss:FHN} and \ref{ss:ExampleSystems-t}, we consider versions of the example systems formulated and studied in \cite{HM94,HM97a,HM97b,Tan03,TN94}. We show how criteria (\ref{e:cond-longtrans}) and (\ref{e:cond-longtrans-t-Intro}) immediately provide analytical insight in the instability of planar interfaces for $(x,y) \in \RR^2$. We test the outcome of our analysis against numerical evaluations of the spectral curve $\la_c(\ell)$ in the context of a FitzHugh-Nagumo type model based on \cite{HM94,HM97a,HM97b} in section \ref{ss:FHN} and compare our results in section \ref{sss:Cylindrical} to those of \cite{Tan03,TN94} in which cylindrical spatial domains are considered. Moreover, in section \ref{sss:InStabBifTrav} we apply (\ref{e:cond-longtrans-t-Intro}) to the counter-propagating fronts initiated by the bifurcation into traveling waves in a special case of the model in \ref{ss:FHN} to show that these interfaces inherit the instability against long wavelength transversal perturbations from the original stationary interface, while $\tM_\ast(\vmu; \ttau)$ has changed sign for the stationary front. We conclude the paper with a section in which we discuss both the ecological and mathematical implications of our findings.
\\ \\
Finally, we briefly remark on the writing style of the paper. Based on the present state of the art literature on geometric singularly perturbation theory (for the existence problem) and the analysis of singularly perturbed spectral operators (for the stability), all our findings can be established rigorously. However, for transparency of presentation we have chosen not to formulate our main results as theorems, and thus not to provide all details of the proofs of these theorems. Instead, we have chosen to present our work in the present hands-on setting, with a central role for section \ref{s:Models} in which we explicitly show how our methods work in the context of a number of diverse and explicit problems. In a way, our paper can be seen as a `user guide': we have written the paper such that it can directly be used for the analysis of (the stability of) planar interfaces in explicit systems of the type (\ref{e:RDE}).

\begin{remark}
\label{r:FG-noeps}
\rm
We only consider $F(U,V)$, $G(U,V)$ that do not explicitly depend on $\eps$. This is a natural assumption from the ecological point of view \cite{Doe22}, however, there is a large body of literature on singularly perturbed reaction-diffusion equations in which $F(U,V)$ and $G(U,V)$ do depend on $\eps$ -- such as the classical Gray-Scott/Gierer-Meinhardt-type models (cf. \cite{Doe22,Ward18} and the references therein).
\end{remark}

\section{The case $\tau = \O(1)$}
\label{s:Tau=O1}

In this section we first set up the (geometrical) construction of the travelling slow-fast-slow front solutions of (\ref{e:RDE})/(\ref{e:ODE}) and next consider its  stability as 2-dimensional planar interface against against transversal long wavelength perturbations. For clarity, we list the assumptions that will be imposed on the existence problem in this section (and the next):
\\ \\
$\bullet$ {\bf (E-I)} There are two homogeneous background states $(\bU^-,\bV^-)$ and $(\bU^+,\bV^+)$ that are stable as solutions of (\ref{e:RDE}).
\\
$\bullet$ {\bf (E-II)} The background states $(\bU^-,\bV^-)$ and $(\bU^+,\bV^+)$ correspond to critical points on two different slow manifolds $\M^-$ and $\M^+$ associated to singularly perturbed system (\ref{e:ODE})/(\ref{e:DS}). \\
$\bullet$ {\bf (E-III)} The projections on the $v$-axis of the manifolds $\M^-$ and $\M^+$ have non-empty overlap.
\\ \\
We refer to the additional assumption {\bf (E-IV)} in upcoming section \ref{ss:construction}. We note that similar front patterns as the ones (to be) constructed here -- and in section \ref{ss:ExStab-SmallTau} for $\tau = \O(\eps)$ -- have been established in \cite{Fife88,NF87,Tan03,TN94,TK88}, also for general systems of the type (\ref{e:RDE}), but typically under stronger/more explicit conditions on the reaction terms $(F(U,V), G(U,V))$ of (\ref{e:RDE}) than in the present setting -- see also section \ref{sss:Cylindrical} for some more details. However, the precise nature of the conditions imposed in the literature are for a large part for practical reasons, the methods used (and developed) can typically also be applied to more general systems (see again section \ref{sss:Cylindrical}). We also refer to the concluding paragraph of (upcoming) section \ref{ss:construction} and to Remark \ref{r:normallyhyperbolic} for an extended discussion on the assumptions under which slow-fast-slow fronts can be constructed.  
\\ \\
There is one crucial additional assumption for the subsequent stability analysis:
\\ \\
$\bullet$ {\bf (S-I)} The constructed slow-fast-slow heteroclinic traveling front is stable against longitudinal perturbations, i.e. perturbations that only depend on $x$.
\\ \\
We refer to section \ref{ss:stability} for a more precise statement and a brief discussion.

\subsection{The construction of slow-fast-slow interfaces}
\label{ss:construction}
There are 2 equivalent ways to write (\ref{e:ODE}) as a 4-dimensional spatial dynamical system,
\begin{equation}
\label{e:DS}
({\rm fast}) \; \;
\left\{
\begin{array}{rcl}
u_\xi & = & p
\\
p_\xi & = & - F(u,v) - c \tau p
\\
v_\xi & = & \eps q
\\
q_\xi & = & \eps [- G(u,v)- \eps c q]
\end{array}
\right.
\; \; ({\rm slow}) \; \;
\left\{
\begin{array}{rcl}
\eps u_X & = & p
\\
\eps p_X & = & - F(u,v) - c \tau p
\\
v_X & = & q
\\
q_X & = & - G(u,v)- \eps c q
\end{array}
\right.
\end{equation}
with $X = \eps \xi$ the slow (spatial) coordinate. The fronts/interfaces of interest correspond to heteroclinic orbits $\ga_h(\xi)$ in (\ref{e:DS}) that connect its critical points $(\bU^-,0,\bV^-,0)$ and $(\bU^+,0,\bV^+,0)$ for a critical value $c_\ast$ of the speed $c$. The orbits $\ga_h(\xi)$ can be constructed explicitly by the methods of geometric singular perturbation theory. To do so, we reintroduce the two slow manifolds $\M^\pm_0$,
\beq
\label{d:Mpm0}
\M^\pm_0 = \{(u,p,v,q) \in \mathbb{R}^4: u=f^\pm(v), p=0, v^\pm_- < v < v^\pm_+ \}
\eeq
with $F(f^\pm(v),v) \equiv 0$ and `endpoints' $v^\pm_\pm$ at which $F_u(f^\pm(v^\pm_\pm),v^\pm_\pm) = 0$ (where we note that $\M^\pm_0$ may be unbounded, i.e. that $v^\pm_+$ may be $\infty$ and/or $v^\pm_-$ may be $-\infty$). As stated in {\bf (E-III)}, it is necessary for the construction of $\ga_h(\xi)$ to assume that the projections of $\M^\pm_0$ on the $v$-axis overlap, i.e. that $(v^-_-, v^-_+) \cap (v^+_-, v^+_+) \neq \emptyset$. Naturally, $(\bU^\pm,0,\bV^\pm,0) \in \M^\pm_0$ are critical points of the full system (\ref{e:DS}) and $(\bV^\pm,0)$ of the reduced slow systems,
\beq
\label{e:slow-red}
v_{XX} + G(f^\pm(v),v)  =  0,
\eeq
where we note that these (planar) systems are integrable (Hamiltonian) under the implicit assumption that $\eps c \to 0$ as $\eps \to 0$ (cf. (\ref{e:slow-red-t}) for $\tau = \O(\eps)$ where $\eps c = \tc = \O(1)$). The assumption that both $(\bU^\pm,\bV^\pm)$ are stable as homogeneous backgrounds states of (\ref{e:RDE}) has a direct impact on the nature of the critical points $(\bU^\pm,0,\bV^\pm,0)$ of (\ref{e:DS}), of the critical points $(\bV^\pm,0)$ of (\ref{e:slow-red}) and of the slow manifolds $\M^\pm_0$. To see this, we note that the states $(\bU^\pm,\bV^\pm)$ are stable if the eigenvalues $\la^\pm_{1,2}(k,\ell)$ of
\beq
\label{e:stabbUVpm}
\tau \la^2 - \left[(\bF_u^\pm + \tau \bG_v^\pm) - \frac{(\tau + \eps^2)(k^2 + \ell^2)}{\eps^2} \right] \la + \left[(\bF_u^\pm \bG_v^\pm - \bF_v^\pm \bG_u^\pm) - \frac{(\bF_u^\pm + \eps^2 \bG_v^\pm)(k^2 + \ell^2)}{\eps^2} + \frac{(k^2 + \ell^2)^2}{\eps^2}\right] = 0
\eeq
satisfy Re$\left(\la^\pm_{1,2}(k,\ell)\right) < 0$ for all $k, \ell \in \RR$ (where $\bF_u^\pm = \frac{\partial F}{\partial U}(\bU^\pm,\bV^\pm)$, etc., and $k, \ell$ are the (Fourier) wavenumbers in the $x, y$-directions) \cite{Doe19}.
For $0 < \eps \ll 1$, this is the case if (and only if),
\beq
\label{e:cond-stabbUVpm}
\bF_u^\pm + \tau \bG_v^\pm < 0, \; \bF_u^\pm \bG_v^\pm - \bF_v^\pm \bG_u^\pm > 0, \; \bF_u^\pm < 0
\eeq
\cite{Doe22}, where we note that the last 2 conditions come from the observation that the $\la$-independent term in (\ref{e:stabbUVpm}) must be positive for all $k, \ell \in \RR$ (consider $k = \ell = 0$ and $k^2 + \ell^2 = \O(\eps)$). Together, $\bF_u^\pm \bG_v^\pm - \bF_v^\pm \bG_u^\pm > 0$ and $\bF_u^\pm < 0$ imply that $(\bV^\pm,0)$ is a saddle of the slow flow (\ref{e:slow-red}) -- where we note that $(f^\pm)'(\bV^\pm) = - \bF_v^\pm/\bF_u^\pm$. Moreover, the critical point $(\bU^\pm,0)$ of the fast reduced system (\ref{e:fast-red}) at the level $v_0 = \bV^\pm$ also must be a saddle ($\bF_u^\pm < 0$). Since $0 < \eps \ll 1$ this implies that the eigenvalues $\la^\pm_j$, $j = 1,...,4$ associated to the critical points $(\bU^\pm,0,\bV^\pm,0)$ of (\ref{e:DS}) satisfy $\la^\pm_1 < \la^\pm_2 < 0 < \la^\pm_3 < \la^\pm_4$ (with $\la^\pm_{1,4} = \O(1)$ determined by (\ref{e:fast-red}) and $\la^\pm_{2,3} = \O(\eps)$ by (\ref{e:slow-red})). Moreover, the various branches $u = f^j(v)$ of the manifold $F(u,v) = 0$ are separated from each other by transition points at which $F_u(u,v) = 0$ (Fig. \ref{f:manifold&interface}(a)), i.e. $F_u(f^\pm(v^\pm_\pm), v^\pm_\pm)=0$ in the case of the curves $f^\pm(v)$. As a consequence, the sign of $F_u(f^\pm(v),v)$ does not change along $\M^\pm_0$: the stability of $(\bU^-,0,\bV^-,0)$ implies that $F_u(f^\pm(v),v) < 0$ on $\M^\pm_0$. This means that $\M^\pm_0$ are both {\it normally hyperbolic}, so that they persist as invariant manifolds $\M^\pm$ for the flow of (\ref{e:DS}) (for $\eps$ sufficiently small \cite{Fen79,Jon95}) -- see also Remark \ref{r:normallyhyperbolic}.
\\
\begin{figure}[t]
\centering
\includegraphics[width=\linewidth]{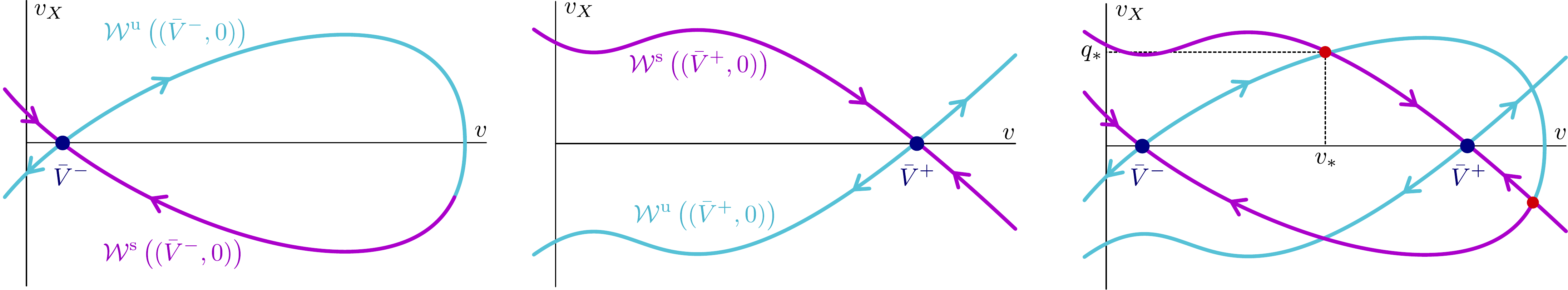}
\caption{\small{Sketches of possible slow reduced flows (\ref{e:slow-red}) on $\M^-_0$ -- with saddle $(\bV^-,0)$ and its unstable manifold $\W^{\rm u}((\bV^-,0))$ (left) -- and $\M^+_0$ -- with saddle $(\bV^+,0)$ and stable manifold $\W^{\rm s}((\bV^+,0))$ (middle). On the right, their combined projection with `jumping point' $(\va,\qa) \in \W^{\rm u}((\bV^-,0)) \cap \W^{\rm s}((\bV^+,0))$ explicitly indicated. Note that $\W^{\rm u}((\bV^-,0)) \cap \W^{\rm s}((\bV^+,0))$ is not uniquely determined in this case: it contains two points (indicated in red).}}
\label{f:slowflows}
\end{figure}
\\
The `skeleton structure' of the heteroclinic orbit $\ga_h(\xi)$ can now be built by combining the reduced slow flows (\ref{e:slow-red}) with an intermediate fast jump governed by the fast reduced system (\ref{e:fast-red}). During its first slow part, $\ga_h(\xi)$ will follow, i.e. will be exponentially close to, the unstable manifold $\W^{\rm u}((\bV^-,0))$ of the saddle $(\bV^-,0)$ on $\M^-_0$ (Fig. \ref{f:slowflows}). Likewise, it will follow the stable manifold $\W^{\rm s}((\bV^+,0))$ of the saddle $(\bV^+,0)$ on $\M^+_0$ during its second slow part. Neither the $v$- nor the $v_X = q$-component of $\ga_h(\xi)$ will change (at leading order in $\eps$) during the fast jump  (\ref{e:fast-red}). Thus, the $\ga_h$-skeleton can be constructed if the unstable manifold $\W^{\rm u}((\bV^-,0))$ on $\M^-_0$ intersects the stable manifold $\W^{\rm s}((\bV^+,0))$ on $\M^+_0$ in the combined projections of the reduced slow flows on $\M^\pm_0$ (Fig. \ref{f:slowflows}): (the skeleton of) $\ga_h(\xi)$ `takes off' from $\M^-_0$ at an intersection $(v_\ast,q_\ast) \in \W^{\rm u}((\bV^-,0)) \cap \W^{\rm s}((\bV^+,0))$, follows the fast field, and `touches down' on $\M^+_0$ again at $(v_\ast,q_\ast)$. Or more precisely, $\ga_h(\xi)$ follows the fast flow from the point $(u^-_\ast,0,v_\ast,q_\ast) \in \M^-_0$ to $(u^+_\ast,0,v_\ast,q_\ast) \in \M^+_0$, with $u^\pm_\ast = f^\pm(v_\ast)$, so that $(u^\pm_\ast,0)$ are critical points of (\ref{e:fast-red}) with $v_0 = v_\ast$ -- see Fig. \ref{f:manifold&interface}(a) and upcoming Fig. \ref{f:BCDE}(a) with a 3-dimensional sketch in the context of ecosystem model (\ref{e:RDE-BCDE}).
\\
%\begin{figure}[t]
%\centering
%\includegraphics[width=9cm]{"Fig3-sketch"}
%\caption{\small{The fast jump.}}
%\label{f:fastflow}
%\end{figure}
\\
The speed $c_\ast$ of the traveling interface is determined by the passage through the fast field (at $v = v_\ast$, $q = q_\ast$). In general, i.e. for general $c$, the unstable manifold of the saddle $(u^-_\ast,0)$ of (\ref{e:fast-red}) will not be connected to the stable manifold of $(u^-_\ast,0)$ -- where we recall that the assumed stability of background states $(\bU^\pm,\bV^\pm)$ implies that the critical points of (\ref{e:fast-red}) are saddles. The situation is now completely similar to that of constructing heteroclinic traveling fronts for scalar equation (\ref{e:RDE-scalar}) -- in fact, by (re)defining $F(u,\va)$ as $H(u)$, (\ref{e:fast-red}) is identical to the (spatial) traveling wave planar dynamical system associated to (\ref{e:RDE-scalar}). Thus, for a given $\va$ we can (by phase plane techniques) determine a unique value $\ca$ of $c$ for which there is a (planar) heteroclinic connection $(u_\ast(\xi),p_\ast(\xi))$ between the saddles $(u^\pm_\ast,0)$. By construction, $p_\ast(\xi) = u_{\ast,\xi}(\xi)$ is either strictly positive or strictly negative:  $(u_\ast, p_\ast)$ forms the (monotonic) fast component of the above skeleton structure (Fig. \ref{f:manifold&interface}). (It should be noted that although the fast jump $\ua(\xi)$ is monotonic, the $u$-component $u_h(\xi)$ of the full slow-fast-slow connecting front typically is not monotonic -- see Fig. \ref{f:manifold&interface}.)
\\ \\
In the setting of a given model it may be a technical undertaking to explicitly set up the construction of the skeleton of the heteroclinic orbit $\ga_h(\xi) = (u_h(\xi),p_h(\xi),v_h(\xi),q_h(\xi))$ -- see section \ref{s:Models}. Here, in the general setting, we assume that this has been done. By definition, the piecewise approximations of $v_h(\xi)$ are given by $v^-_u(X)$ near $\M^-_0$ and $v^+_s(X)$ near $\M^+_0$, while $v_h(\xi) = v_\ast$, i.e. constant, over the fast jump; by defining $u^-_u(X) = f^-(v^-_u(X))$ and $u^+_s(X) = f^+(v^+_s(X))$ we obtain similar leading order approximations of $u_h(\xi)$ in the slow fields. Naturally, $u_h(\xi)$ is approximated by $u_\ast(\xi)$ during the fast jump -- see Figs. \ref{f:manifold&interface} and \ref{f:slowflows}. Note that $p_h(\xi) = u_{h, \xi}(\xi)$ and $q_h(X) = v_{h, X}(X)$ by (\ref{e:DS}). From the geometrical point of view, the construction of the full slow-fast-slow front is now based on a transversal intersection of the (3-dimensional) stable and unstable manifolds $\W^{\rm u}(\M^-_0)$ and $\W^{\rm s}(\M^+_0)$ (in the $\eps \to 0$ limit) and it thus follows by the methods of geometric singular perturbation theory that asymptotically close to this skeleton there indeed is an orbit -- the orbit $\ga_h(\xi)$ -- that is a heteroclinic connection between the critical points $(\bU^-,0,\bV^-,0)$ and $(\bU^+,0,\bV^+,0)$. This orbit corresponds by construction to the (planar) traveling interface $(U(x,y,t),V(x,y,t)) = (U_h(x-c_\ast t), V_h(x-c_\ast t)) = (u_h(\xi),v_h(\xi))$ -- see Fig. \ref{f:manifold&interface}.
\\ \\
Wrapping up the existence analysis, we return to the conditions under which the traveling interface $(U_h(\xi), V_h(\xi))$ can be constructed. Additional to {\bf (E-I)} - {\bf (E-III)} stated at the beginning of this section, that did not consider the slow flows on $\M^\pm_0$, we now add the fourth condition:
\\ \\
$\bullet$ {\bf (E-IV)} The unstable manifold $\W^{\rm u}((\bV^-,0))$ of the slow flow (\ref{e:slow-red}) on $\M^-_0$ intersects the stable manifold $\W^{\rm s}((\bV^+,0))$ of the slow flow (\ref{e:slow-red}) on $\M^+_0$ in their combined projection. 
\\ \\
Note that this certainly is not `automatically' the case, on the other hand, $\W^{\rm u}((\bV^-,0)) \cap \W^{\rm s}((\bV^+,0))$ may consist of several points $(v_{\ast, i}, q_{\ast, i})$ -- see Fig. \ref{f:slowflows} -- so that there may exist several distinct connecting fronts for a given choice of parameters $\vmu$. For a given explicit system, there also is a `hidden' condition under which the skeleton structure can be built: the take off/touch down points $(u^\pm_\ast,0,v_\ast,q_\ast)$ must be on $\M^\pm_0$, i.e. $\va \in (v^-_-, v^-_+) \cap (v^+_-, v^+_+)$ . Although this must be the case by construction, in a given system a manifold $\M^\pm_0$ may extend beyond one of its end/transition points $(f^\pm(v^\pm_\pm), v^\pm_\pm)$ as non-normally hyperbolic manifold -- see for instance (\ref{d:Mpm0-FOTM}) in section \ref{sss:FOTM} -- and then it may happen that (the projections of) $\W^{\rm u}((\bV^-,0))$ and $\W^{\rm s}((\bV^+,0))$ intersect at a $v$-coordinate $\va \notin (v^-_-, v^-_+) \cap (v^+_-, v^+_+)$. In that case at least one of the points $(u^\pm_\ast,0,v_\ast,q_\ast)$ will not be a saddle for the fast flow (\ref{e:fast-red}) so that the fast jump between $\M^\pm_0$ cannot be made (by a monotonic fast connection $\ua(\xi)$ -- see Remark \ref{r:normallyhyperbolic}). 
\\ \\
We conclude that the existence problem for (\ref{e:RDE}) differs from that of (\ref{e:RDE-scalar}) in the sense that there always is a uniquely determined traveling connecting (monotonic, planar) interface in (\ref{e:RDE-scalar}) (when there are two stable background states), while we found that this certainly is not the case for singularly perturbed 2-component model (\ref{e:RDE}). Nevertheless, the conditions on $F(U,V)$ and $G(U,V)$ are quite mild and the slow-fast-slow interfaces $(U_h(\xi),U_h(\xi))$ can be constructed for many models -- as we shall see in section \ref{s:Models}.

\begin{remark}
\label{r:normallyhyperbolic}
\rm
As in \cite{Doe22}, we consider a slow manifold of (\ref{e:DS}) to be normally hyperbolic if it is normally hyperbolic for all $c$ with $c \tau = \O(1)$. The dashed parts of $\{F(U,V) = 0\}$ in Fig. \ref{f:manifold&interface}(a) are normally hyperbolic for $c \neq 0$, but elliptic for $c = 0$, and thus not considered as normally hyperbolic in the present analysis. Orbits that (for $c \neq 0$) either take off or touch down from branches of $\{F(U,V) = 0\}$ that are not normally hyperbolic in the present sense, i.e. that are elliptic for $c = 0$, will oscillate and thus correspond to fronts with non-monotonic fast jumps. Such fronts are expected to be (longitudinally) unstable. Note that by the above arguments, the critical points associated to the stable background states $(\bU^\pm,\bV^\pm)$ necessarily must be on invariant manifolds that are normally hyperbolic for all $c = \O(1)$.
\end{remark}

\begin{remark}
\label{r:slowpatterns}
\rm
The savanna ecosystem model of \cite{vLetal03} only has one (ecologically relevant) normally hyperbolic slow manifold that typically contains critical points (for the slow flow) associated to several different stable homogeneous background states \cite{Doe22}. Thus, the model of \cite{vLetal03} violates assumption {\bf (E-II)} and cannot exhibit the slow-fast-slow interfaces considered here. However, if {\bf (E-II)} does not hold, there will be heteroclinic connections between saddle points on the slow manifold that do not depart from the slow manifold. The existence, stability and bifurcations of such `fully slow' localized patterns is studied in \cite{Doe22}.
\end{remark}

\subsection{Transversal instability of the interface}
\label{ss:stability}

To determine the (spectral) stability of the traveling interface $(U_h(x-c_\ast t), V_h(x-c_\ast t))$ as solution of (\ref{e:RDE}), we consider perturbations of the form,
\beq
\label{e:stab-pert}
(U(x,y,t),V(x,y,t)) = (U_h(\xi) + \bu(\xi) e^{\la t + i \ell y}, V_h(\xi) + \bv(\xi) e^{\la t + i \ell y}), \; \; \la=\la(\ell) \in \CC, \ell \in \RR,
\eeq
and derive the linearized system,
\begin{equation}
\label{e:Lin-ODE}
\left\{	
\begin{array}{rcrclcccccc}
\tau \la \bu &=& \bu_{\xi\xi} & - & \ell^2 \bu & + & c_\ast \tau \bu_\xi & + & F_u(u_h(\xi),v_h(\xi)) \bu & + & F_v(u_h(\xi),v_h(\xi))\bv\\
\la \bv &=& \frac{1}{\eps^2}\bv_{\xi\xi} & - & \frac{\ell^2}{\eps^2} \bv & + & c_\ast \bv_{\xi} & + & G_u(u_h(\xi),v_h(\xi)) \bu & + & G_v(u_h(\xi),v_h(\xi)) \bv
\end{array}
\right.
\end{equation}
which can equivalently be written as
\beq
\label{e:Lin-Operator}
\LL(\xi)
\left(
\begin{array}{c}
\bu
\\
\bv
\end{array}
\right)
=
\la
\left(
\begin{array}{c}
\bu
\\
\bv
\end{array}
\right)
+
\ell^2
\left(
\begin{array}{c}
\frac{1}{\tau}\bu
\\[1mm]
\frac{1}{\eps^2}\bv
\end{array}
\right)
\eeq
where the $2 \times 2$ matrix operator $\LL(\xi)$ is given by
\beq
\label{d:LL}
\LL(\xi) =
\left(
\begin{array}{cc}
\L_u(\xi) & \frac{1}{\tau} F_v(u_h(\xi),v_h(\xi))
\\
G_u(u_h(\xi),v_h(\xi)) & \L_v(\xi)
\end{array}
\right)
\eeq
with Sturm-Liouville operators
\beq
\label{d:Luv}
\L_u(\xi) = \frac{1}{\tau}\left(\frac{d^2}{d \xi^2} + c_\ast \tau \frac{d}{d \xi} +  F_u(u_h(\xi),v_h(\xi)) \right), \;
\L_v(\xi) = \frac{1}{\eps^2} \frac{d^2}{d \xi^2} + c_\ast \frac{d}{d \xi} +  G_v(u_h(\xi),v_h(\xi)).
\eeq
We know by assumption {\bf (E-I)} on the stability of the endstates $(\bU^\pm,\bV^\pm)$ of the interface that the essential spectrum of $\LL(\xi)$ is in the stable complex half-plane, moreover we know by the translational invariance of (\ref{e:RDE}) that $\la(0) = 0$ is an eigenvalue of (\ref{e:Lin-Operator}) at $\ell =0$ with eigenfunction $(\bu(\xi),\bv(\xi)) = (u_{h, \xi}(\xi), v_{h, \xi}(\xi))$, the `derivative of the wave'. The interface $(U_h(\xi),V_h(\xi))$ is spectrally stable if (\ref{e:Lin-ODE})/(\ref{e:Lin-Operator}) only has integrable (nontrivial) solutions $(\bu(X), \bv(X))$ for Re$(\la(\ell)) \leq 0$ \cite{KP13}.
\\ \\
As stated in the beginning of this section (and in the Introduction), we also impose assumption {\bf (S-I)}: the interface is stable against 1-dimensional perturbations, i.e. longitudinal perturbations that do not depend on $y$. In other words, we assume that all elements of the discrete spectrum of $\LL(\xi)$ have Re$(\la(0)) < 0$, except for the -- under this assumption -- critical translational eigenvalue $\la_c(0) = 0$. Moreover, we additionally impose the (generic) assumption that the dimension of the kernel of $\LL$ is equal to 1, i.e. that the parameters $\vmu$ of (\ref{e:RDE}) are chosen away from their (co-dimension 1) bifurcation sets -- see also section \ref{s:Disc}. Naturally, these assumptions need to be validated in the setting of an explicit model -- something that typically can be done by exploiting the singularly perturbed nature of (\ref{e:RDE}) by the methods of \cite{BCD19,CdRS16,DIN04,DV15,NF87,NMIF90,Ward18} and the references therein (although one should not underestimate the technicalities involved).
\\ \\
To analyze the impact of the $y$-direction, and thus the possibility of transversal instabilities, we consider the $\ell$-dependence of the critical eigenvalue $\la_c(\ell)$. As explained in the Introduction, we only consider long wavelength perturbations along the interface, i.e. $0 < \ell \ll 1$. Recalling that the spectral problem (\ref{e:Lin-ODE})/(\ref{e:Lin-Operator}) varies in $\ell^2$ instead of $\ell$, we expand
\beq
\label{d:la2bu2bv2}
\la_c(\ell) = \la_{c,2} \ell^2 + \O(\ell^4),
\; \;
\left(
\begin{array}{c}
\bu_c(\xi;\ell^2)
\\
\bv_c(\xi;\ell^2)
\end{array}
\right)
=
\left(
\begin{array}{c}
u_{h, \xi}(\xi)
\\
v_{h, \xi}(\xi)
\end{array}
\right)
+
\left(
\begin{array}{c}
\bu_{c,2}(\xi)
\\
\bv_{c,2}(\xi)
\end{array}
\right)
\ell^2
+
\O(\ell^4).
\eeq
By substitution of (\ref{d:la2bu2bv2}) into (\ref{e:Lin-Operator}) we find at $\O(\ell^2)$ the inhomogeneous equation,
\[
\LL \left(
\begin{array}{c}
\bu_{c,2}(\xi)
\\
\bv_{c,2}(\xi)
\end{array}
\right)
=
\la_{c,2}
\left(
\begin{array}{c}
u_{h, \xi}(\xi)
\\
v_{h, \xi}(\xi)
\end{array}
\right)
+
\left(
\begin{array}{c}
\frac{1}{\tau} u_{h, \xi}(\xi)
\\[1mm]
\frac{1}{\eps^2} v_{h, \xi}(\xi)
\end{array}
\right)
\]
Since $\LL(\xi)$ has a nontrivial kernel, the inhomogeneous term needs to satisfy a solvability/Fredholm condition \cite{KP13}. Clearly $\LL$ is not self-adjoint (and neither are its diagonal elements $\L_u$ and $\L_v$ (\ref{d:Luv})). We therefore define the adjoint operator $\LL^A(\xi)$,
\beq
\label{d:LLA}
\LL^A(\xi) =
\left(
\begin{array}{cc}
\frac{1}{\tau}\left(\frac{d^2}{d \xi^2} - c_\ast \tau \frac{d}{d \xi} +  F_u(u_h(\xi),v_h(\xi)) \right) & G_u(u_h(\xi),v_h(\xi))
\\
\frac{1}{\tau} F_v(u_h(\xi),v_h(\xi)) &  \frac{1}{\eps^2} \frac{d^2}{d \xi^2} - c_\ast \frac{d}{d \xi} +  G_v(u_h(\xi),v_h(\xi))
\end{array}
\right)
\eeq
and its eigenfunction $(\bu_c^A(\xi),\bv_c^A(\xi))$ at $\la_c^A = 0$, i.e. the integrable (nontrivial) solution of $\LL^A(\bu,\bv) = (0,0)$. It follows that,
\beq
\label{e:solvability}
\left<
\la_{c,2}
\left(
\begin{array}{c}
u_{h, \xi}(\xi)
\\
v_{h, \xi}(\xi)
\end{array}
\right)
+
\left(
\begin{array}{c}
\frac{1}{\tau} u_{h, \xi}(\xi)
\\[1mm]
\frac{1}{\eps^2} v_{h, \xi}(\xi)
\end{array}
\right)
,
\left(
\begin{array}{c}
\bu^A_c(\xi)
\\
\bv^A_c(\xi)
\end{array}
\right)
\right>
= 0,
\eeq
which yields an explicit expression for $\la_{2,c}$,
\beq
\label{e:la2c}
\la_{2,c} = -\frac{\int_\RR \left(\frac{1}{\tau} u_{h, \xi}(\xi) \bu^A_c(\xi) + \frac{1}{\eps^2} v_{h, \xi}(\xi) \bv^A_c(\xi) \right) \, d \xi}{\int_\RR \left(u_{h, \xi}(\xi) \bu^A_c(\xi) + v_{h, \xi}(\xi) \bv^A_c(\xi) \right) \, d \xi}.
\eeq
The planar interface $(U_h(\xi),V_h(\xi))$ is unstable with respond to long wavelength perturbations along the interface if $\la_{2,c} > 0$. It follows from (\ref{e:la2c}) that we can explicitly characterize this potential destabilization mechanism by determining (a leading order approximation in $\eps$ of) $(\bu_c^A(\xi),\bu_c^A(\xi))$.
\\ \\
To do so, we rewrite $\LL^A(\bu,\bv) = (0,0)$ in a form similar to (\ref{e:Lin-ODE}) and split this equation up into two slow parts and one fast part. The slow parts are at leading order given by,
\beq
\label{e:Lin-ODE-slow}
\left\{	
\begin{array}{rcccrcr}
\eps^2 \bu_{XX} & - & \eps c_\ast \tau \bu_X & + & F_u(u^{\pm}(X),v^{\pm}(X)) \bu & = & - \tau G_u(u^{\pm}(X),v^{\pm}(X))\bv\\
\bv_{XX} & - & \eps c_\ast \bv_{X} & + & G_v(u^{\pm}(X),v^{\pm}(X)) \bv & = & - \frac{1}{\tau} F_v(u^{\pm}(X),v^{\pm}(X)) \bu
\end{array}
\right.
\eeq
where we have approximated the slow parts $(u_h(X),v_h(X)$ by $(u^{\pm}(X),v^{\pm}(X))$ and suppressed the $u$-, $s$-subscripts (i.e. $u^\pm(X)$ is either $=u^-_u(X)$ or $u^+_s(X)$, etc.). The fast part is (at leading order) given by
\beq
\label{e:Lin-ODE-fast}
\left\{	
\begin{array}{rcccrcr}
\bu_{\xi\xi} & - & c_\ast \tau \bu_\xi & + & F_u(u_\ast(\xi),v_\ast) \bu & = & - \tau G_u(u_\ast(\xi),v_\ast)\bv\\
\bv_{\xi\xi} & = & \eps^2 \left[ c_\ast \bv_{X} \right.& - & G_v(u_\ast(\xi),v_\ast) \bv & - & \left.\frac{1}{\tau} F_v(u_\ast(\xi),v_\ast) \bu \right]
\end{array}
\right.
\eeq
now with $u_h(X) = u_\ast(\xi)$, $v_h(X) \equiv v_\ast$ as leading order approximations during the fast jump. By taking $\eps \to 0$ in the $\bu$-equation of (\ref{e:Lin-ODE-slow}) we see that $\bu$ is determined by $\bv$ in the slow fields,
\[
\bu = - \tau \frac{G_u(u^{\pm},v^{\pm})}{F_u(u^{\pm},v^{\pm})} \bv
\]
(at leading order). By $(f^\pm)'(v) = - F_v(f^{\pm}(v),v)/F_u(f^{\pm}(v),v)$, we thus have
\[
\bv_{XX} +  \left[G_u(u^{\pm},v^{\pm}) + G_v(u^{\pm},v^{\pm})(f^\pm)'(v^{\pm})\right]\bv = 0
\]
in the slow fields (as $\eps \to 0$), which implies by (\ref{e:slow-red}) that $\bv^{\pm}(X) = \al^\pm v^{\pm}_X(X)$ for some constants $\al^\pm$. Taking $\eps \to 0$ in the fast field yields that $\bv(\xi) \equiv \bv_\ast$, a constant, and that
\[
\bu_{\xi\xi} - c_\ast \tau \bu_\xi + F_u(u_\ast(\xi),v_\ast) \bu  = \L_u^A \bu = - \tau G_u(u_\ast(\xi),v_\ast)\bv_\ast
\]
(\ref{d:Luv}). Since the kernel of $\L_u$ is spanned by $u_{\ast,\xi}$ this implies that
\beq
\label{e:solvLuA}
\tau \bv_\ast \int_\RR G_u(u_\ast(\xi),v_\ast) u_{\ast,\xi} \, d \xi =  \tau \bv_\ast \left[G(u_\ast^+,v_\ast) - G(u_\ast^-,v_\ast) \right] = 0
\eeq
(at leading order in $\eps)$. Thus, under the (generic) assumption that $G(u_\ast^+,v_\ast) \neq G(u_\ast^-,v_\ast)$, we conclude that $\bv$ must be $\O(\eps)$, so that $\bv^{\pm} = \eps \bal^\pm v^{\pm}_X = \bal^\pm v^{\pm}_\xi$. As a consequence, $\L_u^A \bu = 0$ (at leading order), which implies that
\beq
\label{e:bu-ff}
\bu(\xi) = \al_\ast u_{\ast, \xi}(\xi) \, e^{c_\ast \tau \xi}
\eeq
for some $\al_\ast$. By matching over the fast field, we can eliminate 2 of the 3 free factors $\bal^\pm$ and $\al_\ast$. First we note that $\bv$ must be continuous over the fast field so that $\bal^+ = \bal^-= \bal$. Since $v^-_X(0) = q_\ast = v^+_X(0)$ (Fig. \ref{f:slowflows}) it also follows that $\bv_\ast = \eps \bal q_\ast$. Recalling that $\bv = \O(\eps)$, we measure the accumulated change in $\bv_\xi$ over the fast field by integrating $\bv_{\xi \xi}$,
\[
\Delta_f \bv_\xi = - \eps^2 \frac{\al_\ast}{\tau} \int_\RR F_v(u_\ast(\xi),v_\ast) u_{\ast, \xi}(\xi) e^{c_\ast \tau \xi} \, d\xi
\]
(\ref{e:Lin-ODE-fast}), (\ref{e:bu-ff}), which must bridge the jump in the $\bv_\xi$'s coming from the slow flows,
\[
\Delta_s \bv_\xi =
\lim_{X \downarrow 0} \bv_\xi - \lim_{X \downarrow 0} \bv_\xi =
\bal \left[\lim_{X \downarrow 0} v^+_{\xi\xi} - \lim_{X \uparrow 0} v^-_{\xi\xi}\right] =
\eps^2 \bal \left[\lim_{X \downarrow 0} v^+_{XX} - \lim_{X \uparrow 0} v^-_{XX}\right] =
- \eps^2 \bal \left[G(u_\ast^+,v_\ast) - G(u_\ast^-,v_\ast) \right]
\]
(\ref{e:slow-red}), so that
\beq
\label{e:alast-tal}
\frac{\bal}{\al_\ast} = \frac{1}{\tau}\frac{\int_\RR F_v(u_\ast,v_\ast) u_{\ast, \xi} \, e^{c_\ast \tau \xi} \, d\xi}{G(u_\ast^+,v_\ast) - G(u_\ast^-,v_\ast)}.
\eeq
(again assuming that $G(u_\ast^+,v_\ast) \neq G(u_\ast^-,v_\ast)$). Now we can evaluate (\ref{e:la2c}), although we need to be a bit more precise about our choice for the fast field,
\[
\int_\RR  v_{h, \xi} \bv^A_c \, d \xi =
\left(\int_{-\infty}^{-\frac{1}{\sqrt{\eps}}} + \int_{-\frac{1}{\sqrt{\eps}}}^{\frac{1}{\sqrt{\eps}}} + \int_{-\frac{1}{\sqrt{\eps}}}^{\infty} \right) v_{h, \xi} \bv^A_c \, d \xi =
\bal \int_{-\infty}^{-\frac{1}{\sqrt{\eps}}} (v^-_\xi)^2 \, d \xi +
\eps^2 \bal \int_{-\frac{1}{\sqrt{\eps}}}^{\frac{1}{\sqrt{\eps}}} q_\ast^2 \, d \xi +
\bal \int_{\frac{1}{\sqrt{\eps}}}^{\infty} (v^+_\xi)^2 \, d \xi
\]
i.e.
\beq
\label{e:intvvA}
\int_\RR  v_{h, \xi} \bv^A_c \, d \xi =
\eps \bal \left(\int_{-\infty}^0 (v^-_X)^2 \,dX + \int_0^{\infty} (v^+_X)^2 \,dX \right) + \O(\eps \sqrt{\eps}).
\eeq
In contrast, the integral over $u_{h, \xi} \bu^A_c$ is dominated by its fast part and is of $\O(1)$,
\beq
\label{e:intuuA}
\int_\RR  u_{h, \xi} \bu^A_c \, d \xi =
\al_\ast \int_\RR (u_{\ast, \xi})^2 e^{c_\ast \tau \xi} \,d \xi + \O(\eps)
\eeq
(\ref{e:bu-ff}). Thus, we find by (\ref{e:la2c}) and (\ref{e:alast-tal}),
\beq
\label{e:la2c-expl}
\la_{2,c} = -\frac{1}{\eps \tau} \, \frac{\int_\RR F_v(u_\ast,v_\ast) u_{\ast, \xi} \, e^{c_\ast \tau \xi} \, d\xi}{G(u_\ast^+,v_\ast) - G(u_\ast^-,v_\ast)} \, \frac{\int_{-\infty}^0 (v^-_X)^2 \,dX + \int_0^{\infty} (v^+_X)^2 \,dX}{\int_\RR (u_{\ast, \xi})^2 e^{c_\ast \tau \xi} \,d \xi}
\eeq
(at leading order in $\eps$), which indeed yields (\ref{e:cond-longtrans}) (with $\Fa$ and $\Ga$ as in (\ref{d:FGast})).

\section{Taking $\tau$ asymptotically small: transitions at $\tau = \eps \ttau$}
\label{s:SmallTau}

In the introduction, the impact of having $\tau = \O(\eps)$ was briefly sketched. This sketch will be worked out in more detail in the beginning of section \ref{ss:ExStab-SmallTau}. Now that we have derived the leading order approximation of $\la_{2,c}$ (\ref{e:la2c-expl}), we can also sketch how the stability analysis of section \ref{ss:stability} is affected by setting $\tau = \eps \ttau$. A crucial ingredient in the derivation of (\ref{e:la2c-expl}) is the observation that $\bv^A_c$ must be $\O(\eps)$, which follows from the application of the leading order solvability condition (\ref{e:solvLuA}). However, as $\tau$ becomes $\O(\eps)$, this condition is satisfied for $\bv^A_c = \O(1)$. As a consequence, the integral over $v_{h, \xi} \bv^A_c$ is no longer $\O(\eps)$, which has a significant impact on (\ref{e:la2c-expl}) and thus on (\ref{e:cond-longtrans}), as we shall see in the upcoming section.
\\ \\
Unlike for $\tau = \O(1)$, i.e. criterion (\ref{e:cond-longtrans}), instability criterion (\ref{e:cond-longtrans-t-Intro}) -- with $\tM_\ast(\ttau)$ to be derived in the upcoming subsection -- has become an explicit function of $\ttau$. Thus, the sign of $\la_{2,c}$ may change as $\ttau$ is varied. In section \ref{ss:BifTrav-SmallTau} we show that such a change of sign generates a bifurcation.
\begin{remark}
\label{r:normhypsaddles}
\rm
The three subconditions of condition (\ref{e:cond-stabbUVpm}) that establish the stability of the background states $(\bU^\pm,\bV^\pm)$ of (\ref{e:RDE}) reduces to two for $\tau = \eps \ttau$, :  $\bF_u^\pm \bG_v^\pm - \bF_v^\pm \bG_u^\pm > 0$ and  $\bF_u^\pm < 0$. These conditions are equivalent to the condition that the associated critical points of $(\bU^\pm,0,\bV^\pm,0)$ of (\ref{e:DS}) correspond to saddle points for the slow reduced flows (\ref{e:slow-red}) of the slow manifolds $\M^\pm_0$. In other words, if $\tau = \eps \ttau$ then each saddle point for a slow flow on (normally hyperbolic slow manifold) $\M^\pm_0$ necessarily corresponds to a stable homogeneous background state (which is not the case for $\tau = \O(1)$).
\end{remark}

\subsection{Existence and transversal (in)stability}
\label{ss:ExStab-SmallTau}

\begin{figure}[t]
\centering
\includegraphics[width =0.9 \linewidth]{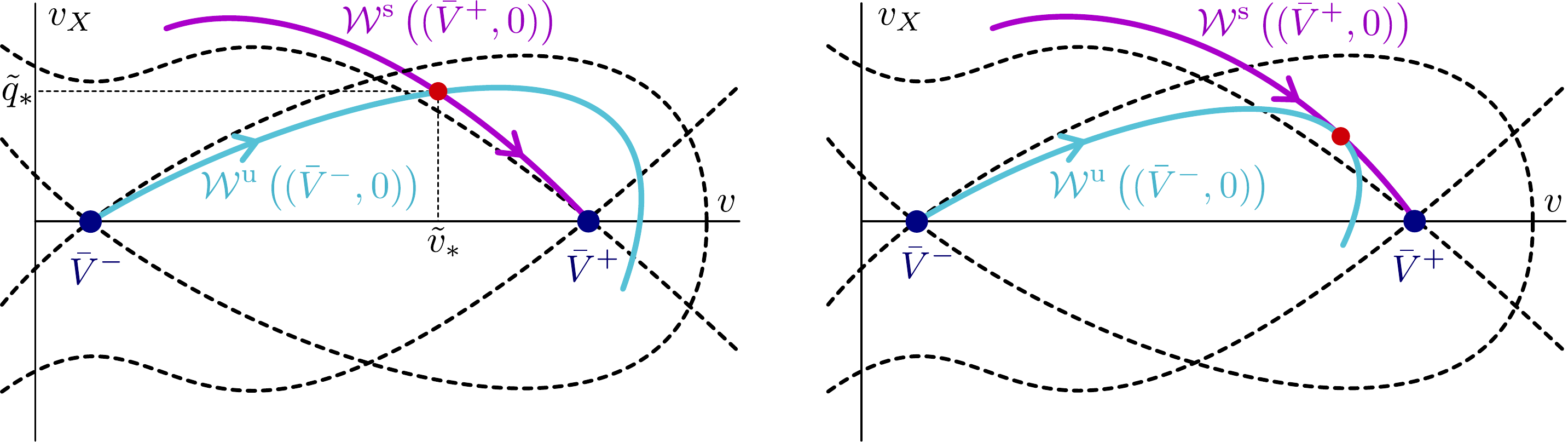}
\caption{\small{Sketches of combined projected slow flows (\ref{e:slow-red-t}) on $\M^\pm_0$ for $\tau = \eps \ttau$, $c = \tc/\eps$ and $\tc \neq 0$, with the $\tau = \O(1)$, i.e. $\tc = 0$, case of Fig. \ref{f:slowflows} indicated by dashed black lines. On the right, $\tc$ is at a critical bifurcational value of $\tc$: beyond this value $\W^{\rm u}((\bV^-,0))$ and $\W^{\rm s}((\bV^+,0))$ no longer intersect.}}
\label{f:slowflows-t}
\end{figure}

As argued in the Introduction, we follow $\tau = \eps \ttau$ by setting $c = \tc/\eps$, so that the `friction term' in  (\ref{e:fast-red}) remains $\O(1)$,
\beq
\label{e:fast-red-t}
u_{\xi\xi} + \tc \ttau u_\xi +  F(u,v_0)  =  0.
\eeq
For a given $v_0$ there thus still is a well-defined value of $\tc = \C(v_0)/\ttau$ for which there is an heteroclinic orbit between the saddles $(f^\pm(v_0),0)$ of (\ref{e:fast-red-t}), i.e. for which the fast system enables a jump between the slow manifolds $\M^-_0$ and $\M^+_0$. However, the slow reduced problem is no longer integrable in this case,
\beq
\label{e:slow-red-t}
v_{XX} + \tc v_X + G(f^\pm(v),v)  =  0 \; \; {\rm on} \; \; \M^\pm_0
\eeq
(cf. (\ref{e:slow-red})), so that the coordinate $(\tv_\ast,\tq_\ast)$ of the intersection of the unstable manifold $\W^{\rm u}((\bV^-,0))$ and the stable manifold $\W^{\rm s}((\bV^+,0))$ of the projected flows on $\M^\pm_0$ -- see Fig. \ref{f:slowflows-t}(a) -- becomes $\tc$-dependent (unlike for $\tau = \O(1)$, Fig. \ref{f:slowflows}): $\tv_\ast = \V(\tc)$. Thus, $(\tc_\ast,\tv_\ast)$ is now determined by
\beq
\label{e:cast-t}
\tc \tau = \C\left(\V(\tc)\right), \; \tv = \V(\tc),
\eeq
which may have more than 1 solution. As for $\tau = \O(1)$, $\tu_\ast(\xi)$ is defined as the (monotonic) heteroclinic orbit of (\ref{e:fast-red-t}) at $\tc = \tc_\ast$ (with $\tu^\pm_\ast = f^\pm(\tv_\ast)$). As in section \ref{s:Tau=O1}, we assume that  these solutions can all be constructed (and refer once more to the examples in section \ref{s:Models}). The components based on $\tv^-(X)$ (determined by slow flow on $\M^-_0$), $\tu_\ast(\xi)$ (the fast jump) and $\tv^+(X)$ (on $\M^-_0$), together form the skeleton of the traveling interface $(U(x,y,t),V(x,y,t)) = (\tU_h(x-\tc_\ast t/\eps), \tV_h(x-\tc_\ast t/\eps)) = (\tu_h(\xi),\tv_h(\xi))$.
\\
\\
For the stability analysis we again assume that the front $(\tU_h(x-\tc_\ast t/\eps), \tV_h(x-\tc_\ast t/\eps))$ is longitudinally stable. To investigate the stability with respect to $y$-dependent perturbations, we define
\beq
\label{d:LL-t}
\tLL(\xi) =
\left(
\begin{array}{cc}
\frac{1}{\eps \ttau}\left(\frac{d^2}{d \xi^2} + \tc_\ast \ttau \frac{d}{d \xi} +  F_u(\tu_h,\tv_h) \right) &  \frac{1}{\eps \tau} F_v(\tu_h,\tv_h)
\\
G_u(\tu_h,\tv_h) &  \frac{1}{\eps^2} \frac{d^2}{d \xi^2} + \frac{\tc_\ast}{\eps} \frac{d}{d \xi} +  G_v(\tu_h,\tv_h)
\end{array}
\right)
\eeq
and derive the equivalent of (\ref{e:la2c}),
\beq
\label{e:la2c-t}
\tla_{2,c} = -\frac{\int_\RR \left(\frac{1}{\eps \ttau} \tu_{h, \xi}(\xi) \bu^A_c(\xi) + \frac{1}{\eps^2} \tv_{h, \xi}(\xi) \bv^A_c(\xi) \right) \, d \xi}{\int_\RR \left(\tu_{h, \xi}(\xi) \bu^A_c(\xi) + \tv_{h, \xi}(\xi) \bv^A_c(\xi) \right) \, d \xi}
\eeq
with now $(\bu^A_c, \bv^A_c)$ a solution of,
\beq
\label{e:Lin-ODE-slow-t}
\left\{	
\begin{array}{rcccrcr}
\eps^2 \bu_{XX} & - & \eps \tc_\ast \ttau \bu_X & + & F_u(\tu^{\pm},\tv^{\pm}) \bu & = & - \eps \ttau G_u(\tu^{\pm},\tv^{\pm})\bv\\
\bv_{XX} & - & \tc_\ast \bv_{X} & + & G_v(\tu^{\pm},\tv^{\pm}) \bv & = & - \frac{1}{\eps \ttau} F_v(\tu^{\pm},\tv^{\pm}) \bu
\end{array}
\right.
\eeq
(cf. (\ref{e:Lin-ODE-slow})) in the slow fields and of,
\beq
\label{e:Lin-ODE-fast-t}
\left\{	
\begin{array}{rcccrcl}
\bu_{\xi\xi} & - & \tc_\ast \ttau \bu_\xi & + & F_u(\tu_\ast,\tv_\ast) \bu & = & - \eps \ttau G_u(\tu_\ast,\tv_\ast)\bv\\
\bv_{\xi\xi} & = & \eps \left[ \tc_\ast \bv_{\xi} \right.& - & \frac{1}{\ttau} F_v(\tu_\ast,\tv_\ast) \bu  & + & \left. \O(\eps) \right]
\end{array}
\right.
\eeq
(cf. (\ref{e:Lin-ODE-fast})) in the fast field. Following the approach of the $\tau = \O(1)$-case we find (at leading order),
\[
\bv^\pm = \tal \tv^\pm_X \, e^{\tc_\ast X}, \;
\bu = \tal_\ast \tu_{\ast, \xi}(\xi) \, e^{\tc_\ast \ttau \xi}
\]
in the slow fields and confirm that indeed $\tal = \O(1)$ (and not $\O(\eps)$ as for $\tau = \O(1)$). The matching procedure also proceeds as for $\tau = \O(1)$ and yields an outcome almost identical to (\ref{e:alast-tal}), $\tal_\ast \tF_\ast = \tal \ttau \tG_\ast$ with
\beq
\label{d:tFGast}
\tF_\ast = \int_\RR F_v(\tu_\ast,\tv_\ast) \tu_{\ast, \xi} \, e^{\tc_\ast \ttau \xi} \, d\xi, \;
\tG_\ast = G(\tu_\ast^+,\tv_\ast) - G(\tu_\ast^-,\tv_\ast).
\eeq
To further simplify notation, we introduce
\beq
\label{d:tIsf}
\tI_f = \int_\RR (\tu_{\ast, \xi})^2 e^{\tc_\ast \ttau \xi} \,d \xi, \; \;
\tI_s = \int_{-\infty}^0 (\tv^-_X)^2 e^{\tc_\ast X} \,dX + \int_0^{\infty} (\tv^+_X)^2e^{\tc_\ast X} \,dX,
\eeq
and the Melnikov-type expression introduced in the Introduction,
\beq
\label{d:tMa}
\tM_\ast(\vmu; \ttau) = \tF_\ast(\vmu; \ttau) \tI_s(\vmu; \ttau) + \ttau \tG_\ast(\vmu; \ttau) \tI_f(\vmu; \ttau).
\eeq
Substitution of all leading order approximations into (\ref{e:la2c-t}) leads to a leading order expression that differs significantly from (\ref{e:la2c-expl}),
\beq
\label{e:la2c-expl-t}
\tla_{2,c}(\vmu; \ttau) = -\frac{1}{\eps^2} \, \frac{\tF_\ast \tI_s}{\tM_\ast},
\eeq
which indeed yields instability criterion (\ref{e:cond-longtrans-t-Intro}). Nevertheless, (\ref{e:la2c-expl-t}) and (\ref{e:cond-longtrans-t-Intro}) do take the form of (\ref{e:la2c-expl}) and (\ref{e:cond-longtrans}) for $\ttau \gg 1$ (\ref{d:tFGast}), (\ref{d:tIsf}). Moreover, if sign$(\tF_\ast)$ $=$ sign$(\tG_\ast)$, i.e. if the front is stable against transversal long wavelength perturbations for $\tau = \O(1)$ (\ref{e:cond-longtrans}), then (\ref{e:cond-longtrans-t-Intro}) indicates that it will remain so for any value of $\tau$/$\ttau$ (see section \ref{sss:Cylindrical}).
\\ \\
However, in a given system it may not be possible to take $\ttau$ sufficiently small, since a bifurcation may take place as $\ttau$ is decreased. The origin of this bifurcation may be geometrical: it follows from (\ref{e:fast-red-t}) that $|\tc|$ must become large as $\ttau$ becomes small to have a jump through the fast field. However, the intersection of unstable manifold $\W^{\rm u}((\bV^-,0))$ and stable manifold $\W^{\rm s}((\bV^+,0))$ of the projected flows on $\M^\pm_0$ may only persist up to a certain maximal value of $|\tc|$ -- see Fig. \ref{f:slowflows-t}(b): the interface may for instance have `disappeared' in a saddle-node bifurcation. Moreover, in the case that we know that the interface is unstable with respect to transversal perturbations for $\tau = \O(1)$, i.e. if sign$(\tF_\ast)$ $\neq$ sign$(\tG_\ast)$ for $\ttau$ sufficiently large (\ref{e:cond-longtrans}), then the sign of $\tla_{2,c}$ may also change at the pole of leading order approximation (\ref{e:la2c-expl-t}) associated to $\M_\ast = \tF_\ast \tI_s + \ttau \tG_\ast \tI_f$ crossing through 0 for decreasing $\ttau$ (\ref{e:la2c-expl-t}). We shall consider this in the upcoming section. (Here it should be noted that a pole in the leading order approximation (\ref{e:la2c-expl-t}) does not imply that $\la_{2,c}$ itself is non-smooth: one needs to perform a higher order analysis to resolve the leading order character of $\la_{2,c}$ near the zero of the denominator in (\ref{e:la2c-expl-t}).) 

\subsection{Bifurcation into traveling waves}
\label{ss:BifTrav-SmallTau}
As was already noted in the Introduction, the expression $\tM_\ast(\vmu; \ttau)$ (\ref{d:tMa}) has a Melnikov function-like character, especially for $\tc_\ast = 0$ (see (\ref{d:tFGast}), (\ref{d:tIsf})). Therefore, it is natural to expect that a zero of $\tM_\ast(\vmu; \ttau)$ may be related to a `structural geometrical change' in the existence analysis. In this subsection we consider the special case of stationary fronts, i.e. fronts for which the existence problem does not depend on $\ttau$. We show that in this case a change of the sign of $\tM_\ast(\vmu; \ttau)$ indeed has a direct impact on the existence problem: it initiates a bifurcation into traveling waves, i.e. as $\tM_\ast(\vmu; \ttau)$ crosses through 0 -- for instance as function of $\ttau$ -- a pair of counter-moving traveling fronts branches off from the stationary fronts.
\\ \\
Existence problem (\ref{e:DS}) with $\tau = \eps \ttau$ and $c = \tc/\eps$ has a stationary front solution if (\ref{e:fast-red-t}) has a heteroclinic orbit for $\tc =0$ for a value of $v_0 = \tva$ which also is the $v$-coordinate of the intersection of the unstable manifold $\W^{\rm u}((\bV^-,0))$ and the stable manifold $\W^{\rm s}((\bV^+,0))$ of the projected flows (\ref{e:slow-red-t}) on $\M^\pm_0$ with $\tc = 0$. In other words, $\C(\tv)$ and $\V(\tc)$ of (\ref{e:cast-t}) must be so that $\C\left(\V(0)\right) = 0$, which imposes a co-dimension 1 condition on the parameters $\vmu$ of (\ref{e:RDE}) -- see below. We assume that $\vmu$ indeed has values so that $\C\left(\V(0)\right) = 0$ and thus use the fact that (\ref{e:RDE}) has a stationary front solution for any $\tau = \eps \ttau > 0$. Given this set-up, we ask the question: under which additional condition can there be slowly traveling fronts, i.e. fronts traveling with speed $\tc = \hc \delta$ with (artificial) small parameter $0 < \eps \ll \delta \ll 1$ (that is independent of $\eps$), such that these fronts merge with the stationary front as $\delta \to 0$?
\\ \\
We first consider the fast problem (\ref{e:fast-red-t}). We have assumed that there is a heteroclinic orbit $\tua(\xi)$ between the critical points $(\tua^-,0)$ and $(\tua^+,0)$ -- with $\tua^\pm = f^\pm(\tva)$ -- in the integrable case $\hc = 0$. Thus, we define the family of Hamiltonians,
\beq
\label{d:Ham-fast}
\tHH_f(u,p; v_0) = \frac12 p^2 + \int_{\tua^-(v_0)}^u F(s,v_0) \, ds,
\eeq
parameterized by $v_0$. The assumption on the existence of a stationary front implies
\beq
\label{d:statfronts}
\tHH_f(\tua^-,0; \va) = \tHH_f^\ast(\tua^+,0; \va) = \int_{\tua^-}^{\tua^+} F(s,\tva) \, ds = 0,
\eeq
see also \cite{NF87,TN94} and note that (\ref{d:statfronts}) is the explicit version of the expression $\C(\tva) = 0$ (\ref{e:cast-t}). Naturally, one expects that if a traveling front appears from the stationary front, it will not make its fast jump at exactly the same value of $\tva$, therefore we introduce $\hv_1$ by $\tva(\tc) = \tva + \delta \hv_1$ (with $\hv_1 = \hv_1(\hc)$ and recall that $\tc = \delta \hc$). As a consequence, the $u$-coordinates of the 2 saddle points $(\tu^\pm_\ast(\tc), 0) = (f^\pm(\tva(\tc)),0)$ of (\ref{e:fast-red-t}) also undergo an $\O(\delta)$ shift,
\[
f^\pm(\tva(\tc)) = \tua^\pm + \delta (f^\pm)'(\tva) \hv_1 + \O(\delta^2) =
\tua^\pm - \delta \frac{F_v(\tua^\pm,\tva)}{F_u(\tua^\pm,\tva)} \hv_1 + \O(\delta^2),
\]
where we note that $\tva(0) = \tva$ and $\tu^\pm_\ast(0) = \tu^\pm_\ast$. By (\ref{d:Ham-fast}) and (\ref{d:statfronts}) we have,
\[
\begin{array}{rcl}
\tHH_f(f^-(\tva(\tc)),0;\tva(\tc)) & = & \int_{\tua^-(\tc)}^{\tua^-(\tc)} F(s,\tva(\tc)) \, ds = 0
\\
\tHH_f(f^+(\tva(\tc)),0;\tva(\tc)) & = & \int_{\tua^-(\tc)}^{\tua^+(\tc)} F(s,\tva + \delta \hv_1 + \O(\delta^2)) \, ds 
\\
& = & \int_{\tua^-(0)}^{\tua^+(0)} F(s,\tva + \delta \hv_1 + \O(\delta^2)) \, ds + \O(\delta^2) 
\\
& = & \int_{\tua^-}^{\tua^+} F(s,\tva) \, ds + \delta \hv_1 \int_{\tua^-}^{\tua^+} F_v(s,\tva) \, ds + O(\delta^2)
\\
& = & \delta \hv_1 \int_{\tua^-}^{\tua^+} F_v(s,\tva) \, ds + O(\delta^2).
\end{array}
\]
The $\tc \neq 0$-heteroclinic jump between $(f^-(\tva(\tc)),0)$ and $(f^+(\tva(\tc)),0)$ thus needs to accumulate a total change in `energy' of
\[
\Delta \tHH_f = \delta \hv_1 \int_{\tua^-}^{\tua^+} F_v(s,\tva) \, ds + O(\delta^2).
\]
Naturally, at $v_0 = \tva(\hv)$, $\tHH_{f, \xi} =  - \delta \hc \ttau u^2_\xi$ (\ref{e:fast-red-t}), (\ref{d:Ham-fast}), so that
\[
\Delta \tHH_f = - \delta \hc \ttau \int_{\RR} (\tu_{\ast, \xi})^2 \, d \xi + \O(\delta \sqrt{\delta}),
\]
where $\tu_{\ast}(\xi)$ is the heteroclinic orbit at $\tva = \tva(0)$. Thus, we conclude by (\ref{d:tFGast}), (\ref{d:tIsf}) that,
\beq
\label{e:tv1-fast}
\hv_1(\hc) = - \frac{\ttau \int_{\RR} (\tu_{\ast, \xi})^2 \, d \xi}{\int_{\tua^-}^{\tua^+} F_v(s,\tva) \, ds} \, \hc + \O(\sqrt{\delta}) = - \frac{\ttau \tI_f(\tc = 0)}{\tF_{\ast}(\tc = 0)} \, \hc + \O(\sqrt{\delta}).
\eeq
Since the slow reduced flows on manifold $\M^\pm_0$ are also integrable as $\delta \to 0$, we define the slow Hamiltonians,
\beq
\label{d:Ham-slow}
\tHH_s^\pm(v,q) = \frac12 q^2 + \int_{\bV^\pm}^v G(f^\pm(s), s) \, ds,
\eeq
so that $\tHH_s^\pm(v,q) = 0$ at the saddle points $(\bV^\pm,0)$ and at the intersection $(\tva,\tqa)$. Note that $\tva$ is determined by,
\beq
\label{d:V(0)}
\int_{\bV^-}^{\tva} G(f^-(s), s) \, ds = \int_{\bV^+}^{\tva} G(f^+(s), s) \, ds,
\eeq
which thus yields the explicit version of $\V(0)$ (\ref{e:cast-t}): together, (\ref{d:statfronts}) and (\ref{d:V(0)}) determine the co-dimension 1 condition on the parameters $\vmu$, $\C(\V(0)) = 0$. Since again $\tHH_{s, X}^\pm = - \hc \delta q^2$ (\ref{e:slow-red-t}), we  observe by a similar Melnikov argument as above that the unstable manifold $\W^{\rm u}((\bV^-,0))$ on $\M^-_0$ and the stable manifold $\W^{\rm s}((\bV^+,0))$ on $\M^+_0$ of the ($\delta$-)perturbed $\hc \neq 0$-flows intersect the lines $v = \tva (\subset \M^\pm_0)$ at the level sets
\beq
\label{d:DeltatHspm}
\begin{array}{rcccr}
\tHH_s^-(\tva,\tqa + \delta \hq^-) & \stackdef & - \delta \Delta \tHH_s^- & = & - \delta \hc \int_{-\infty}^0 (\tv^-_X(X))^2 \, dX + \O(\delta^2)
\\
\tHH_s^+(\tva,\tqa + \delta \hq^+) & \stackdef & \delta \Delta \tHH_s^+ & = &  \delta \hc \int_0^{\infty} (\tv^+_X(X))^2 \, dX + \O(\delta^2)
\end{array}
\eeq
Thus, the $\O(\delta)$ corrections to the $q$-coordinates of the intersections of $\W^{\rm u}((\bV^-,0))$ and $\W^{\rm s}((\bV^+,0))$ with the vertical $\{v=\tva \}$ are determined by,
\[
\hq^\pm = \pm \frac{\Delta \tHH_s^\pm}{\tqa} \hc + \O(\delta).
\]
As a consequence, we know that $\O(\delta)$ close to $(\tva,\tqa)$, $\W^{\rm u}((\bV^-,0))$ and $\W^{\rm s}((\bV^+,0))$ are up to $\O(\delta^2)$ corrections given by the linear approximations,
\[
\Gamma^\pm(s^\pm) = \left\{
\left(
\begin{array}{c}
\tva \\ \tqa \pm \frac12 \delta \frac{\Delta \tHH_s^\pm}{2 \tqa} \hc
\end{array}
\right)
+
\delta s^\pm
\left(
\begin{array}{c}
\tqa \\ -G(f^\pm(\tva),\tva)
\end{array}
\right),
\; s^\pm = \O(1)
\right\}.
\]
Naturally, these two lines intersect and their intersection gives the leading order approximation of the intersection of (the projections of) $\W^{\rm u}((\bV^-,0))$ and $\W^{\rm s}((\bV^+,0))$ for the perturbed flows (\ref{e:slow-red-t}) with $\tc = \delta \hc$. By the above definition, this intersection has $v$-coordinate $\tva(\hc) = \tva + \delta \hv_1$ and we conclude by (\ref{d:tFGast}), (\ref{d:tIsf}) that,
\beq
\label{e:tv1-slow}
\hv_1 = \frac{\Delta \tHH_s^+ + \Delta \tHH_s^-}{G(\tua^+,\tva) - G(\tua^-,\tva)} \, \hc =
\frac{\tI_s(\tc = 0)}{\tG_\ast} \, \hc
\eeq
(at leading order in $\delta$). Combining (\ref{e:tv1-slow}) with (\ref{e:tv1-fast}) indeed yields $\hc \tM_\ast(\ttau)  = \O(\delta)$ (\ref{d:tMa}): traveling fronts with non-zero speed $\tc = \delta \tc$ branch off for $\ttau$ such that $\tilde{\lambda}_{2,c}$ changes sign by crossing through the pole in the leading order approximation (\ref{e:la2c-expl-t}).
\\ \\
Note that the analysis so far only provides a first approximation: the line $\ttau = \ttau_\ast = - \tF_\ast \tI_s/\tG_\ast \tI_f$ in the $(\ttau, \tc)$-plane only gives a leading order vertical approximation of the (parabolic) curve of bifurcating traveling waves (again in the $(\ttau, \tc)$-plane) that is expected for a non-degenerate bifurcation. More specifically, the orientation of this curve determines whether the bifurcation is sub- or supercritical: with the present analysis we do not know yet the precise nature of the bifurcation (and under which condition(s) it is non-degenerate). This more detailed information can be determined explicitly by a higher order asymptotic analysis (see for instance \cite{Doe22}), however, we refrain from going into this here. Moreover, it is natural to ask whether the bifurcating traveling interfaces will be stable against transversal (long wavelength) perturbations or not: clearly the sign of $\tM_\ast$ changes for the basic stationary pattern as $\ttau$ passes through $\ttau_\ast$, but the sign of $\tM_\ast$ is a priori not clear for the bifurcating traveling fronts. Such a higher order analysis can also be performed as we show in section \ref{sss:InStabBifTrav} in the context of a specific model.

\section{Ecosystem models and other example systems}
\label{s:Models}
The analysis of this paper was inspired by the insights of \cite{RBBKBD21} on the positive impact of spatial patterns on the resilience of ecosystems -- see section \ref{s:Disc} -- and of \cite{FOTM19} on potentially reversing desertification through a fingering instability of an invasion front. Therefore, we exemplify the general approach by first considering two explicit ecosystem models in sections \ref{sss:BCDE} and \ref{sss:FOTM}. However, the ecosystem models considered typically do not have a (small) $\tau$-pre-factor so that we next consider a FitzHugh-Nagumo-type model as  studied in \cite{HM94,HM97a,HM97b}: first in section \ref{ss:FHN} for $\tau = \O(1)$, next for $\tau = \eps \ttau$ in section \ref{sss:InStabBifTrav}. In section \ref{sss:Cylindrical}, we consider the relation between the present paper and \cite{Tan03,TN94}, in which the stability of planar interfaces on cylindrical domains is studied for models of the type (\ref{e:RDE}), especially by analyzing the explicit prototypical example model considered in \cite{Tan03,TN94}.

\subsection{Fingering invasion fronts in ecosystems}
\label{ss:fingering}

\subsubsection{Invasion fronts in a Klausmeier-type dryland ecosystem model}
\label{sss:BCDE}
\begin{figure}[t]
\hspace{.02\textwidth}
\begin{subfigure}{.52 \textwidth}
\centering
\includegraphics[width=1\linewidth]{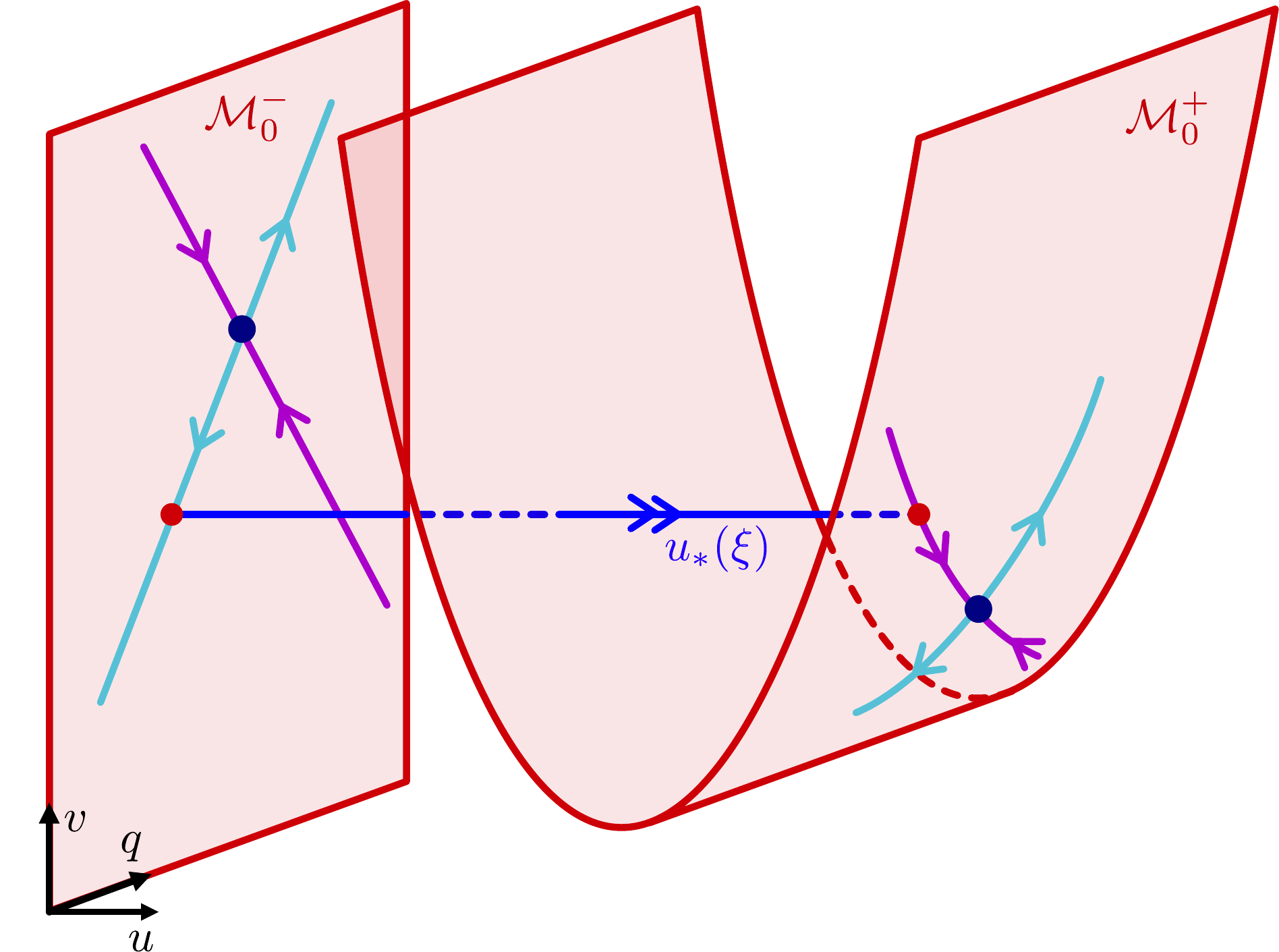}
\caption{}
\end{subfigure}
\hspace{.05\textwidth}
\begin{subfigure}{.38 \textwidth}
\centering
\includegraphics[width=1\linewidth]{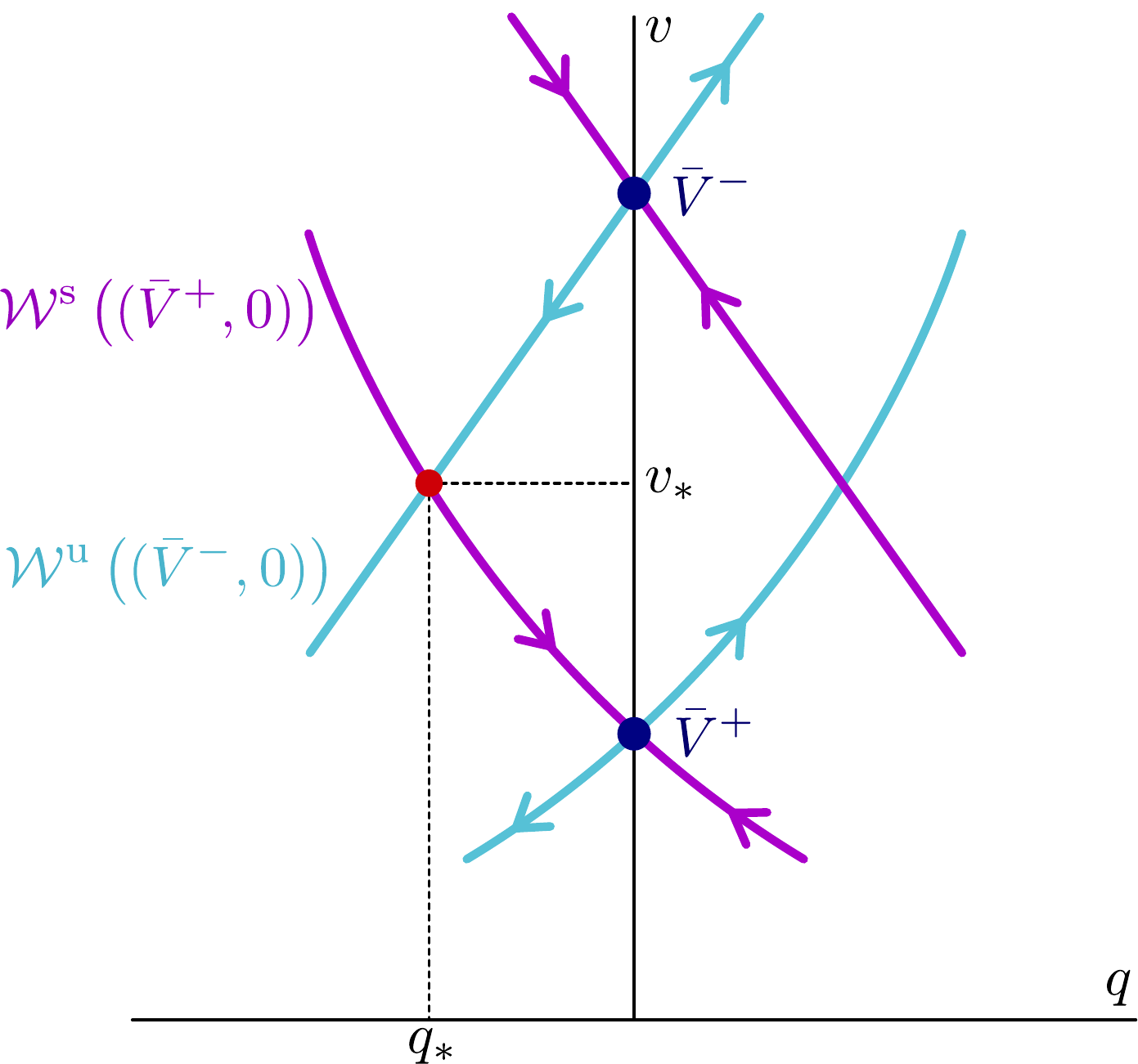}
\caption{}
\end{subfigure}
\caption{\small{The `skeleton' for the invasion front in dryland ecosystem model (\ref{e:RDE-BCDE}) as sketch in 3-dimensional $(u,v,q)$ space (a); (b): the `vertical' combined projections of slow flows (\ref{e:slow-red-BCDE}) on $\M^\pm_0$ in the $(q,v)$-plane with the explicit indication of the (uniquely determined) jumping point $(\qa,\va) \in \W^{\rm u}((\bV^-,0)) \cap \W^{\rm s}((\bV^+,0))$   -- cf. Fig. \ref{f:slowflows}.}}
\label{f:BCDE}
\end{figure}

On flat terrains, the ecosystem model proposed and studied in \cite{BCD19} coincides with the model of \cite{Eig21} reduced to one (instead of two) species,
\begin{equation}
\label{e:RDE-BCDE}
\left\{	
\begin{array}{rcl}
U_t &=& \, \, \, \, \, \, \Delta U -\mu_1 U + U^2 (1- \mu_2 U) V \\
V_t &=& \frac{1}{\varepsilon^2} \Delta V + \mu_3 - V - U^2 V
\end{array}
\right.
\end{equation}
i.e. $F(U,V;\vmu) = -\mu_1 U + U^2 (1- \mu_2 U) V$ and $G(U,V;\vmu) = \mu_3 - V - U^2 V$ in (\ref{e:RDE}). In this model, $U(x,y,t)$ represents vegetation or biomass and $V(x,y,t)$ water (or `the limiting resource' \cite{Eig21}) -- that diffuses much faster (i.e. the assumption $0 < \eps \ll 1$ is natural in this model \cite{RvdK08}). As in the Klausmeier model \cite{Kla99,Sher10} and its reaction-diffusion modification introduced in \cite{SDHR13}, $\mu_1 > 1$ represents the mortality rate of the vegetation and $\mu_3 > 0$ the (yearly) precipitation rate. Through its zero $U = 1/\mu_2$ (i.e. $\mu_2 > 0$), the additional factor $(1- \mu_2 U)$ in the $U$-equation represents the carrier capacity of the vegetation (in absence of mortality).
\\ \\
From the ecological point of view it is natural to consider invasion fronts of the bare soil state, i.e. $(U(x,y,t),V(x,y,t)) \equiv (\bU^-,\bV^-) = (0, \mu_3)$, into a homogeneously vegetated state $(U(x,y,t),V(x,y,t)) \equiv (\bU^+,\bV^+)$ -- or vice versa. The fast reduced system for traveling waves, i.e. the equivalent of (\ref{e:fast-red}), is given by the cubic, planar system,
\beq
\label{e:fast-red-BCDE}
u_{\xi\xi} + c u_\xi - \mu_1 u + v_0 u^2 (1- \mu_2 u) =  0, \; \;
v = v_0, q = q_0
\eeq
that defines two normally hyperbolic slow manifolds,
\beq
\label{d:Mpm0-BCDE}
\M^-_0 = \{u=f^-(v)\equiv 0, p= 0\}, \; \; \M^+_0 = \{u=f^+(v)=u^+(v), p=0, v> 4 \mu_1 \mu_2 \}
\eeq
with $u^\pm(v)$ the 2 solutions of $\mu_2 v_0 u^2 - v_0 u + \mu_1 = 0$,
\beq
\label{d:upm-BCDE}
u^\pm(v) = \frac{1}{2\mu_2} \pm \frac{1}{2\mu_2} \sqrt{1 - \frac{4 \mu_1 \mu_2}{v}}.
\eeq
Clearly, the (reduced) slow flow on $\M^-_0$ is linear, $v_{XX} + \mu_3 - v = 0$, while the flow on $\M^+_0$ is given by
\beq
\label{e:slow-red-BCDE}
v_{XX} + \mu_3 + \frac{\mu_1}{\mu_2} = \left(1 + \frac{1}{2 \mu_2^2}\right)v + \frac{1}{2 \mu_2^2} \sqrt{v^2 - 4\mu_1\mu2 v} = 0.
\eeq
Since the right hand side of (\ref{e:slow-red-BCDE}) is an increasing function of $v$, it follows that (\ref{e:slow-red-BCDE}) has 1 critical point $(v,q) = (\bV^+,0)$ of saddle type if $\mu_3 > (1 + 4 \mu_2^2) \frac{\mu_1}{\mu_2}$. Naturally, this saddle point corresponds to the stable homogeneously vegetated state $(\bU^+,\bV^+)$ (with $\bU^+ = f^+(\bV^+)$). Note that it follows from evaluating the right hand side of (\ref{e:slow-red-BCDE}) at $v = \mu_3$ that $\bV^+ < \mu_3 = \bV^-$ (see Fig. \ref{f:BCDE}(a)): since the vegetation takes up water, the water level of the vegetated state is lower than that of the bare soil state. In their combined projections, the (linear) unstable manifold $\W^{\rm u}((\mu_3,0))$ on $\M^-_0$ must clearly intersect the stable manifold $\W^{\rm s}((\bV^+,0))$ on $\M^+_0$ (see Fig. \ref{f:BCDE}(b)), thereby defining the fast jumping point $(v_\ast,q_\ast) \in \W^{\rm u}((\bV^-,0)) \cap \W^{\rm s}((\bV^+,0))$ with $4\mu_1\mu_2 < \bV^+ < v_\ast < \bV^- = \mu_3$ and $q_\ast < 0$.
\\
\\
It is possible to write down an explicit expression for $v_\ast$. Since this expression plays a crucial role in the companion paper \cite{BCDL22} -- see also section \ref{s:Disc} -- we also briefly formulate it here. By the linearity of the flow on $\M^-_0$, we know that
$(v_\ast,q_\ast) \in \W^{\rm u}((\bV^-,0))$ is given by $(v_\ast, v_\ast - \mu_3)$. By the Hamiltonian associated to integrable system (\ref{e:slow-red-BCDE}) on $\M^+_0$ we thus find,
\[
\frac12(\mu_3-\va)^2 + \int_{\bV^+}^{v_\ast} \left[\left(\mu_3 + \frac{\mu_1}{\mu_2}\right) - \left(1 + \frac{1}{2 \mu_2^2}\right)v - \frac{1}{2 \mu_2^2} \sqrt{v^2 - 4\mu_1\mu_2 v}\right] dv = 0,
\]
which can be simplified into,
\beq
\label{e:condvast-BCDE}
(\mu_3-\bV^+)^2 = \frac{1}{\mu_2^2} \int_{\bV^+}^{v_\ast} \left[\left(v - 2 \mu_1 \mu_2\right) +  \sqrt{(v - 2 \mu_1 \mu_2)^2 - (2\mu_1\mu_2)^2}\right] dv.
\eeq
So far, we have thus explicitly obtained the two slow components of the skeleton structure of the planar slow-fast-slow front $(U_h(\xi),V_h(\xi))$ between the homogeneous background states $(\bU^-,\bV^-)$ and $(\bU^+,\bV^+)$. Its fast component is spanned by the connecting orbit $u_\ast(\xi)$ between the saddle points $(u^-_\ast, 0) = (f^-(v_\ast),0) = (0,0)$ and $(u^+_\ast,0) = (f^+(v_\ast),0)$ in the fast reduced system (\ref{e:fast-red-BCDE}) with $v_0 = v_\ast$ that exists for a specific choice $c_\ast$ of $c$. By the cubic nature of (\ref{e:fast-red-BCDE}), both $c_\ast$ and $u_\ast(\xi)$ can be expressed explicitly in terms of $v_\ast$. Since we will also encounter cubic fast reduced systems in several of the upcoming sections, we present a brief derivation of the formulae for $c_\ast$ and $u_\ast(\xi)$ in a general setting. Let $u_\ast(\xi)$ be a heteroclinic orbit of
\beq
\label{d:fastcubic-gen-2nd}
u_{\xi \xi} + c u_\xi - \al (u-\be_-)(u-\be_c)(u-\be_+) =0, \; \; \al > 0, \; \be_- < \be_c < \be_+
\eeq
between the critical points $(\be_-,0)$ and $(\be_+,0)$ and assume that $\ua(\xi)$ also satisfies,
\beq
\label{d:fastcubic-gen-1st}
u_\xi = - K (u -\be_-)(u-\be_+), \; \; K > 0,
\eeq
so that indeed $\lim_{\xi \to \pm \infty} \ua(\xi) = \be_{\pm}$. By taking the derivative of (\ref{d:fastcubic-gen-1st}) and substitution of (\ref{d:fastcubic-gen-1st}) into (\ref{d:fastcubic-gen-2nd}), we obtain two polynomial expressions (in $u$) for $u_{\xi\xi}$ that coincide for,
\beq
\label{e:Kcast-gen}
K = \frac12 \sqrt{2 \al}, \; \; c_\ast = \sqrt{2 \al} \left(\be_c - \frac12(\be_+ + \be_-) \right).
\eeq
By (\ref{d:fastcubic-gen-1st}) we also have an explicit heteroclinic solution of (\ref{d:fastcubic-gen-2nd}),
\beq
\label{e:uast-gen}
\ua(\xi) =\frac12(\be_+ + \be_-) + \frac12(\be_+ - \be_-) \tanh \left(\frac14 \sqrt{2 \al} (\be_+ - \be_-) (\xi - \xi_0) \right).
\eeq
Applying this general result to fast reduced equation (\ref{e:fast-red-BCDE}) yields $c_\ast = \left(\sqrt{v_\ast}- 3 \sqrt{v_\ast - 4 \mu_1 \mu_2}\right)/\left(2 \sqrt{2 \mu_2}\right)$, which implies that the interface is stationary for $v_\ast = v_\ast(\mu_1,\mu_2,\mu_3) = \frac92 \mu_1 \mu_2$ and that the desert invades the savanna for $\vmu \in \RR^3$ such that $v_\ast < \frac92 \mu_1 \mu_2$ (and vice versa for $v_\ast > \frac92 \mu_1 \mu_2$).
\\ \\
We conclude by the methods of geometric singular perturbation theory that a traveling planar invasion front $(U_h(\xi),V_h(\xi))$ between the bare soil state $(\bU^-,\bV^-)$ and savanna state $(\bU^+,\bV^+)$ indeed exists -- see Fig. \ref{f:BCDE}(c). Here, we do not go into the details of the spectral stability of $(U_h(\xi),V_h(\xi))$ with respect to perturbations that only depend on $x$ (or $\xi$) -- as usual we note that the stability can be established by the methods developed in \cite{BCD19,CdRS16,DIN04,DV15,NF87,NMIF90,Ward18} and the references therein. Thus, as in section \ref{ss:stability} we assume that the front is stable with respect to perturbations that do not depend on $y$ and that the critical eigenvalue curve $\la_c(\ell)$ has $\la_c(0) = 0$. By section \ref{ss:stability} we know that the local parabolic character of $\la_c(\ell)$ is given by (\ref{d:la2bu2bv2}) and (\ref{e:la2c-expl}) and that the (in)stability of the invasion front $(U_h(\xi),V_h(\xi))$ with respect to long wavelength transversal perturbations is determined by (\ref{e:cond-longtrans}).
\\
\begin{figure}[t]
\begin{minipage}{.245\textwidth}
		\centering
		\includegraphics[width =\linewidth]{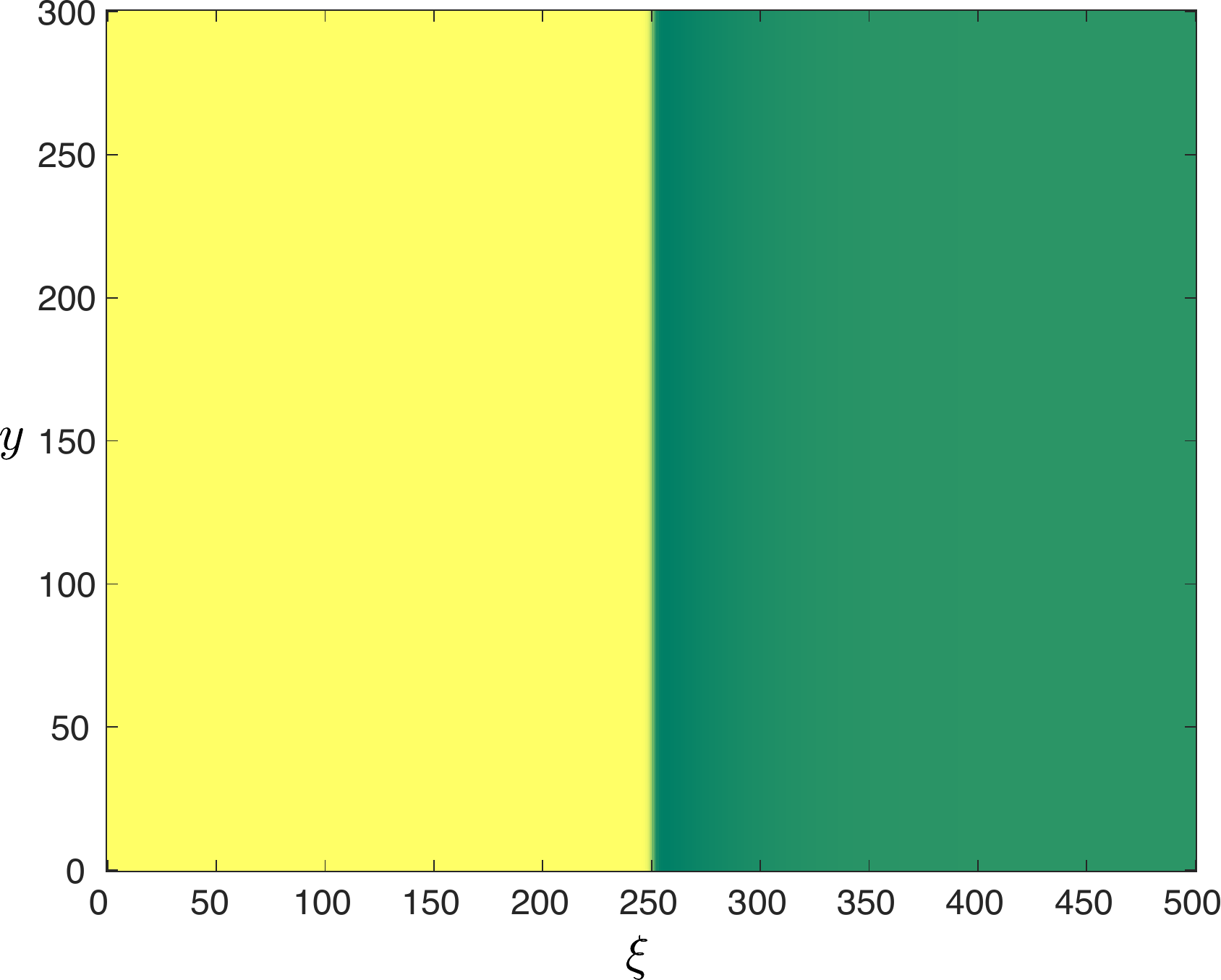}
\end{minipage}
%\hspace{.05cm}
\begin{minipage}{.245\textwidth}
		\centering
		\includegraphics[width =\linewidth]{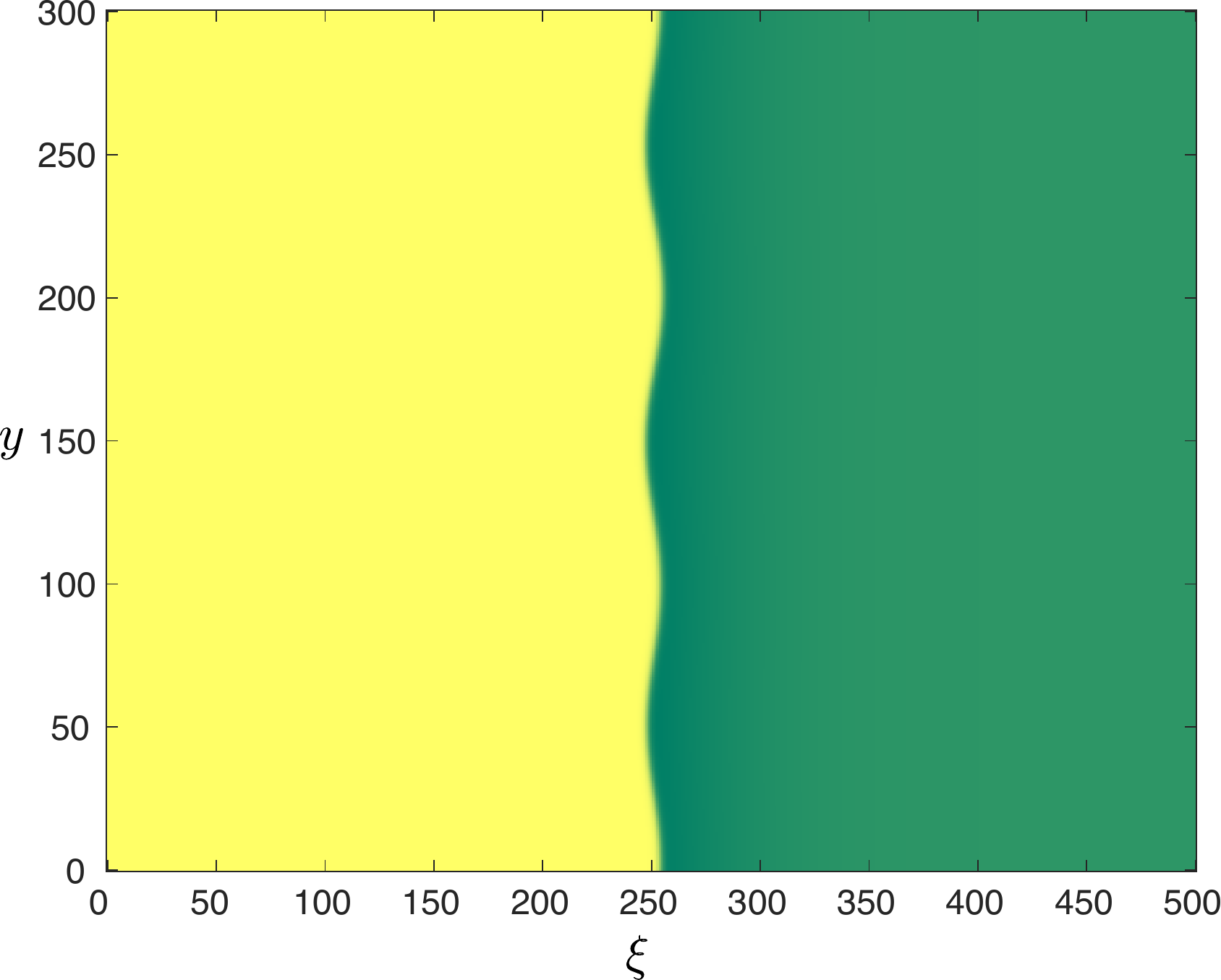}
\end{minipage}
%\hspace{.05cm}
\begin{minipage}{.245\textwidth}
		\centering
		\includegraphics[width =\linewidth]{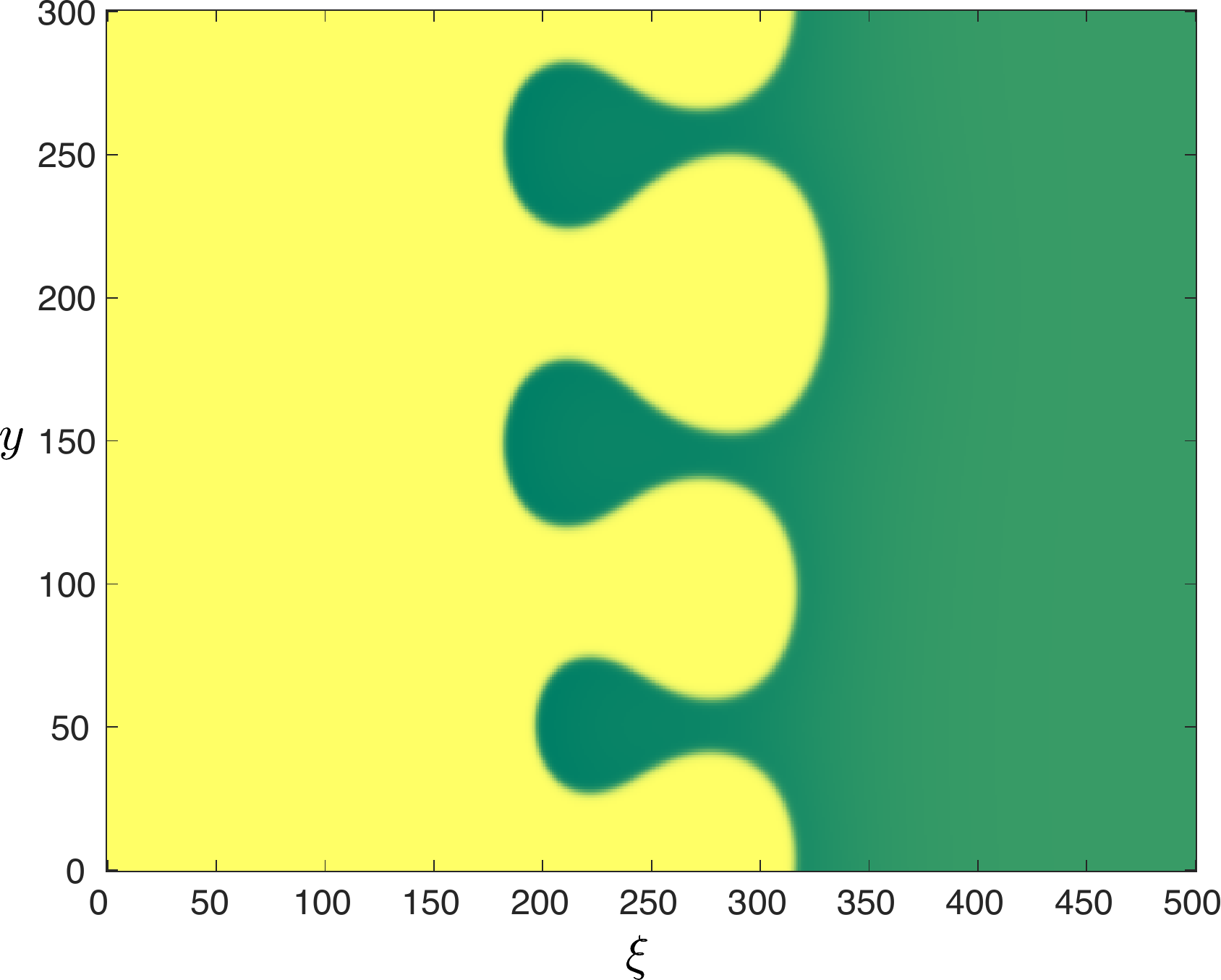}
\end{minipage}
%\hspace{.05cm}
\begin{minipage}{.245\textwidth}
		\centering
		\includegraphics[width =\linewidth]{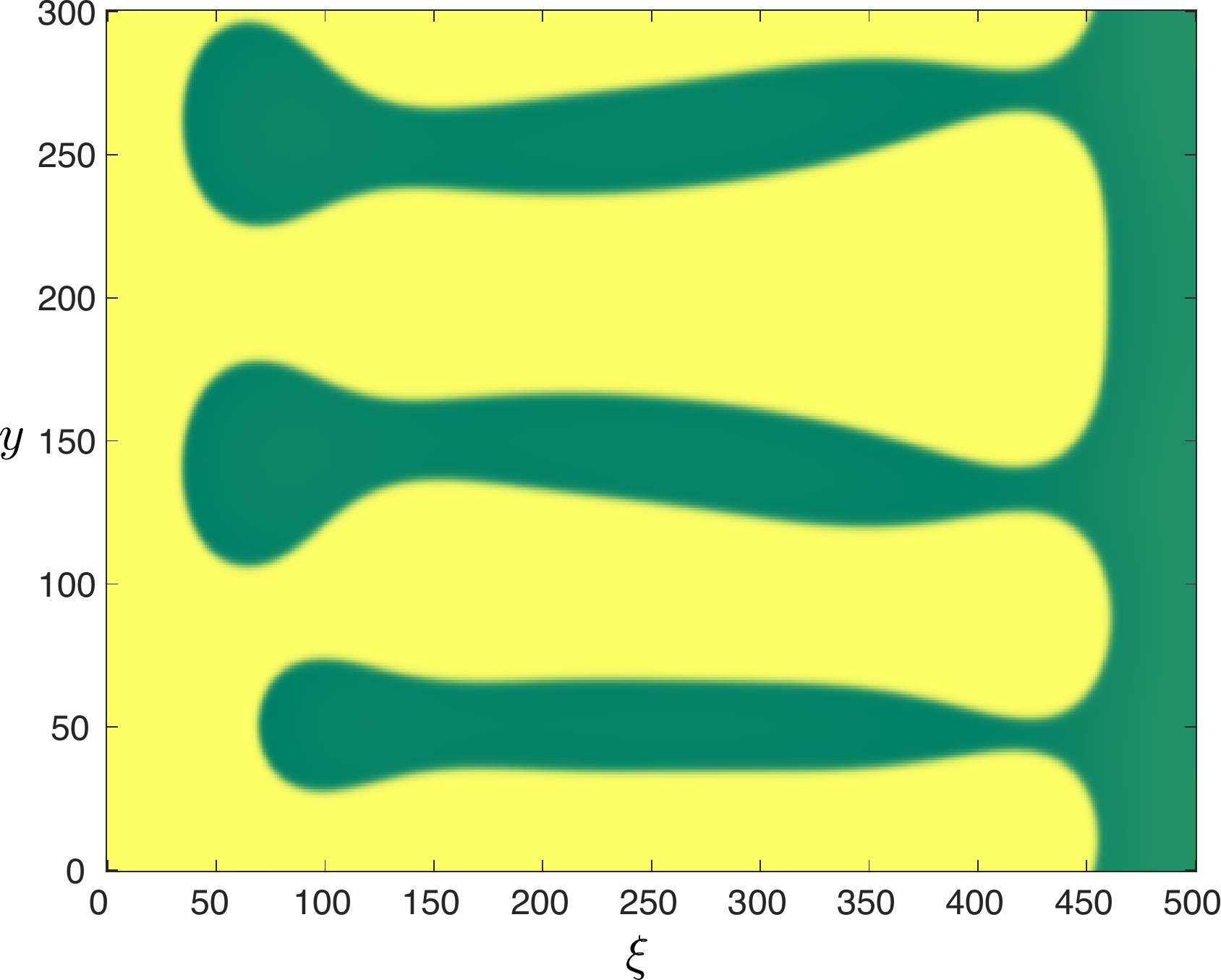}
\end{minipage}
%\\
%\vspace{.2cm}
\\
\begin{minipage}{.245\textwidth}
		\centering
		\includegraphics[width =\linewidth]{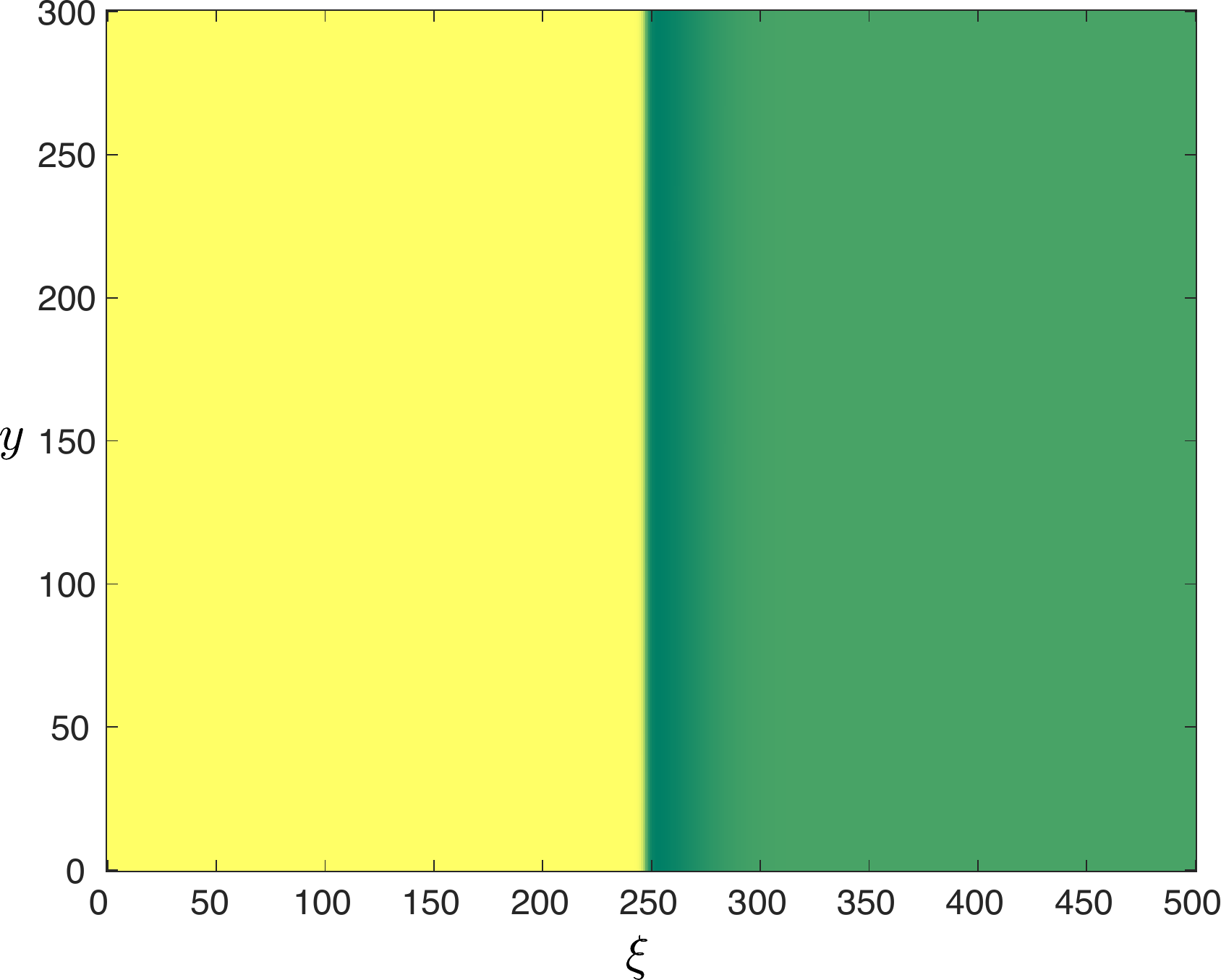}
\end{minipage}
%\hspace{.05cm}
\begin{minipage}{.245\textwidth}
		\centering
		\includegraphics[width =\linewidth]{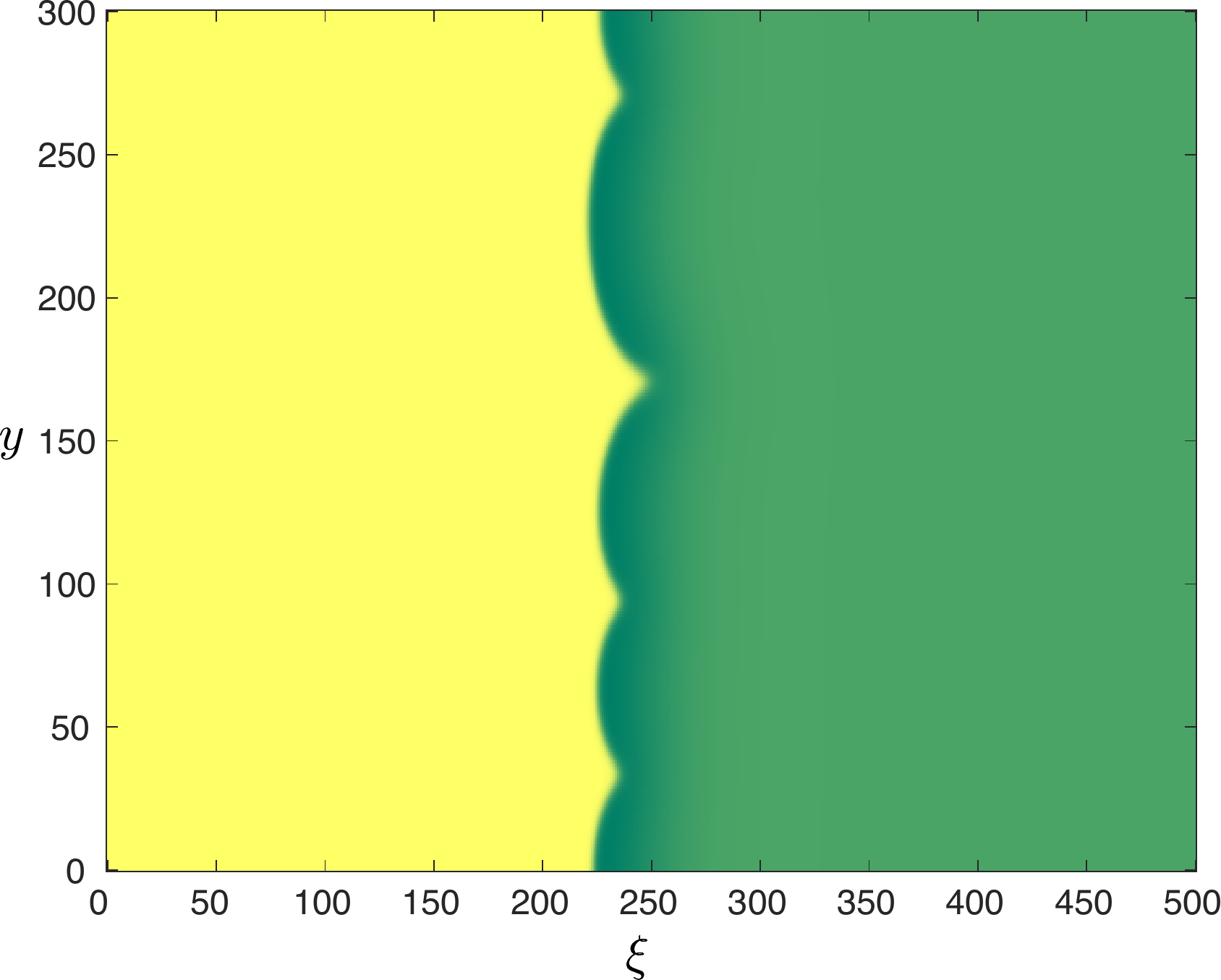}
\end{minipage}
%\hspace{.05cm}
\begin{minipage}{.245\textwidth}
		\centering
		\includegraphics[width =\linewidth]{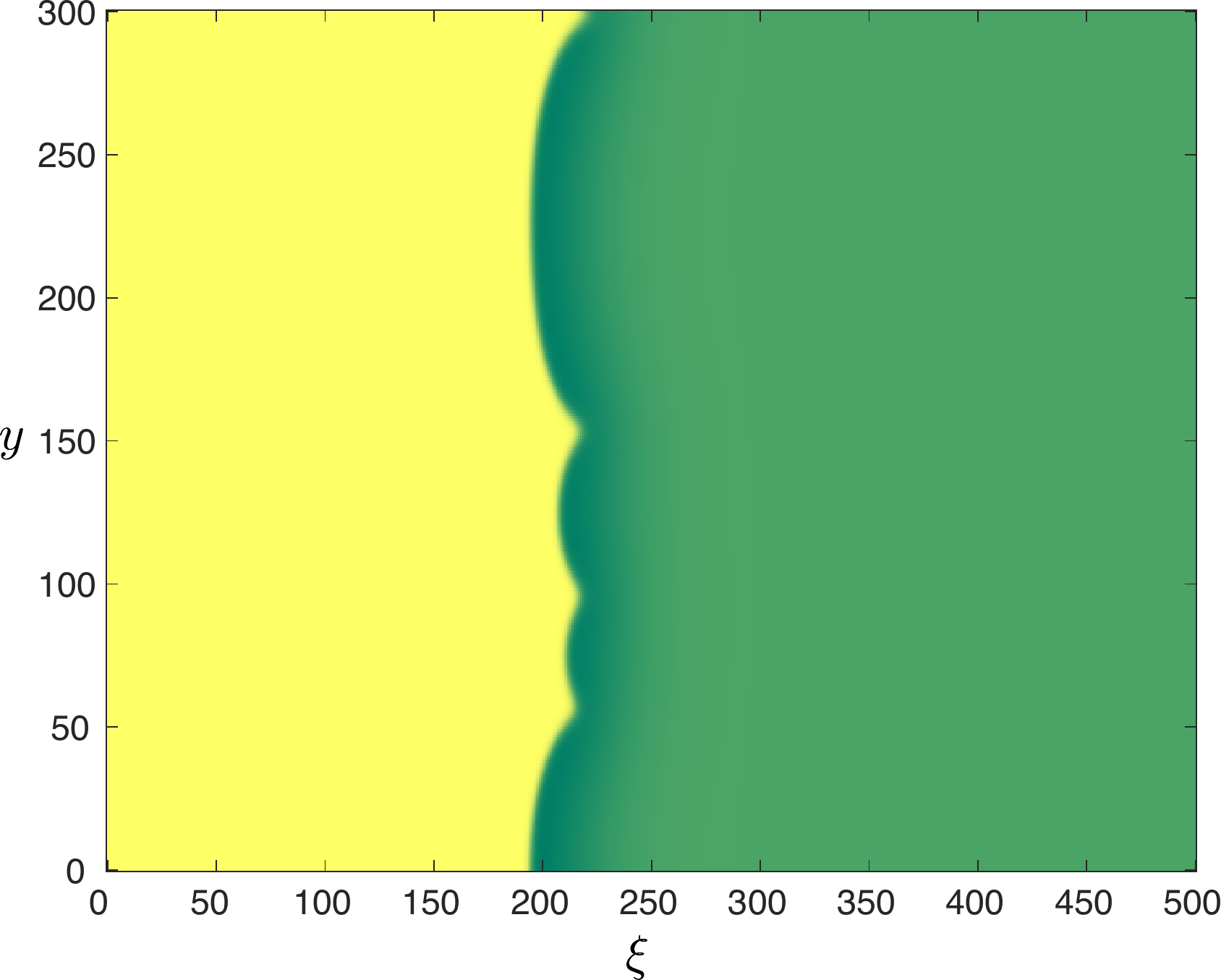}
\end{minipage}
%\hspace{.05cm}
\begin{minipage}{.245\textwidth}
		\centering
		\includegraphics[width =\linewidth]{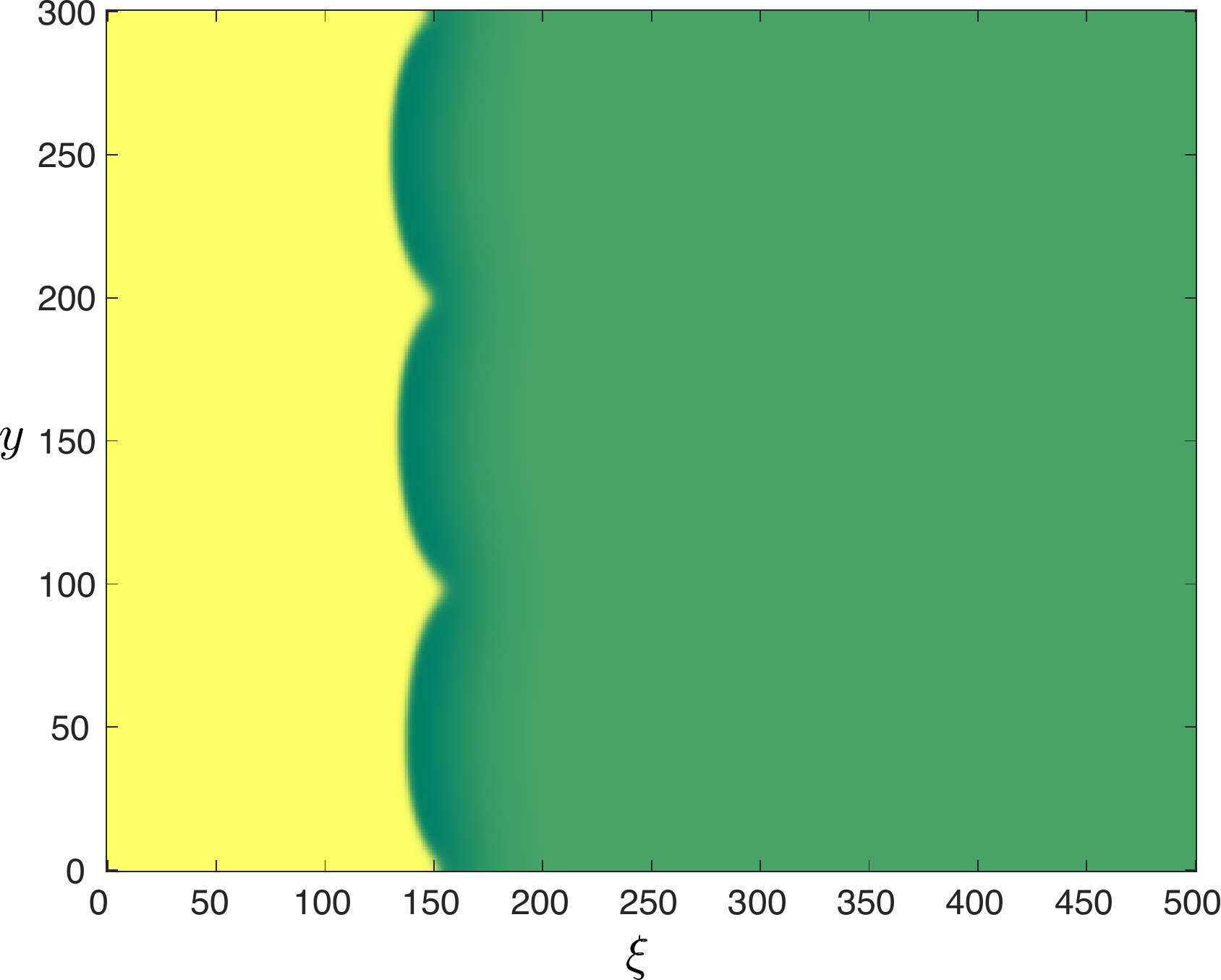}
\end{minipage}
\caption{\small{Results of direct numerical simulations in~(\ref{e:RDE-BCDE}), using finite differences and MATLAB's ode15s routine for time integration. Top row: four snapshots (from left to right, $t=300,2100, 3000, 4200$) of an invading planar desertification front in ecosystem model (\ref{e:RDE-BCDE}) based on \cite{BCD19,Eig21} (with $\vmu =(1.2,1.0,6.2)$, $\eps = 0.02$). As in model (\ref{e:RDE-FOTM}), the transversally unstable interface initiates the counter-invasion of vegetation patterns (\cite{FOTM19}, Fig \ref{f:FingIntro}). Bottom row: four snapshots (from left to right, $t=150,600,900,2100$) of an invading planar vegetation front in (\ref{e:RDE-BCDE}) (with $\vmu =(0.1,0.1,2.0)$, $\eps = 0.01$). The unstable interface does not finger but generates a dynamic interplay of bounded cusps similar to the Kuramoto-Sivashinsky dynamics associated to a sideband instability \cite{DSSS09,HN86}. }}
\label{f:FingCusp-BCDE}
\end{figure}
\\
With $F(U,V;\vmu)$ and $G(U,V;\vmu)$ as in (\ref{e:RDE-BCDE}), we immediately find by the above analysis that,
\[
F_v(u_\ast(\xi), v_\ast) u_{\ast,\xi}(\xi) = \left[u_\ast^2(\xi) (1 - \mu_2 u_\ast(\xi))\right] u_{\ast,\xi}(\xi) > 0,
\]
since $0 < u_\ast(\xi) < u_\ast^+ = f^+(v_\ast) = u^+(v_\ast) < \frac{1}{\mu_2}$ (\ref{d:upm-BCDE}) and $u_{\ast,\xi}(\xi) > 0$ by construction. Moreover,
\[
G(u^+_\ast, v_\ast) - G(u^-_\ast, v_\ast) = \left[\mu_3 - v_\ast - (u^+_\ast)^2 v_\ast \right] - \left[\mu_3 - v_\ast \right] = - (u^+_\ast)^2 v_\ast < 0,
\]
from which we conclude that the invasion fronts in dryland ecosystem model (\ref{e:RDE-BCDE}) cannot be stable with respect to transversal perturbations. Like in model (\ref{e:RDE-FOTM}), the invading desertification front is typically counter-invaded by fingering vegetation patterns. However, the transversally unstable invasion fronts do not necessarily finger: an invading vegetation front may generate dynamically evolving bounded cusps typical for Kuramoto-Sivashinsky-type dynamics near a sideband instability \cite{DSSS09,HN86} -- see Fig. \ref{f:FingCusp-BCDE}.

\subsubsection{Reversing desertification by front instabilities}
\label{sss:FOTM}

In \cite{FOTM19} a novel mechanism is proposed through which desertification may be reversed. It is shown that an invasion front of a bare soil state into a homogeneously vegetated that is stable with respect to longitudinal perturbations may be unstable with respect to transversal perturbations as planar interface. The vegetation fingers driven by this instability grow back into the bare soil state. Thus, instead of the collapse of the savanna state into a desert state that is predicted by the 1-dimensional point of view, (realistically) allowing for 2-dimensional perturbations along the interface predicts the counter-invasion of vegetation patterns into the desert, with as endstate a fully vegetated patterned domain (see \cite{FOTM19} and Fig. \ref{f:FingIntro}).
\\ \\
The model considered in \cite{FOTM19} reads,
\begin{equation}
\label{e:RDE-FOTM}
\left\{	
\begin{array}{rcl}
U_t &=& \, \, \, \, \, \, \Delta U - U + U (1+ \mu_1 U)^2 (1- U) V \\
V_t &=& \frac{1}{\varepsilon^2} \Delta V + \mu_2 - \frac{\mu_3 V}{1 + \mu_4 U} - \mu_5 U (1+ \mu_1 U)^2 V
\end{array}
\right.
\end{equation}
(cf. \cite{ZMB15}). Note that model (\ref{e:RDE-BCDE}) can be seen as a simplification of this (somewhat) more realistic model. The evaporation rate is scaled to $-1$ in (\ref{e:RDE-BCDE}), while here it is unscaled (parameter $\mu_3$) and takes into account the effect vegetation has on evaporation ($\mu_4$). Moreover, following \cite{Kla99}, the nonlinear mechanism by which plants increase water infiltration is modeled (and scaled) as $U^2V$ in (\ref{e:RDE-BCDE}), while in (\ref{e:RDE-FOTM}) the fact that plants extend spatially through their roots is included in the infiltration term, which yields the extended term $U (1+ \mu_1 U)^2 V$ \cite{Mer18}. Apart from that the systems are only scaled in slightly different ways. As we shall see, the upcoming analysis is also quite similar -- but slightly more technical and less explicit -- to that of the simplified Klausmeier-type model (\ref{e:RDE-BCDE}). Moreover, the outcome is in essence identical.
\\ \\
The fast reduced system for traveling waves is now given by,
\beq
\label{e:fast-red-FOTM}
u_{\xi\xi} + c u_\xi - u + v_0 u (1+ \mu_1 u)^2 (1- u) =  0, \; \;
v = v_0, q = q_0
\eeq
and the associated two normally hyperbolic slow manifolds of interest are thus given by,
\beq
\label{d:Mpm0-FOTM}
\M^-_0 = \{u=f^-(v)\equiv 0, p= 0, v < 1 \}, \; \; \M^+_0 = \{u=f^+(v), p=0, v> v_{\rm m}\},
\eeq
with $u=f^+(v)$ the largest solution of
\beq
\label{d:vF-FOTM}
v_F(u) = v = \frac{1}{(1-u)(1+\mu_1 u)^2},
\eeq
that intersects the $\{u=0\}$-axis at $v = 1$ and has a minimum at,
\beq
\label{d:uvm}
(u_{\rm m}, v_{\rm m}) = \left( \frac{2 \mu_1 - 1}{3 \mu_1}, \frac{27 \mu_1}{4 (1 + \mu_1)^3} \right).
\eeq
Note that in order to have heteroclinic front solutions jumping from $\M^-_0$ to $\M^+_0$ we need to impose that $u_{\rm m} > 0$ and $v_{\rm m} < 1$, i.e. that $\mu_1 > \frac12$ -- see Fig. \ref{f:manifold&interface}(a). However, we need more: we also need that the critical points of the reduced slow flows associated on the homogeneous background states $(\bU^\pm,\bV^\pm)$ are on $\M^\pm_0$, i.e. we need that $0 < \bV^- < 1$ and $\bV^+ < v_{\rm m}$. Moreover, since these background states are assumed to be stable, we also need that the associated critical points on $\M^\pm_0$ are saddles -- see section \ref{ss:construction}. Naturally, the homogeneous background states are determined by the intersections of
\beq
\label{d:vG-FOTM}
v_G(u) = v = \frac{\mu_2(1 + \mu_4 u)}{\mu_3 + \mu_5 u (1 + \mu_4 u)(1+ \mu_1 u)^2},
\eeq
either with $u=0$ or with $v = v_F(u)$ (\ref{d:vF-FOTM}). The former yields the bare soil state $(\bU^-,\bV^-) = (\mu_2/\mu_3,0)$, and we thus impose the condition $\mu_2/\mu_3 < 1$. Since for $u > 0$ $v_G'(u)$ must have a unique zero and $\lim_{u \to \infty} v_G(u) = 0$, we conclude that $v_G(u)$ has a maximum at $u = u_{\rm M} > 0$ (Fig. \ref{f:manifold&interface}(a)). Hence, the intersection $\{v = v_F(u)\} \cap \{v = v_G(u)\}$ has at most 2 elements that both correspond to homogeneously vegetated savanna states (see also Remark \ref{r:spectralinstab-FOTM}). Although the soil moisture level below a bare soil state can both be higher or lower than that below a vegetated state \cite{KZBM14,RKBHO00}, model (\ref{e:RDE-FOTM}) has been derived from a more extended model under assumptions under which the water level in the bare soil state is necessarily higher than that in the vegetated state (cf. \cite{Mer18}). Therefore, we need to impose the condition $\bV^+ < \bV^-$, i.e. $v_G(\bU^+) < \mu_2/\mu_3$ or 
\beq
\label{e:eco-cond-FOTM}
\mu_3 \mu_4 - \mu_5 (1 + \mu_4 \bU^+)(1+ \mu_1 \bU^+)^2 < 0.
\eeq
It now follows that only the largest element of $\{v = v_F(u)\} \cap \{v = v_G(u)\}$ can correspond to a critical point on $\M^+_0$, and thus to $(\bU^+,\bV^+)$: the other one is an unstable state. Note that this situation is similar to that of the (simplified) Klausmeier model (\ref{e:RDE-BCDE}).
\\ \\
As for (\ref{e:RDE-BCDE}), the slow flow on $\M^-_0$ is linear, $v_{XX} + 1 - v = 0$, while the flow on $\M^+_0$ is given by
\beq
\label{e:slow-red-FOTM}
v_{XX} + \mu_2 - \frac{\mu_3 v}{1 + \mu_4 f^+(v)} - \mu_5 \frac{f^+(v)}{1 - f^+(v)} = 0,
\eeq
where we have used (\ref{d:vF-FOTM}) to simplify the expression. Since we assume that the background state $(\bU^+,\bV^+)$ is stable as solution of (\ref{e:RDE-FOTM}), we know that the associated critical point $(\bV^+,0)$ of (\ref{e:slow-red-FOTM}) -- its only critical point -- must be a saddle (section \ref{ss:construction}). Thus it follows -- by arguments completely similar to those in the previous section -- that the (linear) unstable manifold $\W^{\rm u}((\bV^-,0))$ on $\M^-_0$ must intersect the stable manifold $\W^{\rm s}((\bV^+,0))$ on $\M^+_0$ (in their combined projection) -- exactly as for model (\ref{e:RDE-BCDE}, see Fig. \ref{f:BCDE}(b) -- thereby defining the fast jumping point $(v_\ast,q_\ast) \in \W^{\rm u}((\bV^-,0)) \cap \W^{\rm s}((\bV^+,0))$ (with $v_{\rm m} < \bV^+ < v_\ast < \bV^- = \mu_2/\mu_3$ and $q_\ast < 0$).
\\ \\
To complete the skeleton structure, we only need to establish that there is a $\ca$ for which there is a heteroclinic connection $\ua(\xi)$ in (\ref{e:fast-red-FOTM}) between the critical points $(\ua^-,0) = (0,0)$ and $(\ua^+,0) = (f^+(\va),0)$. Although the nonlinearity of (\ref{e:fast-red-FOTM}) is quartic and we thus do not have explicit expression for $\ca$ and $\ua(\xi)$ as in the cubic case of (\ref{e:fast-red-BCDE}), it indeed follows by phase plane arguments based on the integrable nature of (\ref{e:fast-red-FOTM}) at $c = 0$ that $\ca$ and $\ua(\xi)$ exist and are uniquely determined. (In fact, an explicit condition on the value of $\va$ at the Maxwell point, i.e. the parameter combination for which the front is stationary, can still be obtained (cf. section \ref{ss:BifTrav-SmallTau}).)
\\ \\
Therefore, we again conclude by the methods of geometric singular perturbation theory that a traveling planar invasion front $(U_h(\xi),V_h(\xi))$ between the desert state $(\bU^-,\bV^-)$ and savanna state $(\bU^+,\bV^+)$ exists. Moreover, we assume that the front is stable with respect to longitudinal perturbations and thus that the critical eigenvalue curve $\la_c(\ell)$ has $\la_c(0) = 0$ -- as has become usual by now. By comparing (\ref{e:RDE-FOTM}) to (\ref{e:RDE}) we see that,
\[
F_v(u_\ast(\xi), v_\ast) u_{\ast,\xi}(\xi) = \left[u_\ast(\xi)(1 + \mu_1 \ua(\xi))^2 (1 - u_\ast(\xi))\right] u_{\ast,\xi}(\xi) > 0,
\]
(since $0 < u_\ast(\xi) < 1$ (cf. Fig. \ref{f:manifold&interface}(a))), and that,
\[
G(u^+_\ast, v_\ast) - G(u^-_\ast, v_\ast) = \left[\mu_3 \mu_4 - \mu_5 (1 + \mu_4 \ua^+)(1+ \mu_1 \ua^+)^2 \right] \frac{\ua^+ \va}{1 + \mu_4 \ua^+}.
\]
Since $\bU^+ < \ua^+$ (see again Fig. \ref{f:manifold&interface}(a)), it follows by  (\ref{e:eco-cond-FOTM}), i.e. the fact that model (\ref{e:RDE-FOTM}) has been derived under assumptions that imply that the water level below the bare soil state is higher than below the vegetated state, that $G(u^+_\ast, v_\ast) - G(u^-_\ast, v_\ast) < 0$. Hence, under the overall assumption underlying our approach that $\eps$ is sufficiently small, we may conclude by (\ref{e:cond-longtrans}) that all planar invasion fronts covered by model (\ref{e:RDE-FOTM}) are also unstable with respect transversal perturbations: an invading planar desertification front will typically be counter-invaded by a patterned, labyrinthine, vegetated state -- see the simulations in \cite{FOTM19} and Fig. \ref{f:FingIntro}.
\\ \\
Finally we note that it follows by the same arguments as above that a planar invasion front in the (a priori) un-ecological case that $\bV^+ > \bV^-$, i.e. a setting for which (\ref{e:RDE-FOTM}) cannot be derived, is not unstable with respect to long wavelength transversal perturbations, in the sense that $\la_{2,c} < 0$ (\ref{e:la2c-expl}), (\ref{e:cond-longtrans}), since $\ua^+ < \bU^+$ in that case. Thus, it could be stable (although one first needs to analyze the spectral stability of the interface for all $\ell \in \RR$ to arrive at that conclusion).

\begin{remark}
\label{r:spectralinstab-FOTM}
\rm
In \cite{FOTM19}, the existence and stability of `weak fronts' is studied by an amplitude equation approach near the co-dimension 2 point at which all 3 background states collide (i.e. asymptotically close to $(\mu_{1,c},\mu_{2,c}) = ((\mu_3 + \mu_5)/(2 \mu_3), \mu_3)$ for which $(\bU,\bV) = (1,0)$ is a triple zero of $F(U,V) = G(U,V) = 0$ -- if $\mu_4 = 0$, which is the case considered in \cite{FOTM19}). By deriving an expression for the curvature driven normal velocity of the interface, an explicit leading order upperbound $\eps_c$ is obtained below which the front is unstable with respect to transversal perturbations -- thus, diffusion coefficient $\eps$ is considered to be an $\O(1)$ parameter in this analysis. Since the assumption that $\eps$ is sufficiently small is underlying all our (asymptotic) analysis, the result of \cite{FOTM19} agrees with our analysis.  On the other hand, in \cite{FOTM19} a simulation is shown of a `linearly stable desertification front' that is `unstable to finite-amplitude transverse modulations' (Fig. 4 in \cite{FOTM19} in which $\vmu = (3.5, 1.15, 3.2, 1.0, 0.5) \in \RR^5$). This precise parameter combination is at the boundary of the applicability of our existence analysis, nevertheless, our stability analysis predicts that this interface must be spectrally (and thus linearly) unstable with respect to transversal long wavelength perturbations -- if $\eps$ is sufficiently small. A careful simulation that is allowed to run for a sufficiently long time with $\eps = 1/\sqrt{150} \approx 0.08$ as in \cite{FOTM19} shows that this value of $\eps$ indeed is sufficiently small: also this front is unstable, i.e. it develops fingers (like all other simulations shown in \cite{FOTM19}).
\end{remark}

\begin{remark}
\label{r:JDCBM}
\rm
The model considered in \cite{JDCBM20} can be seen as `intermediate' between models (\ref{e:RDE-BCDE}) and (\ref{e:RDE-FOTM}) -- like (\ref{e:RDE-FOTM}), it is also based on \cite{ZMB15}. Application of our methods to this model goes along the very same lines as here and yields completely similar results.
\end{remark}

\subsection{A FitzHugh-Nagumo model}
\label{ss:FHN}

\begin{figure}[t]
\hspace{.02\textwidth}
\begin{subfigure}{.45 \textwidth}
\centering
\includegraphics[width=1\linewidth]{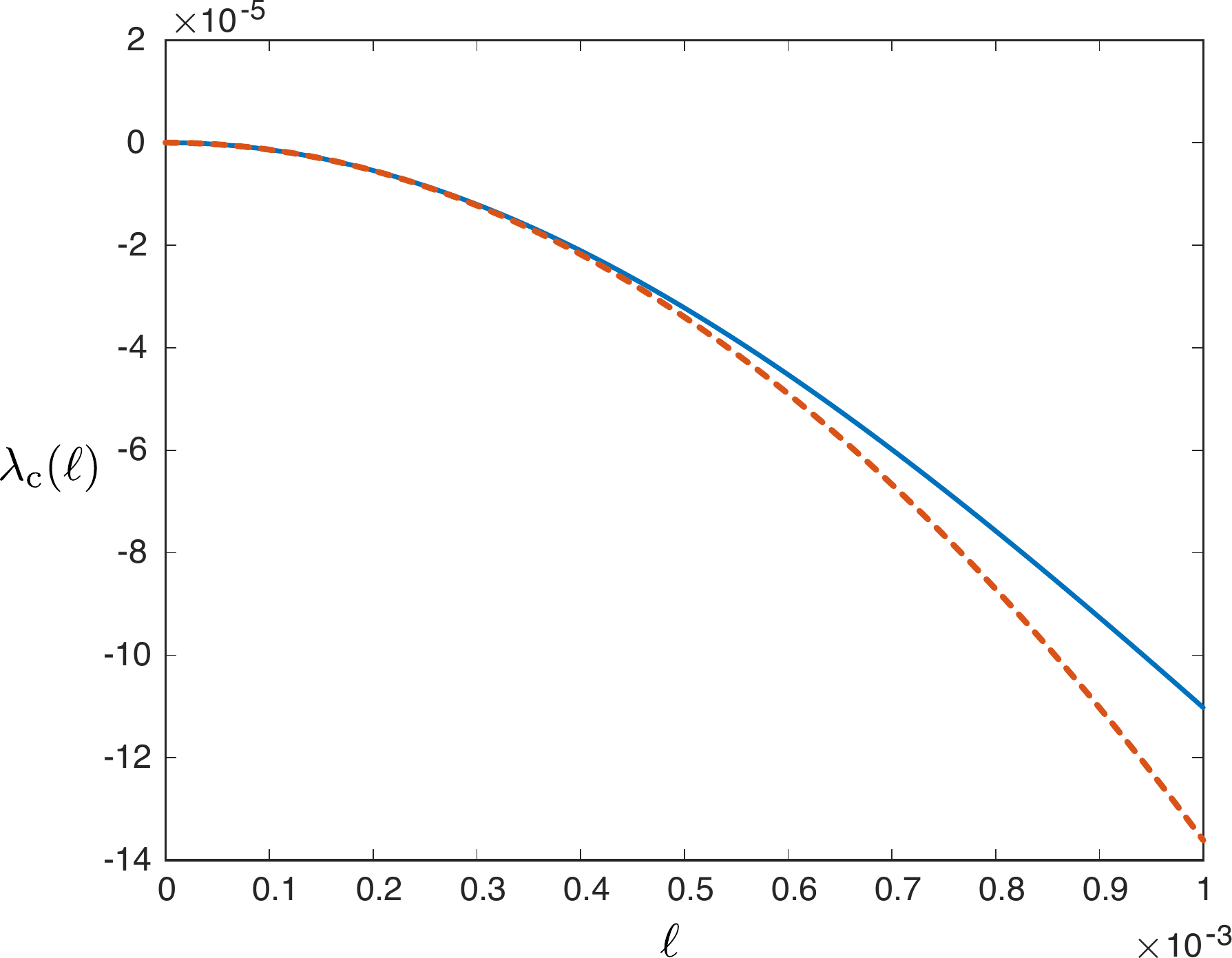}
\caption{}
\end{subfigure}
\hspace{.05\textwidth}
\begin{subfigure}{.45 \textwidth}
\centering
\includegraphics[width=1\linewidth]{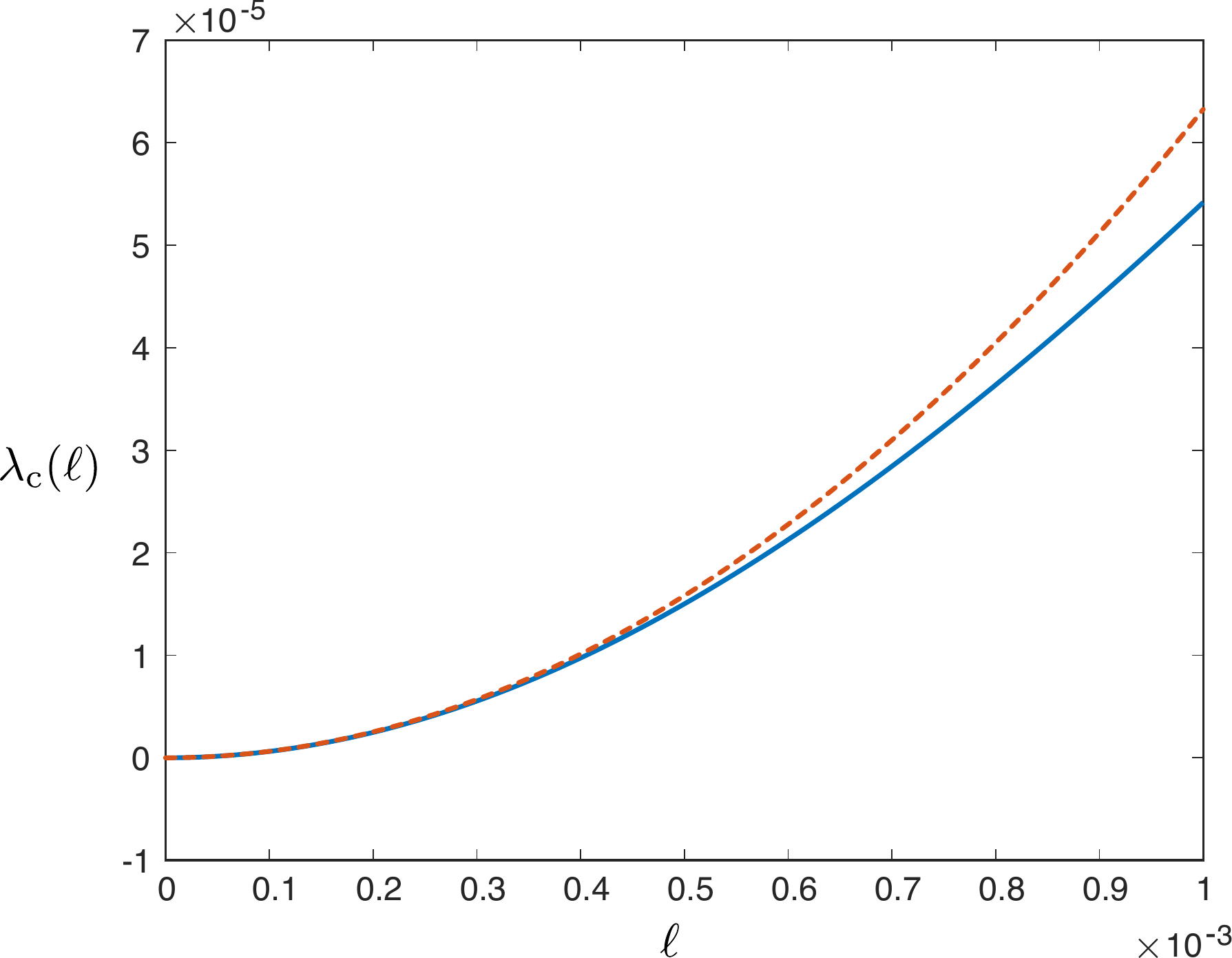}
\caption{}
\end{subfigure}
\caption{\small{Comparisons between the numerical evaluations (solid blue) and the parabolic approximation (dashed red) by (\ref{d:la2bu2bv2}) with (\ref{e:la2c-expl-FHN}) of the critical curves $\la_c(\ell)$ associated to the stationary front of (\ref{e:RDE-FHN}) with $f^c(v) \equiv 0$, $\tau=1$, $\mu_1=4$, $\eps = 0.001$, and $f^+(v) = 1+\mu_2 v$. (a) $(f^+)'(0) = \mu_2 = 1 > 0$ so that $\la_{2,c} = -\frac{1}{3 \eps \sqrt{6}}< 0$ and (b) $\mu_2 = -1 < 0$ with $\la_{2,c}= \frac{1}{5\eps \sqrt{10}} > 0$. (Note that the remaining integrals in (\ref{e:la2c-expl-FHN}) can also be evaluated explicitly in this case -- see section \ref{sss:InStabBifTrav}.)}}
\label{f:FHN-front}
\end{figure}

We consider a modified version of the FitzHugh-Nagumo models studied in \cite{HM94,HM97a,HM97b},
\begin{equation}
\label{e:RDE-FHN}
\left\{	
\begin{array}{rcccl}
\tau U_t &=& \Delta U  & - & (U-f^-(V))(U-f^c(V))(U-f^+(V))\\
V_t &=& \frac{1}{\varepsilon^2}\Delta V & + & U - \mu_1 V
\end{array}
\right.
\end{equation}
(for $\mu_1 \geq 0$). In this section we consider $\tau = \O(1)$, in section \ref{sss:InStabBifTrav} we explicitly study  the impact of the pre-factor $\tau$ in (\ref{e:RDE}) in the context of specific choices of the functions $f^{\pm, c}(v)$. To mimic the systems studied in \cite{HM94,HM97a,HM97b}, we restrict ourselves to the symmetric case $f^-(-v) = -f^+(v)$ and we note that (\ref{e:RDE-FHN}) corresponds to the systems in \cite{HM94,HM97a,HM97b} with $u = f^{\pm, c}(v)$ the 3 solutions of $u-u^3 = v$. For $v_0$ such that $f^-(v_0) = -f^+(-v_0) < f^c(v_0) < f^+(v_0)$, the critical points $(f^\pm(v_0), 0)$ are saddles of the fast reduced systems (\ref{e:fast-red}), so that the (normally hyperbolic) slow manifolds $\M^\pm_0$ associated to the travelling wave problem (\ref{e:DS}) indeed are given by (\ref{d:Mpm0}). On $\M^\pm_0$, the reduced slow flows are, for $\tau = \O(1)$, given by
\beq
\label{e:slowred-FHN}
v_{XX} + f^\pm(v) - \mu_1 v =0.
\eeq
Note that these flows are by construction symmetrical under $v \to -v$, $X \to -X$. It thus follows (by the symmetry) that if the slow reduced flow (\ref{e:slowred-FHN}) on $\M^+_0$ has a saddle $(\bV^+,0)$ with $\bV^+$ such that $-f^+(-\bV^+) < f^c(\bV^+) < f^+(\bV^+)$ and such that its stable manifold $\W^{\rm s}((\bV^+,0))$ intersects the $v_X$-($= q$)-axis (on $\M^+_0$), then the unstable manifold $\W^{\rm u}((\bV^-,0))$ of the slow flow on $\M^-_0$ necessarily intersects the $q$-axis at the same point. Thus, by the general approach of section \ref{ss:construction}, there exists a traveling slow-fast-slow front $(U_h(x-c_\ast t), V_h(x-c_\ast t))$ in (\ref{e:RDE-FHN}) that connects the stable background states $(\bU^-,\bV^-) = (f^-(\bV^-),\bV^-) = (-f^+(-\bV^+),-\bV^+)$ and $(\bU^+,\bV^+) = (f^+(\bV^+),\bV^+)$ (where we note that necessarily $\mu_1 \geq 0$ (\ref{e:cond-stabbUVpm})). Due to the cubic character (in $U$) of $F(U,V)$ in (\ref{e:RDE-FHN}), we have explicit expressions both for $c_\ast$ and $u_\ast(\xi) = u_\ast(x - c_\ast t)$,
\beq
\label{e:FHNcuast}
c_\ast = \frac{\sqrt{2} f^c(0)}{\tau}, \; u_\ast(\xi) = f^+(0) \tanh \frac12 \sqrt{2} f^+(0) \, \xi
\eeq
(\ref{e:Kcast-gen}), (\ref{e:uast-gen}). In the special (stationary) case that $f^c(0) = 0$, that includes the models of \cite{HM94,HM97a,HM97b} we obtain in a straightforward manner,
\beq
\label{e:FGIf-c=0-FHN}
F_\ast =  \frac43 (f_+(0))^3 \left( (f^+)'(0) - (f^c)'(0) \right), \;
G_\ast = 2 f_+(0), \;
\int_\RR (u_{\ast, \xi})^2 d \xi = \frac23 \sqrt{2} (f_+(0))^3
\eeq
(\ref{d:FGast}), which yields by (\ref{e:la2c-expl}),
\beq
\label{e:la2c-expl-FHN}
\la_{2,c} = - \frac{(f^+)'(0) - (f^c)'(0)}{\tau \eps \sqrt{2} f^+(0)} \left(\int_{-\infty}^0 (v^-_X)^2 \,dX + \int_0^{\infty} (v^+_X)^2 \,dX \right).
\eeq
Under the assumption that the front is longitudinally stable -- which can be checked analytically -- we conclude that the interface will be destabilized by transversal perturbations if $(f^+)'(0) < (f^c)'(0)$ -- so that for these $f^{+,c}(v)$ fingering may typically be expected (for $\tau = \O(1)$). In Fig. \ref{f:FHN-front} we validate our asymptotic analysis by two numerical evaluations of the critical curve $\la_c(\ell)$ for FitzHugh-Nagumo(-type) model (\ref{e:RDE-FHN}): one with $\la_{2,c}<0$ and the other with $\la_{2,c}>0$. We conclude from this figure that our analysis indeed has the expected accuracy. We refer to \cite{BCDL22} for similar numerical validations -- especially for radially symmetric (multi-)front patterns -- in the context of ecosystem model (\ref{e:RDE-BCDE}).
\\ \\
The models considered in \cite{HM94,HM97a,HM97b} have $u-u^3 = v$, so that
$f'(v) = 1/(1-3 f^2(v))$. Since in this case $f^c(0)=0$ and $f^+(0)=1$ we conclude that $(f^c)'(0) = 1 > (f^+)'(0) = -\frac12$: for $\tau = O(1)$ the fronts are unstable with respect to transversal perturbations and in fact finger -- as can be seen in numerical simulations. 
\\ \\
Lastly, we note that in the (general) case $f^c(0) \neq 0$, taking $f^c(0) \in (-f^+(0),f^+(0))$ as parameter yields by a slightly more involved calculus that $\la_{2,c}(f^c(0))$ changes sign at a certain value of $f^c(0) \in (-f^+(0),f^+(0))$ if $(f^c)'(0)/(f^+)'(0) \in (-2,1)$.

\subsection{Example systems with $\tau = \eps \ttau$}
\label{ss:ExampleSystems-t}

\subsubsection{On cylindrical domains}
\label{sss:Cylindrical}

The stability of planar interfaces in a general class of singularly perturbed 2-component reaction-diffusion systems has
also been considered in \cite{Tan03,TN94}. The class of models considered in \cite{Tan03,TN94} is more restricted than the models considered here. Several simplifying assumptions are formulated on the functions $F(U,V)$ and $G(U,V)$, it is for instance assumed that the slow flows on $\M^\pm_0$ do not have additional critical points between the saddles at $(\bV^\pm,0)$ and the intermediate `jumping line' $\{v = \va\}$. Especially in the setting of ecological models this indeed is a restriction \cite{Doe22,JDCBM20}. Moreover, ecological models also typically have $\tau = \O(1)$ -- in fact $\tau = 1$ -- while in \cite{Tan03,TN94} only the case $\tau = \eps \ttau$ is considered. However, we have seen that $\tau = \O(1)$ is more simple than $\tau = \O(\eps)$ and the nature and strength of the methods of \cite{Tan03,TN94} -- that are based on \cite{NF87,NMIF90} -- strongly suggest that the more analytical approach of these works (compared to the present more geometrical set-up) can also be successfully applied beyond the class of models considered in \cite{Tan03,TN94}. In fact, we expect that any model that can be studied along the lines of this paper can also be analyzed by the methods of \cite{NF87,NMIF90,Tan03,TN94}.
\\ \\
A much more important difference between our present study and \cite{Tan03,TN94} is that in the latter works the spatial domain is considered to be cylindrical: instead of $(x,y) \in \RR^2$ considered here, $(x,y) \in \RR \times \Omega_y$ with $\Omega_y = (0,L)$ bounded (in fact, in \cite{TN94}, $(x,y) \in \Omega_x \times \Omega_y$ with also $\Omega_x$ bounded (as in \cite{NF87}), while $\Omega_y \subset \RR^{N-1}$ in \cite{Tan03} -- these differences can be overcome though). The (in)stability criteria obtained here are against long wavelength transversal perturbations: if an interface is unstable in this sense, it will only be unstable on a cylindrical domain $\RR \times (0,L)$ for $L$ sufficiently large. So, in a way \cite{Tan03,TN94} consider the more complex issue of the destabilization of the interface for a given size $L$ of the domain $\Omega_y$ -- in fact, a large part of the analytical focus of \cite{Tan03,TN94} is on the nature of the fastest growing unstable `eigenmode', a topic we do not consider here.
\\ \\
A clear advantage of the more simple setting of the present work is that it enables us to also obtain simple instability results for the models considered in \cite{Tan03,TN94} under the assumption that $L$ is sufficiently large. Instead of considering the models of \cite{Tan03,TN94} in their full generality, we focus on an example model also considered in \cite{NF87,Tan03,TN94},
\begin{equation}
\label{e:RDE-Cylindrical}
\left\{	
\begin{array}{rcl}
\tau U_t &=& \, \, \, \, \, \, \Delta U + U^2 (1- U) - U V  \\
V_t &=& \frac{1}{\varepsilon^2} \Delta V + \mu_1 UV + \mu_2 V - \mu_3 V^2
\end{array}
\right.
\end{equation}
Clearly, the normally hyperbolic slow manifolds $\M^\pm_0$ are given by,
\[
\M^-_0 = \{ u = p = 0, v > 0 \}, \; \; \M^+_0 = \{ u = u^+(v) = \frac12 + \frac12 \sqrt{1 - 4 v}, p = 0, v < \frac14\},
\]
and the reduced slow flows on $\M^\pm_0$ by,
\[
\M^-_0: v_{XX} + \mu_2 v - \mu_3 v^2 = 0, \; \; \M^+_0: v_{XX} + \mu_1 v u^+(v) + \mu_2 v - \mu_3 v^2 = 0.
\]
Both flows have $(0,0)$ as critical point and we want this to be a center in both cases: $\mu_2 > 0$ and $\mu_1 + \mu_2 > 0$. We also want the slow flows to both have exactly one additional critical point on $\M^\pm_0$, i.e. $\mu_3 > {\rm max}\{0,2(\mu_1+2\mu_2)\}$, which then must be saddles: $(\bV^\pm ,0)$ with $\bV^\pm > 0$. Since one of these saddles always must be inside the homiclinic loop attached to the other it follows that $\W^{\rm u}((\bV^- ,0))$ must intersect $\W^{\rm s}((\bV^+ ,0))$ at $(\va,\qa)$ with $\va > {\rm min}\{\bV^\pm\} > 0$. Since $\va < {\rm max}\{\bV^\pm\} < \frac14$ for $\mu_1 > 0$ -- the case considered \cite{Tan03,TN94} -- the existence of the traveling front follows immediately; for $\mu_1 < 0$, there may be 3 critical points on $\M^-_0$ in which case a fast heteroclinic jump may not always be possible so that  additional (technical) conditions need to be added to ensure the existence of a traveling front. Since by construction $\ua(\xi) > 0$, $u_{\ast, \xi}(\xi) > 0$ and,
\beq
\label{e:tFatGa-Cylindrical}
F_v(\ua,\va) = - \ua, \; \; G(\ua^+,\va) - G(\ua^-,\va) = \mu_1 \ua^+,
\eeq
it follows immediately by (\ref{e:cond-longtrans}) that the planar traveling interface must be unstable in the case $\mu_1 > 0$ considered in \cite{Tan03,TN94}: even if it is stable with respect to longitudinal perturbations, it is unstable against transversal long wavelength perturbations.
\\ \\
However, we so far only considered $\tau = \O(1)$, while $\tau = \eps \ttau = \O(\eps)$ in \cite{Tan03,TN94}. Nevertheless, this does impact neither the instability result itself nor its simplicity (for $\mu_1 >0$). As in \cite{TN94} and section \ref{ss:BifTrav-SmallTau}, we focus on the case of a stationary planar interface, i.e. we assume that $\vmu \in \RR^3$ (in (\ref{e:RDE-Cylindrical})) is such that $\tHH_f(\tua^-,0) = \tHH_f(\tua^+,0) = 0$ (cf. (\ref{d:Ham-fast}) and \cite{NF87,TN94}). Since the fast field associated to (\ref{e:RDE-Cylindrical}) is cubic, we note
\[
\tva = \frac29, \; \tua^+ = \frac23, \; \tF_{\ast} = - \frac29, \; \tG_\ast = \frac{4}{27}, \; \tI_f = \frac{2}{81} \sqrt{2}
\]
((\ref{d:tFGast}), (\ref{d:tIsf}), (\ref{d:fastcubic-gen-1st}), (\ref{e:Kcast-gen}), (\ref{e:tFatGa-Cylindrical})), which implies that the planar interface is unstable with respect to transversal perturbations for,
\[
\ttau > \ttau_\ast = \frac{243}{8} \tI_s(\vmu),
\]
at which it undergoes a bifurcation into traveling waves (by section \ref{ss:BifTrav-SmallTau}). Thus, also for $\tau = \O(\eps)$, the straightforward observations (\ref{e:tFatGa-Cylindrical}) immediately yield that a longitudinally stable stationary front on $\RR^1$ must be unstable as planar interface on $\RR^2$ up to the critical value $\ttau_\ast$ at which it (also) destabilizes as 1-dimensional front. In the setting of \cite{Tan03,TN94}, this implies that the interface is always unstable on cylindrical domains $\RR \times (0,L)$ for $L$ sufficiently large. For a given $L$, the instability condition is much more involved and needs to be evaluated numerically, see \cite{TN94} (where we note that in \cite{TN94} also the $x$-domain is bounded, which may impact the nature of the front).
\\ \\
Finally, we note that for $\mu_1 < 0$, the stationary interface -- if it exists and is longitudinally stable -- will always be stable against long wavelength transversal perturbations -- again by simple observation (\ref{e:tFatGa-Cylindrical}). In fact, since neither $\tF_\ast$, $\tG_\ast$, $\tI_f$ nor $\tI_s$ depends on $\ttau$, it follows that $\la_{2,c} \to -1/\eps^2$ as $\ttau \ll 1$ (\ref{e:la2c-expl-t}).

\subsubsection{The (in)stability of bifurcating traveling planar interfaces}
\label{sss:InStabBifTrav}

Finally, we consider the case $\tau = \eps \ttau$ in the setting of FitzHugh-Nagumo model (\ref{e:RDE-FHN}). Naturally, bifurcations may take place due to the geometry of the now asymmetrical slow reduced flows (for $\tc \neq 0$, cf. Fig. \ref{f:slowflows-t}(b)). However, here we focus on the bifurcation into traveling waves as discussed in section \ref{ss:BifTrav-SmallTau}. In this section we consider the question of the (in)stability of the traveling interfaces that bifurcate off the stationary front with respect to transversal perturbations. Note that it follows from the analysis in section \ref{ss:BifTrav-SmallTau} that $\tla_{2,c}$ associated to the primary stationary interface has changed sign in the bifurcation (\ref{e:la2c-expl-t}), however, the existence/Melnikov analysis of section \ref{ss:BifTrav-SmallTau} does not provide insight in the transversal stability of the bifurcating waves, and more specifically not in the associated sign of $\tla_{2,c}$.
\\ \\
As in \ref{ss:BifTrav-SmallTau}, this issue can be studied in the general setting of equation (\ref{e:RDE}). Especially since this is a somewhat technical affair, we for simplicity restrict ourselves to a simple case for which explicit analysis is possible: we consider the FitzHugn-Nagumo model (\ref{e:RDE-FHN}) in the most simple case in which the slow reduced flows on $\M^\pm_0$ are linear (cf. the approximations in \cite{HM94}). Thus, we set $f^+(v) = 1 + \mu_2 v$ and thus $f^-(v) = -1 + \mu_2 v$, so that
\beq
\label{e:slowflow-FHN-t-lin}
v_{XX} + \tc v_X \pm 1 - (\mu_1 -\mu_2)v =0
\eeq
on $\M^\pm_0$ (cf. (\ref{e:slow-red-t}), (\ref{e:slowred-FHN})). Naturally, we need to assume that $\mu_1 -\mu_2 > 0$ (so that the critical point of (\ref{e:slowflow-FHN-t-lin}) is a saddle). It follows that the slow components of the skeleton structure that determines the front are given by,
\beq
\label{e:tvpm-FHN}
\tv^\pm(X) = \left(\frac{\tc}{(\mu_1 -\mu_2)\sqrt{\tc^2 + 4(\mu_1 -\mu_2)}} \mp \frac{1}{\mu_1 -\mu_2} \right) e^{-\frac12\left(\tc \pm \sqrt{\tc^2 + 4(\mu_1 -\mu_2)}\right)X} \pm \frac{1}{\mu_1 -\mu_2},
\eeq
and thus that,
\beq
\label{e:tvast-FHN}
\tva = \tv^\pm(0) = \frac{\tc}{(\mu_1 -\mu_2)\sqrt{\tc^2 + 4(\mu_1 -\mu_2)}}
\eeq
(at leading order in $\eps$). By the cubic nature of $F(U,V)$, the relation $\C(v_0)$ that determines for which value of `friction term' $\tc \ttau$ there is a heteroclinic orbit in (\ref{e:fast-red-t}) can again be determined explicitly,
\[
\tc \ttau = \C(v_0) = \sqrt{2} \left(f^c(v_0) - f^+(v_0) \right) = \sqrt{2} \left(f^c(v_0) - \mu_2 v_0 \right),
\]
(\ref{e:Kcast-gen}). By (\ref{e:tvast-FHN}), (\ref{e:cast-t}) thus takes the explicit form,
\[
\frac12 \sqrt{2} \tc \ttau =  f^c\left(\frac{\tc}{(\mu_1 -\mu_2)\sqrt{\tc^2 + 4(\mu_1 -\mu_2)}}\right) - \left(\frac{\mu_2 \tc}{(\mu_1 -\mu_2)\sqrt{\tc^2 + 4(\mu_1 -\mu_2)}}\right).
\]
The further simplifying assumption $f^c(v) = \mu_3 v$ yields,
\beq
\label{e:tc-FHN}
\tc = 0 \; \; {\rm or} \; \; \tc^2 = \frac{2 (\mu_3 - \mu_2)^2}{(\mu_1 -\mu_2)^2 \ttau^2} - 4 (\mu_1 -\mu_2),
\eeq
with the additional condition $\mu_3 - \mu_2 > 0$ (since $\ttau > 0$). Thus, the situation is as in section \ref{ss:BifTrav-SmallTau}: independent of the value of $\ttau$ there is a stationary front that undergoes a bifurcation into (counter-propagating) travelling waves/planar interfaces as $\ttau$ decreases through
\beq
\label{d:ttauast}
\ttau_{\ast} = \frac{\mu_3 - \mu_2}{(\mu_1 -\mu_2) \sqrt{2 (\mu_1 -\mu_2)}}.
\eeq
In fact, we have already obtained more than in section \ref{ss:BifTrav-SmallTau}: we know that the bifurcation is supercritical (i.e. the bifurcating waves appear for values of $\ttau$ at which the stationary wave is unstable). Note however that there is also a difference with the general setting of section \ref{ss:BifTrav-SmallTau}: due to the assumed symmetry $f^-(-v) = - f^+(v)$ the stationary interface exists in (\ref{e:RDE-FHN}) without an additional condition on the parameters $\vmu$ ($\in \RR^3$ in (\ref{e:RDE-FHN}) with the present choices of $f^+(v)$ and $f^c(v)$). Note also that we can copy (\ref{e:FGIf-c=0-FHN}) for the stationary fronts and that by (\ref{e:tvpm-FHN}),
\beq
\label{e:tIs-FHN}
\tI_s(\tc) = \frac{8}{\left(\tc^2 + 4(\mu_1 -\mu_2)\right)^{3/2}},
\eeq
so that by (\ref{e:la2c-expl-t}),
\beq
\label{e:la2c-FHN-t}
\tla_{2,c}(\ttau) = \frac{1}{\eps^2} \, \frac{\mu_3 - \mu_2}{(\mu_1 -\mu_2) \ttau \sqrt{2 (\mu_1 -\mu_2)}-(\mu_3 - \mu_2)}.
\eeq
Thus, assuming that the stationary front remains stable against longitudinal perturbations -- which can be established analytically -- we conclude that the associated planar interface is unstable against (long wavelength) transversal perturbations for $\ttau > \ttau_{\ast}$, where $\ttau_{\ast}$ of (\ref{d:ttauast}) indeed coincides with a zero in the denominator of (\ref{e:la2c-FHN-t}) -- as was already established by the general analysis of section \ref{ss:BifTrav-SmallTau}.
\\ \\
Under the assumption that the bifurcating travelling interfaces are stable as solutions of $x \in \RR$, which is natural by the above obtained supercritical nature of the bifurcation (and can be established analytically), we can now consider the question of the stability of the interfaces with respect to long wavelength transversal perturbations. Since we know from section \ref{ss:ExStab-SmallTau} that this is determined by (\ref{e:la2c-expl-t}), we can do this by evaluating $\tF_\ast$, $\tG_\ast$ (\ref{d:tFGast}), $\tI_f$ and $\tI_s$ (\ref{d:tIsf}) for $\tc \neq 0$.
\\ \\
As in section \ref{ss:BifTrav-SmallTau}, we consider the system near the bifurcation: we introduce $\htau > 0$ and second small parameter $\delta$, i.e. $0 < \eps \ll \delta \ll 1$ and independent of $\eps$ -- by setting $\ttau = \ttau_\ast - \de^2 \htau$, so that by (\ref{e:tvast-FHN}) and (\ref{e:tc-FHN}),
\beq
\label{e:hvhc-FHN}
\tc = \tc_\ast = \delta \hc(\htau) = \pm \delta \frac{4 (\mu_1-\mu_2)^{3/2} \sqrt{\htau} \left(1 + \O(\de^2)\right)}{\sqrt{(\mu_3-\mu_2)\sqrt{2(\mu_1-\mu_2)}}}, \;
\tva(\tc) = \delta \hv_1(\htau) = \pm \delta \frac{2 \sqrt{\htau} \left(1 + \O(\de^2)\right)}{\sqrt{(\mu_3-\mu_2)\sqrt{2(\mu_1-\mu_2)}}},
\eeq
with $\tc$, $\hc$ and $\hv_1$ as defined in section \ref{ss:BifTrav-SmallTau}. It follows by the derivation in appendix \ref{a:tMa} that 
\beq
\label{e:denominator-FHN-t}
\tM_\ast(\htau) = \left( \tF_\ast \tI_s + \ttau \tG_\ast \tI_f = \right)  \frac83 \sqrt{2} \htau \de^2 + \O(\de^3)
\eeq
for $\delta$ sufficiently small. By (\ref{e:la2c-expl-t}) and (\ref{e:tIftFa-cneq0}) we now have for the bifurcating traveling waves,
\beq
\label{e:la2c-expl-FHN-t}
\la_{2,c}(\htau) = \frac{\mu_3 - \mu_2}{2\eps^2 \de^2 \htau (\mu_1 - \mu_2)\sqrt{2(\mu_1 - \mu_2)}}\left(1 + \O(\de)\right)
\eeq
(at leading order in $\eps$) and conclude that (in this specific case) the counter-propagating planar interfaces that appear from the bifurcation into traveling waves are unstable with respect to long wavelength transversal perturbations -- like the longitudinally stable stationary interface they originate from (\ref{e:la2c-FHN-t}).
\\ \\
Note that although the stationary front has become longitudinally unstable by the bifurcation into traveling waves, one can still apply (\ref{e:la2c-FHN-t}) to conclude that the curve $\la_{t}(\ell)$ through the translational eigenvalue -- which thus no longer is the critical curve $\la_{c}(\ell)$ -- has a local maximum at $\la_t=0$, i.e. that the orientation of $\la_t(\ell)$ flips as $\ttau$ decreases through $\ttau_\ast$. However, it follows from (\ref{e:la2c-expl-FHN-t}) -- and the analysis in the appendix -- that the orientation of the critical curve $\la_{c}(\ell)$ does not change by the transition from the stationary interface to the bifurcating traveling interfaces. Especially the fact that the sign of (\ref{e:denominator-FHN-t}) does not depend on the parameters in the problem suggests that there may be an underlying mechanism, i.e that there is a direct relation between the sign of $\la_{2,c}$ for the (longitudinally stable) standing interface before the bifurcation and that of the bifurcating traveling interfaces. Such a general relation may be uncovered by a similar analysis as in appendix \ref{a:tMa} in the general setting of section \ref{ss:BifTrav-SmallTau}. Here, we refrain from going further into this -- see also the upcoming discussion. 

\section{Discussion}
\label{s:Disc}

In \cite{Basetal18}, predictions of model studies of vegetation patterns in dryland ecosystems have been validated by observational evidence. These model studies go back to \cite{SDERRS15,Sitetal14} in which the persistence of Turing patterns in a (generalized Klausmeier) model of type (\ref{e:RDE}) under slowly varying parameters was investigated. The main conclusion of \cite{Basetal18} was that the multistability induced by the richness of spatial (Turing) patterns increases the resilience of an ecosystem: instead of a catastrophic collapse, a patterned state may adapt under worsening external circumstances to a nearby adjacent (vegetated and patterned) state. It is postulated in \cite{RBBKBD21} that this is not special for dryland ecosystems: any ecosystem that exhibits spatial patterns is expected to have an increased resilience. In fact, it is argued in \cite{RBBKBD21} that spatial patterning may provide the ecosystem with a mechanism to evade tipping, and thus to circumvent catastrophic collapse.
\\ \\
Although the patterns of \cite{Basetal18,SDERRS15,Sitetal14} that initiated the insights of \cite{RBBKBD21} are all of Turing pattern type, it is observed in \cite{RBBKBD21} that there are other types of patterns -- with backgrounds that are not of Turing type -- that will have the same positive impact on the resilience of an ecosystem. An important role is given to the patterns that originate from the interfaces between homogeneous states. These coexistence states occur naturally in various kinds of ecosystems, moreover, these ecosystems are typically modeled by systems of reaction-diffusion equations like (\ref{e:RDE}) -- see \cite{BR13,Eig21,ES20,FOTM19,vLetal03,ZMB15} and the references therein. Thus, the conditions obtained here on the instability of planar interfaces and the associated onset of fingering patterns -- and especially the insight of section \ref{ss:fingering} on the natural occurrence of this instability in ecosystem models -- have, by \cite{RBBKBD21}, a direct interpretation in the context of the resilience of ecosystems.
\\ \\
Nevertheless, the present results mostly only are a first step towards understanding the relevance of transversal instabilities in ecosystems. 
\\ \\
$\bullet$ {\bf Non-planar fronts.} The present focus on straight, flat interfaces orthogonal to the longitudinal direction of a 1-dimensional front is quite restricted -- especially from the ecological point of view. In the companion paper \cite{BCDL22}, radially symmetric curved fronts are considered -- in the explicit context of ecosystem model (\ref{e:RDE-BCDE}). Using geometric blow-up methods, first a spot pattern -- i.e. a localized (circular) vegetated area surrounded by bare soil -- is constructed rigorously. Gap patterns can be constructed in a very similar fashion and we find that parameter condition (\ref{e:condvast-BCDE}) on the existence of stationary planar fronts determines the transition between the existence of spots and gaps (in parameter space). It is additionally indicated in \cite{BCDL22} how the same geometric approach yields the existence of ring and target patterns. Moreover, the stability of the gap and spot patterns is considered. It is shown that  spots and gaps with a sufficiently large radius are unstable with respect to the same (fingering) instability as encountered here: simulations show that the evolving fingers merge the vegetated and bare soil patches into a fully mixed labyrinthine state.  
\\ \\
$\bullet$ {\bf Multiple fronts and non-flat terrains.} Existence conditions {\bf (E-I)} - {\bf (E-IV)} allow for the construction of double-front homoclinic stripe patterns -- either of gap or of spot type. Moreover, these localized stripes are the `endpoints' of a 1-parameter family of periodic stripe patterns -- see \cite{JDCBM20} for the construction of these patterns in a model of intermediate complexity between (\ref{e:RDE-BCDE}) and (\ref{e:RDE-FOTM}) (Remark \ref{r:JDCBM}). Stripe patterns are typically observed on sloped terrains \cite{DCLBB11,GWIGS18}, therefore it is natural to study these patterns in an extended version of model (\ref{e:RDE}) that includes an advection term `$+\mu_{m+1}V_x$' in the $V$-equation -- which models the downslope flow of water in an ecosystem setting \cite{BCD19,Kla99,SD17,Sher10}. Central questions to consider are: will the multiple-front patterns (also) be unstable with respect to transversal long wavelength perturbations, and what is the impact of parameter $\mu_{m+1}= \mu_{m+1}(\eps)$? (This is the subject of work in progress.)
\\ \\
$\bullet$ {\bf Multi-component systems.} Ecosystem models typically contain more than one type of vegetation and often also separate equations for soil water and overland water -- see for instance \cite{BR13,Eig21,KZBM14,Mer18,RBBKBD21} and the references therein. The 2-component ecosystem models of \cite{BCD19,Basetal18,ES20,FOTM19,JDCBM20,Kla99,vLetal03,ZMB15} are either obtained by reduction from a more extended model or postulated as conceptual, simplified, model. To achieve a next level in ecological validity, the present analysis needs to be extended to multi-component systems. The observation that many ecosystem models have various ordered small parameters opens up the possibility for setting up a geometric multi-scale analysis -- as initial explorations indicate. 
\\ \\
The present analysis also triggers a number of mathematical research questions, of which we mention a few. 
\\ \\
$\bullet$ {\bf Bifurcations of fronts and interfaces.} In section \ref{ss:stability}, stability assumption {\bf (S-I)} was formulated more precisely: apart from assuming that the interface is stable against longitudinal perturbations, it was also imposed that the dimension of the kernel of operator $\LL$ (\ref{d:LL}) is equal to 1. In other words, it is assumed that the parameters $(\tau,\vmu)$ of (\ref{e:RDE}) are chosen away from potential (local) bifurcation sets. This assumption underlies the approach of this section -- and the paper -- since the leading order approximations (\ref{e:la2c-expl}) and (\ref{e:la2c-expl-t}) -- that underlie criteria (\ref{e:cond-longtrans}) and (\ref{e:cond-longtrans-t-Intro}) -- can now be derived by approximating the 1-dimensional kernel of the adjoint operator $\LL^A$ (\ref{d:LLA}). As $(\tau,\vmu)$ or $(\ttau,\vmu)$ passes through a bifurcational (co-dimension 1) set $(\tau_\ast,\vmu_\ast)$/$(\ttau_\ast,\vmu_\ast)$ -- for instance as in section \ref{ss:BifTrav-SmallTau} -- then there is a second eigenvalue curve $\la_b(\ell)$ that passes through $\la_c(\ell)$. Naturally, $\la_b(0) = \la_c(0) = 0$: the dimension of the kernel of $\LL$ is (at least) 2 at the bifurcation set. Thus, one has to extend the approach that led to (\ref{e:la2c-expl}) and (\ref{e:la2c-expl-t}) to determine the local character -- as function of $|\ell| \ll 1$ -- of $\la_b(\ell)$ and $\la_c(\ell)$ as the bifurcational set is crossed. This is an interesting line of research, especially since the results of section \ref{sss:InStabBifTrav} suggest that there may be a relation between the signs of $\la_{2,c}$ of the longitudinally stable interfaces before and after the bifurcation. \\ \\
$\bullet$ {\bf Breakdown of existence conditions.} Although we write at the end of section \ref{ss:construction} `the conditions on $F(U,V)$ and $G(U,V)$ are quite mild and the slow-fast-slow interfaces $(U_h(\xi),U_h(\xi))$ can be constructed for many models', it is certainly not always possible to construct a traveling front solution that corresponds to a heteroclinic orbit of (\ref{e:DS}) if only condition {\bf (E-I)} holds -- unlike in the scalar case (\ref{e:RDE-scalar}). Conditions {\bf (E-III)} and {\bf (E-IV)} may be violated for generic classes within the full family of systems (\ref{e:RDE}) (see Remark \ref{r:slowpatterns} and \cite{Doe22} for {\bf (E-II)}). It is natural to expect, under the fundamental assumption {\bf (E-I)}, that (\ref{e:RDE}) `must' exhibit some kind of localized patterns that approach $(\bU^\pm,\bV^\pm)$ for $x \to \pm \infty$. We are not aware of any systematic studies of this issue. We note that such a study could start with considering the implications of varying $\vmu$ such that the intersection $\W^{\rm u}((\bV^-,0)) \cap \W^{\rm s}((\bV^+,0))$ vanishes (i.e. the parameters are varied beyond the validity {\bf (E-IV)}). The example of \cite{DIN04} shows that this may lead to (spatially localized) finite-time blow-up behavior. 
\\ \\
Finally, we note that (viscous) fingering has already been identified as `archetype for growth patterns' in fluid mechanics \cite{Cou00}, and has for example also been observed in experiments/observations and model simulations of expanding bacterial colonies \cite{GVC15}, tumor growth \cite{BS18} and wound healing \cite{Marketal10}. At the same time, our understanding of the nonlinear dynamics of fingering -- or non-fingering (Fig. \ref{f:FingCusp-BCDE}) -- interfaces is limited. Due to the relative simplicity of 2-component model (\ref{e:RDE}) and its suitability for mathematical analysis -- by its singularly perturbed nature -- a mathematical study of the evolution of `sharp' slow-fast-slow interfaces in (\ref{e:RDE}) that are unstable under `sideband conditions' (\ref{e:cond-longtrans}) or (\ref{e:cond-longtrans-t-Intro}) may be a promising way to gain novel fundamental insights in the fingering phenomenon.  
\\ \\
{\bf Acknowledgment.} The authors gratefully acknowledge Ehud Meron for his input and feedback, both on the analysis of ecosystem model (\ref{e:RDE-FOTM}) and on the impact of (the magnitude of) $\tau$ on the dynamics of (\ref{e:RDE}). KL, EO, and SR were supported by the NSF REU program through the grant DMS-2204758. PC was supported by the NSF through grants DMS-2204758 and DMS-2105816. AD acknowledges the hospitality of Arnd Scheel and the School of Mathematics during his stay at the University of Minneapolis as Ordway Visiting Professor.

\appendix

\section{A derivation of $\tMa(\htau)$ for the bifurcating traveling interfaces of (\ref{e:RDE-FHN})}
\label{a:tMa}

As in section \ref{sss:InStabBifTrav}, we consider  in this appendix FitzHugh-Nagumo system (\ref{e:RDE-FHN}) with $f^+(v) = 1 + \mu_2 v$, $f^c(v) = \mu_3 v$ and $f^-(v) = - f^+(-v)$: the standing planar interface undergoes a bifurcation into traveling waves as $\ttau$ decreases through $\ttau_\ast$ (\ref{d:ttauast}). We assume $0 < \eps \ll \delta \ll 1$, introduce $\htau$ by setting $\ttau = \ttau_\ast - \delta^2 \htau$ and proceed to derive expression (\ref{e:denominator-FHN-t}) for $\tM_\ast(\htau)$, the denominator of $\la_{2,c}(\htau)$ (\ref{e:la2c-expl-t}) for the travelling interfaces (near the bifurcation). 
\\ \\
To evaluate $\tI_f$ with $\tc$ as in (\ref{e:hvhc-FHN}), we first note that $\tc_\ast \ttau = \hc \ttau_\ast \left(\de + \O(\de^3)\right)$ and that $\tu_\ast(\xi)$ solves,
\beq
\label{e:uxi-FHN}
u_\xi = \frac12 \sqrt{2} \left(1 + \hu_1 \de - u \right) \left(1 - \hu_1 \de + u \right)
\eeq
((\ref{d:fastcubic-gen-2nd}), (\ref{d:fastcubic-gen-1st}), (\ref{e:Kcast-gen})) with
\beq
\label{d:hu1}
\hu_1 = \pm \frac{2 \mu_2 \sqrt{\htau}}{\sqrt{(\mu_3-\mu_2)\sqrt{2(\mu_1-\mu_2)}}} + \O(\de^2).
\eeq
(\ref{e:hvhc-FHN}). Moreover, the location of $\xi = 0$ is by construction such that $(\tua(0),\tu_{\ast, \xi}(0)) = (\de \hu_1, \frac12 \sqrt{2})$, so $\tua(\xi) - \de \hu_1$ is odd around $\xi =0$ and $\tu_{\ast, \xi}(\xi)$ even. Thus, we have
\beq
\label{e:tIf-FHN-exp}
\begin{array}{rcl}
\tI_f & = & \int_\RR (\tu_{\ast, \xi})^2 e^{\tc_\ast \ttau \xi} \, d \xi
\\
& = &
\int_\RR (\tu_{\ast, \xi})^2 \left[1 + \de(\hc + \O(\de^2))(\ttau_\ast + \O(\de^2))\xi + \frac12 \de^2 \hc^2 \ttau_\ast^2 \xi^2 + \O(\de^4) \right] d \xi
\\
& = &
\int_\RR (\tu_{\ast, \xi})^2 d \xi +  \de \hc \ttau_\ast \int_\RR \xi (\tu_{\ast, \xi})^2  d \xi +
\frac12 \de^2 \hc^2 \ttau_\ast^2 \int_\RR \xi^2 (\tu_{\ast, \xi})^2 d \xi + \O(\de^3)
\\
& = &
\frac12 \sqrt{2} \int_{-1 + \de \tu_1}^{1 + \de \tu_1}  \left(1 + \hu_1 \de - u \right) \left(1 - \hu_1 \de + u \right) du +
\frac14 \de^2 \hc^2 \ttau_\ast^2 \int_\RR \frac{\xi^2}{\cosh^4 \frac12 \sqrt{2} \xi} \, d \xi + \O(\de^3)
\\
& = & \frac23 \sqrt{2} +
\frac12 \sqrt{2} \hc^2 \ttau_\ast^2 \K_{2,4} \de^2 + \O(\de^3)
\end{array}
\eeq
(\ref{e:FHNcuast}), with
\beq
\label{d:Kij}
\K_{i,j} = \int_\RR \frac{s^i \, ds}{\cosh^j s}, \; \; i \geq 0, j \geq 1.
\eeq
Since $F_v(u,v) = 2 \mu_2 (u - \mu_3 v)(u - \mu_2 v) + \mu_3 (u + 1 - \mu_2 v)(u - 1 - \mu_2 v)$ for (\ref{e:RDE-FHN}) with the present choices of $f^{\pm, c}(v)$, we have
\beq
\label{e:tFa-FHN-expand}
\begin{array}{rcl}
\tF_\ast & = & 2 \mu_2 \int_\RR (\tua - \frac{\mu_3}{\mu_2} \hu_1 \de)(\tua - \hu_1 \de) \tu_{\ast, \xi} \, e^{\tc_\ast \ttau \xi} d\xi + \mu_3 \int_\RR (\tua + 1 - \hu_1 \de)(\tua - 1 - \hu_1 \de) \tu_{\ast, \xi} \, e^{\tc_\ast \ttau \xi} d\xi
\\
& = &
2 \mu_2 \left( \tF_{h,0} + \de \hc \ttau_\ast \tF_{h,1} + \frac12 \de^2 \hc^2 \ttau_\ast^2 \tF_{h,2}\right) - \sqrt{2} \mu_3 \tI_f + \O(\de^3)
\end{array}
\eeq
(\ref{e:uxi-FHN}) with
\[
\begin{array}{rcl}
\tF_{h,0} & = & \int_{-1 + \de \tu_1}^{1 + \de \tu_1} (u - \frac{\mu_3}{\mu_2} \hu_1 \de)(u - \hu_1 \de) d u = \frac23
\\
\tF_{h,1} & = & \int_\RR \xi \left( (\tua - \de \hu_1) + \left(1-\frac{\mu_3}{\mu_2} \right)\de \hu_1 \right) (\tua - \de \hu_1) \tu_{\ast, \xi} d \xi
\\
& = & \int_\RR \xi (\tua - \de \hu_1)^2 \tu_{\ast, \xi} d \xi  - \frac{\mu_3-\mu_2}{\mu_2} \de \hu_1 \int_\RR \xi \tua \tu_{\ast, \xi} d \xi + \O(\de^2)
\\
& = & - \frac{\mu_3-\mu_2}{2 \mu_2} \de \hu_1 \sqrt{2}  \int_\RR \xi \frac{\tanh \frac12 \sqrt{2} \xi}{\cosh^2 \frac12 \sqrt{2} \xi} d \xi + \O(\de^2)
\\
& = & \frac{\mu_3-\mu_2}{2\mu_2} \hu_1 \sqrt{2} \left(3 \K_{2,4}-2 \K_{2,2}\right) \de + \O(\de^2)
\\
\tF_{h,2} & = & \int_\RR \xi^2 \tua^2 \tu_{\ast, \xi} d\xi = \frac12 \sqrt{2} \int_\RR \xi^2 \frac{\tanh^2 \frac12 \sqrt{2} \xi}{\cosh^2 \frac12 \sqrt{2} \xi} d \xi = 2 \left(\K_{2,2}-\K_{2,4}\right)
\end{array}
\]
(\ref{e:FHNcuast}), (\ref{d:Kij}). Since
\[
\hc^2 \ttau_\ast^2 = 4 \sqrt{2} \frac{(\mu_3-\mu_2)\htau}{\sqrt{\mu_1-\mu_2}} + \O(\de), \;
\hc \ttau_\ast \hu_1  = 4 \sqrt{2} \frac{\mu_2 \htau}{\sqrt{\mu_1-\mu_2}} + \O(\de)
\]
((\ref{e:hvhc-FHN}), (\ref{d:hu1})), it follows by (\ref{e:tIf-FHN-exp}), (\ref{e:tFa-FHN-expand}) that
\beq
\label{e:tIftFa-cneq0}
\tI_f  = \frac23 \sqrt{2}  + \frac{4(\mu_3-\mu_2) \K_{2,4} \htau}{\sqrt{\mu_1-\mu_2}} \de^2, \;
\tF_\ast  = -\frac43 (\mu_3 - \mu_2) -  \frac{4 \sqrt{2}(\mu_3-\mu_2)^2 \K_{2,4} \htau}{\sqrt{\mu_1-\mu_2}} \de^2
\eeq
up to $\O(\de^3)$ corrections. Since $G(u,v) = u - \mu_1 v$ we have by (\ref{e:tIs-FHN}), (\ref{e:hvhc-FHN}) that,
\beq
\label{e:tGatIs-cneq0}
\tG_\ast = 2, \; \; \tI_s = \frac{1}{(\mu_1 - \mu_2)^{3/2}} - \frac{3 \sqrt{2} \htau}{\mu_3 - \mu_2} \de^2 + \O(\de^4).
\eeq
Together, (\ref{e:tIftFa-cneq0}) and (\ref{e:tGatIs-cneq0}) indeed yield (\ref{e:denominator-FHN-t}) -- after some remarkable further simplifications.

\end{document}